\documentclass[11pt]{article}

\author{Damian Sercombe}
\title{A family of uniform lattices acting on a Davis complex with a non-discrete set of covolumes}
\date{29 November 2015}

\pdfoutput=1

\usepackage{amssymb,amsmath,tkz-berge,amsthm}
  \usepackage{times}
  \usepackage[T1]{fontenc}
  \usepackage{textcomp}
 \usepackage{lipsum}
\usepackage[margin=1.2in]{geometry}

\usetikzlibrary{decorations.markings}

\let\OLDthebibliography\thebibliography
\renewcommand\thebibliography[1]{
  \OLDthebibliography{#1}
  \setlength{\parskip}{0pt}
  \setlength{\itemsep}{0pt plus 0.3ex}
}

\numberwithin{equation}{section}

\theoremstyle{definition}
\newtheorem{Definition}[equation]{Definition}
\newtheorem{Remark}[equation]{Remark}
\newtheorem{Example}[equation]{Example}
\newtheorem{Theorem}[equation]{Theorem}
\newtheorem{Proposition}[equation]{Proposition}

\newtheorem{Corollary}[equation]{Corollary}

\newtheorem{Condition}[equation]{Condition}

\DeclareMathOperator{\Int}{Int}
\DeclareMathOperator{\CAT}{CAT}
\DeclareMathOperator{\BS}{BS}

\DeclareMathOperator{\Max}{max}

\DeclareMathOperator{\Aut}{Aut}

\DeclareMathOperator{\St}{St}
\DeclareMathOperator{\Lk}{Lk}
\DeclareMathOperator{\SL}{SL}
\newcommand{\N}{\mathbb{N}}
\newcommand{\R}{\mathbb{R}}

\newcommand{\Q}{\mathbb{Q}}

{\catcode`\|=\active
  \gdef\set#1{\mathinner{\lbrace\,{\mathcode`\|"8000%
                                   \let|\midvert #1}\,\rbrace}}
}
\def\midvert{\egroup\mid\bgroup}

\begin{document}

\tikzset{->-/.style={decoration={
  markings,
  mark=at position .5 with {\arrow{>}}},postaction={decorate}}}
\tikzset{-<-/.style={decoration={
  markings,
  mark=at position .5 with {\arrow{<}}},postaction={decorate}}}

\pagenumbering{arabic}

\maketitle
\begin{abstract}
 \noindent Let $(W,S)$ be a Coxeter system with Davis complex $\Sigma$. The polyhedral automorphism group $G$ of $\Sigma$ is a locally compact group under the compact-open topology. If $G$ is a discrete group (as characterised by Haglund--Paulin \cite{HP}), then the set $\mathcal V_u(G)$ of uniform lattices in $G$ is discrete. Whether the converse is true remains an open problem. Under certain assumptions on $(W,S)$, we show that $\mathcal V_u(G)$ is non-discrete and contains rationals (in lowest form) with denominators divisible by arbitrarily large powers of any prime less than a fixed integer. We explicitly construct our lattices as fundamental groups of complexes of groups with universal cover $\Sigma$. We conclude with a new proof of an already known analogous result for regular right-angled buildings.
\end{abstract}

\section{Introduction}

\noindent Let $G$ be a locally compact topological group equipped with an appropriate normalisation of the Haar measure $\mu$. A discrete subgroup $\Gamma$ is a \textit{lattice} in $G$ if the quotient $\Gamma \backslash G$ admits a finite $G$-invariant measure (which we also call $\mu$), and it is \textit{uniform} if $\Gamma \backslash G$ is compact. 

\vspace{2mm}\noindent A basic question is: what is the set of possible (positive real) covolumes $\mu(\Gamma \backslash G)$ of lattices -- uniform or otherwise -- in $G$? This question was first studied in $1945$ by Siegel \cite{Si} for the case where $G=\SL(2,\R)$. It has since been pursued by various researchers for the case of a semi-simple algebraic group $G$ over a local field. In addition, more recently, this question has been studied by Lubotzky \cite{BL,L}, Bass \cite{BL} and Burger and Mozes \cite{BM} for the automorphism group $G$ of a (product of) tree(s).

\vspace{2mm}\noindent Let $(W,S)$ be a Coxeter system with nerve $L$ and Davis complex $\Sigma$ (refer to Section~\ref{Daviscomplex} for the relevant definitions). The group of polyhedral automorphisms $\Aut(\Sigma)$ is naturally a locally compact topological group under the compact-open topology \cite{T1} (in fact this is true for any connected, locally finite polyhedral complex). In this paper we are interested in determining the set $$\mathcal V_u(\Aut(\Sigma)) := \big\{\mu(\Gamma \backslash \Aut(\Sigma)) \hspace{1mm} \big| \hspace{1mm} \Gamma \textnormal{ is a uniform lattice in } \Aut(\Sigma) \big\}$$ Whilst covolumes of non-uniform lattices may be irrational, Serre's theorem \cite{S1} tells us immediately that $\mathcal V_u(\Aut(\Sigma)) \subseteq \Q^{>0}$.

\vspace{2mm}\noindent In the case where $(W,S)$ is the free product of $n \geq 3$ copies of the cyclic group of order two, with Davis complex $\Sigma_n$ isomorphic to the $(n,2)$-biregular tree, the result is already known from Levich and Rosenberg's (Chapter $9$ of \cite{R}) classification of uniform lattices acting on regular and biregular trees. Precisely, $$\mathcal V_u(\Aut(\Sigma_n)) =\big\{a/b \hspace{1mm} \big| \hspace{1mm} a,b \in \N \textnormal{ are coprime, prime divisors of } b \textnormal{ are } <n \big\}$$ We briefly look at this situation in Section~\ref{result1}. In Proposition~\ref{|||} we give a new proof of the fact that $\mathcal V_u(\Aut(\Sigma_n))$ is a non-discrete set and in Corollary~\ref{||||} provide a new proof of the fact that for any prime $p<n$ and any $\alpha \in \N$, there exists a uniform lattice $\Gamma$ in $\Aut(\Sigma_n)$ with covolume $a/b$ (in lowest terms) such that $b$ is divisible by $p^{\alpha}$. Our proofs only use cyclic groups, dihedral groups and direct products of these, and so are combinatorially simpler than previous proofs.

\vspace{2mm}\noindent In the case where $\Sigma$ is not a tree, not much is known about the set $\mathcal V_u(\Aut(\Sigma))$. For certain types of Davis complex $\Sigma$ it is known that $\mathcal V_u(\Aut(\Sigma))$ is a discrete set. Haglund and Paulin \cite{HP} show that the group $\Aut(\Sigma)$ is non-discrete (as a topological space) iff there exists a non-trivial automorphism of the nerve $L$ as a weighted graph that fixes the star of some vertex in $L$. We call $\Sigma$ \textit{flexible} if it satisfies Haglund and Paulin's condition, and \textit{rigid} otherwise. If $\Sigma$ is rigid then $\mathcal V_u(\Aut(\Sigma))$ must be a discrete set. The converse however remains an open problem. That is, for any flexible Davis complex $\Sigma$, must $\mathcal V_u(\Aut(\Sigma))$ be a non-discrete set?

\vspace{2mm}\noindent In Theorem $2$ of \cite{T1}, Thomas found a restriction on the possible covolumes of uniform lattices in the automorphism group of a connected locally finite polyhedral complex with finitely many isometry classes of links (refer to Section~\ref{heHEe} for the definition of the \textit{link} of a vertex of a polyhedral complex or scwol). This restriction is based only on the prime divisors of the orders of the automorphism groups of the links. We refer the reader to \cite{T1} for the full statement of the theorem, but give an example of its application here. Consider the Coxeter system $$(W_0,S_0):= \big\langle s_1,s_2,s_3 \hspace{1mm} \big| \hspace{1mm} s_1^2=s_2^2=s_3^2=(s_1s_2)^3=(s_1s_3)^3=(s_2s_3)^3=1 \big\rangle$$ with Davis complex $\Sigma_0$ a tesselation of the plane into barycentrically subdivided triangles. Observe that there are only three isometry classes of links of vertices in $\Sigma_0$. The link of any vertex of type $\varnothing$ (respectively $\{s_i\}$, $\{s_i,s_j\}$ for $i \neq j$) is a regular hexagon (respectively square, dodecagon), with automorphism group the dihedral group $\mathcal D_{2 \cdot 6}$ (respectively $\mathcal D_{2 \cdot 4}$, $\mathcal D_{2 \cdot 12}$). The lowest common multiple of the orders of these groups is $24=2^3 \cdot 3$. For each link, the subgroup of its automorphism group that fixes a maximal cell pointwise is always trivial. Applying Theorem $2$ of \cite{T1} tells us that if $\Gamma$ is a uniform lattice in $\Aut(\Sigma_0)$ then its covolume $\mu(\Gamma \backslash \Aut(\Sigma_0))=a/b$ for coprime $a,b \in \N$ is such that $b$ is indivisible by any prime other than $2$ or $3$, and is indivisible by $2^4$ and by $3^2$. Hence $\mathcal V_u(\Aut(\Sigma_0))$ is a discrete set. Alternatively, the discreteness of $\mathcal V_u(\Aut(\Sigma_0))$ follows from the fact that $\Sigma_0$ is rigid.

\vspace{2mm}\noindent In order to state the main result of this paper, we introduce some notation. A generator $s$ of a Coxeter system $(W,S)$ is \textit{free} if there are no relators in the Coxeter presentation of $W$ that include $s$ other than $s^2$. The \textit{label} of an edge in a weighted graph is its weighting.

\vspace{3mm}\noindent \textbf{Main Theorem.} Let $(W,S)$ be a Coxeter system with $j \geq 3$ free generators, and with nerve $L$ and Davis complex $\Sigma$. Assume that each connected component of $L$ is either a cycle graph with an even number of edges and constant edge-labels, or corresponds to a spherical subgroup of $(W,S)$. Then the set $\mathcal V_u(\Aut(\Sigma))$ is non-discrete.

\vspace{3mm}\noindent As far as the author is aware, the Main Theorem was only known previously in the case where $\Sigma$ is a tree. However some related results were already known. Thomas \cite{T2} found a subclass of Davis complexes $\Sigma$ for which there exists an infinite family of uniform lattices in $\Aut(\Sigma)$ with a set of covolumes converging to that of a non-uniform lattice in $\Aut(\Sigma)$. For a (different) subclass of Davis complexes $\Sigma$, White \cite{W1} showed that there exist uniform lattices in $\Aut(\Sigma)$ of arbitrarily small covolume. By considering automorphism(s) of the nerve that swap two free generators and fix all other generators, it is not hard to see that the set of Davis complexes in the Main Theorem is a subset of that found by both Thomas \cite{T2} and by White \cite{W1}.

\vspace{2mm}\noindent Combining Theorem~\ref{coolthm} with Remark~\ref{importantremark} and Proposition~\ref{|||} gives us the Main Theorem. We prove Theorem~\ref{coolthm} by constructing a family of uniform lattices as fundamental groups of a family of finite complexes of finite groups each with universal cover isomorphic to $\Sigma$. We extend these constructions to prove Corollary~\ref{|||||} which, combined with Corollary~\ref{||||}, gives us the following result.

\vspace{3mm}\noindent \textbf{Main Corollary.} Let $(W,S)$ be as in the Main Theorem. For any $\alpha \in \N$ and any prime $p <j$, there exists a uniform lattice $\Gamma$ in $\Aut(\Sigma)$ with covolume $\mu(\Gamma \backslash \Aut(\Sigma))=a/b$ (in lowest terms) such that $b$ is divisible by $p^{\alpha}$.

\vspace{3mm}\noindent We conclude this paper by looking at regular right-angled buildings. Say $S:=\{s_1,...,s_n\}$ and associate an integer $p_i \geq 2$ to generator $s_i$ for every $i$. If $(W,S)$ is a right-angled Coxeter system then there exists a unique regular \textit{right-angled} building $X$ of type $(W,S)$ with $p_i $ chambers in each $s_i$-equivalence class, for all $i\in \{1,...,n\}$ (Theorem $8$ of \cite{T4}). Thomas \cite{T3} showed that if $m:=\Max_{1 \leq i \leq n}\{p_i\} >2$ then $$\mathcal V_u(\Aut(X)) =\big\{a/b \hspace{1mm} \big| \hspace{1mm} a,b \in \N \textnormal{ are coprime, prime divisors of } b \textnormal{ are } <m\big\}$$ In Section~\ref{result3} we look at the situation where $(W,S)$ has $2 \leq j \leq n$ free generators, which we call $s_1,...,s_j$, with $\Max_{1 \leq i \leq j}\{p_i\} >2$. In Theorem~\ref{buildingthm} we give a new proof of the fact that $\mathcal V_u(\Aut(X))$ is a non-discrete set and in Corollary~\ref{bldg} provide a new proof of the fact that for any prime $p<\Max_{1 \leq i \leq j}\{p_i\} $ and any $\alpha \in \N$, there exists a uniform lattice $\Gamma$ in $\Aut(X)$ with covolume $a/b$ (in lowest terms) such that $b$ is divisible by $p^{\alpha}$. Our proofs only use direct products of cyclic groups, and so are again combinatorially simpler than previous proofs.

\vspace{2mm}\noindent I would like to thank my supervisor Dr Anne Thomas for introducing me to these interesting questions, and for the immense amount of time she has spent guiding me through this research.

\section{Definitions and background}

\subsection{Polyhedral complexes}

\noindent We define a polyhedral complex in this section, following Part I of \cite{BH}.

\vspace{2mm}\noindent Let $M^n_{\kappa}$ be the complete, simply connected Riemannian $n$-manifold of constant sectional curvature $\kappa \in \R$. In particular, $M^n_{0}$ (resp. $M^n_{1}$, $M^n_{-1}$) is Euclidean space $\R^n$ (resp. the sphere $\mathbb{S}^n$, hyperbolic space $\mathbb{H}^n$) under the standard metric.

\vspace{2mm}\noindent Let $m \leq n$ be integers and let $\kappa \in \R$. A \textit{$m$-plane} is a subspace of $M^n_{\kappa}$ isometric to $M^m_{\kappa}$. A \textit{compact convex polyhedron} $P$ is the convex hull of a finite set of points in $M^n_{\kappa}$. If $\kappa > 0$ we require all such points to lie within an open hemisphere. The \textit{dimension} of $P$ is the dimension of the smallest $m$-plane containing $P$. The \textit{interior} of $P$ is with respect to this $m$-plane.

\noindent \begin{Definition} Let $\kappa \in \R$. A \textit{$M_{\kappa}$-polyhedral complex} $\mathcal X$ is a $CW$-complex such that for any $k$-cell $\sigma \in \mathcal X$, with attaching map $f:\mathbb S^{k-1} \to \mathcal X$:

\vspace{2mm}\noindent $(i)$ $ \Int (\sigma )$ is isometric to the interior of a compact convex polyhedron in $M^k_{\kappa}$; and

\vspace{1mm}\noindent $(ii)$ the restriction of $f$ to each open, codimension-$1$ face of $\sigma$ is an isometry onto an open cell of $\mathcal X$.

\vspace{2mm}\noindent A $M_{0}$ (resp. $M_{1}$, $M_{-1}$)-polyhedral complex is called \textit{piecewise Euclidean} (resp. \textit{spherical}, \textit{hyperbolic}). A \textit{polyhedral complex} is a $M_{\kappa}$-polyhedral complex, for some $\kappa$.
\end{Definition}

\subsection{Locally compact groups}

\noindent In this section we recall some facts about locally compact groups, referring to Section 2.1. of \cite{T2} and Section 2.2. of \cite{T1}. In particular we describe a locally compact topology on the group of automorphisms of a locally finite polyhedral complex, and characterise uniform lattices in this group.

\vspace{2mm}\noindent Let $G$ be a locally compact group (left)-acting on a set $V$ such that for every $v \in V$ the stabiliser $G_v$ is compact and open. If a subgroup $\Gamma$ of $G$ is discrete, then each stabiliser $\Gamma_v$ is finite \cite{T1}. Now assume that the quotient set $G \backslash V$ is finite.

\noindent \begin{Theorem}[\textbf{Serre}, \cite{S1}]\label{ss} There exists a normalisation of the Haar measure $\mu$ on $G$, depending only on the choice of $G$-set $V$, such that for each discrete subgroup $\Gamma$ of $G$, we have $$\mu (\Gamma \backslash G)=\sum\limits_{v \in \Gamma \backslash V} \frac{1}{|\Gamma_{v'}|}$$ where for each $v \in \Gamma \backslash V$ we choose a single representative $v' \in V$ of the orbit of $v$. Note that $|\Gamma_{v'}|$ is independent of the choice of $v'$.
\end{Theorem}

\noindent Now let $X$ be a connected, locally finite polyhedral complex and let $G=\Aut(X)$ be the group of polyhedral isometries of $X$. Assume that the quotient space $G \backslash X$ is compact. Under the compact-open topology, $G$ is a Hausdorff, first-countable locally compact group with compact open $G$-stabilisers of cells in $X$ (see Section 2.1. of \cite{T2}). We take $V$ to be the vertex set of $X$ and apply Theorem~\ref{ss} to give us a normalisation of the Haar measure $\mu$ on $G$, a formula to compute the covolume of a discrete subgroup of $G$, and the following criterion to test whether a subgroup $\Gamma$ of $G$ is a uniform lattice.

\noindent \begin{Corollary}\label{lllllllllll} A subgroup $\Gamma$ of $G$ is a uniform lattice iff the quotient $\Gamma \backslash X$ is compact and the stabiliser $\Gamma _v$ is finite for every vertex $v$ in $X$.
\begin{proof} Refer to Section 2.2. of \cite{T1}.
\end{proof}
\end{Corollary}

\noindent Alternatively, there are other possible choices of $V$ we could use when applying Theorem~\ref{ss}. For example, we could take $V$ to be the set of chambers of a Davis complex.

\subsection{Coxeter systems and the Davis complex}\label{Daviscomplex}

\noindent In this section we construct and recall some facts about the Davis complex of a Coxeter system. We follow Davis \cite{D}.

\vspace{2mm}\noindent A \textit{Coxeter group} is a group $W$ with a presentation of the form $$W = \big\langle s_1,s_2,...,s_n \hspace{1mm} \big| \hspace{1mm} (s_{i}s_{j})^{m_{ij}}=1, \forall i,j \in \{1,...,n\} \big\rangle $$ where $M=\{m_{ij}\}_{i,j}$ is a symmetric $n \times n$ matrix with entries in the set $\{2,3,...\}\cup\{\infty\}$ unless they are on the major diagonal, in which case they equal $1$. We take $g^{\infty}=1$ to mean that the group element $g$ has infinite order. The set $S:=\{s_1,s_2,...,s_n\}$ generates $W$. The pair $(W,S)$ is called a \textit{Coxeter system} of \textit{type} $M$ and \textit{rank} $n$. A Coxeter system $(W,S)$ is \textit{right-angled} if every non-major diagonal entry of $M$ is either $2$ or $\infty$. A generator $s_k$ is \textit{free} if $ m_{ik}=\infty$ for all $i \neq k$, and $(W,S)$ is \textit{free} if all its generators are free. For any subset $T \subset S$, the subgroup $W_T$ of $W$ generated by $T$ is a \textit{special subgroup}. By convention $W_{\varnothing}$ is the trivial group. The pair $(W_{T},T)$ is itself a Coxeter system (see Theorem 4.1.6. of \cite{D}). Both $T$ and $W_T$ are known as \textit{spherical} if $W_T$ is finite.

\vspace{2mm}\noindent Henceforth let $(W,S)$ be a Coxeter system.

\vspace{2mm}\noindent The \textit{nerve} $L$ of $(W,S)$ is the poset of nontrivial, spherical special subgroups of $W$ under set inclusion. The $1$-skeleton of $L$ may be interpreted as a weighted, undirected graph by identifying each generator $s_i$ (or equivalently the special subgroup $ \langle  s_i \rangle  $ of $W$) with a vertex and by joining two distinct vertices $s_i$ and $s_j$ with an edge labelled $m_{ij}$ iff $m_{ij} < \infty $. The \textit{chamber} $K$ of $(W,S)$ is the topological cone on the barycentric subdivision of $L$, where we identify the point of the cone with the empty set $\varnothing$.  Equivalently, $K$ is the geometric realisation of the abstract simplicial complex $L \cup \varnothing$. For any $s_i \in S$, the $s_i$-\textit{mirror} $K_{s_i}$ of $K$ is the union of all closed simplices in $K$ that contain the vertex $s_i$ but not the vertex $\varnothing$.

\noindent \begin{Definition}[\textbf{Davis, \cite{D}}]\label{Davis complex} The \textit{Davis complex} of $(W,S)$ is the quotient space $\Sigma:= W \times K / \sim $, where $(w,k) \sim (w',k')$ iff $k=k'$ and $w^{-1}w' \in W_{\{s_i \in S \mid  k \in K_{s_i} \}}$.
\end{Definition}

\noindent  For $w \in W$, the translates $(w,K)$ are known as \textit{chambers} of $\Sigma$. Each vertex of $\Sigma$ has a \textit{type} induced by the corresponding element of the poset $L \cup \varnothing$. There is a natural type-preserving action of $W$ on $\Sigma$ given by $(g,(w,k)) \mapsto (gw,k)$. This action freely permutes the set of chambers of $\Sigma$. So we may identify the $W$-orbit classes of $\Sigma$ with the type of their component vertices.

\noindent \begin{Theorem}[\textbf{Moussong, 12.3.3. of \cite{D}}]\label{hjh} The Davis complex $\Sigma$ may be endowed with a piecewise Euclidean metric such that it is a complete $\CAT(0)$ space.
\end{Theorem}

\subsection{Simple complexes of groups}\label{heHEe}

\noindent The theory of complexes of groups is a generalisation of Bass and Serre's theory of graphs of groups. Standard references for Bass--Serre theory are \cite{B,S}. We recall some definitions taken from Part III.$\mathcal C$ of \cite{BH} and Section 2.3 of \cite{T2}. 

\begin{Definition} A \textit{small category without loops} (a.k.a. a \textit{scwol}) $X$ is a disjoint union of a set of vertices $V(X)$ and a set of (directed) edges $E(X)$, with each edge $e \in E(X)$ oriented with an initial vertex $i(e)$ and a terminal vertex $t(e)$, such that:

\vspace{2mm}\noindent $(i)$ $i(e) \neq t(e)$ (i.e. no loops); and

\vspace{2mm}\noindent $(ii)$ if $e'$ is another edge with $t(e)=i(e')$ then there exists a composition edge $ee'$ with $i(ee')=i(e)$ and $t(ee')=t(e')$ (i.e. composition of edges).
\end{Definition}

\begin{Example} Let $(W,S)$ be a Coxeter system with chamber $K$ and Davis complex $\Sigma$. Recall that the vertices of $K$ are labelled by the spherical subsets of $S$. We may interpret $K$ as a scwol by orienting the edges in a manner consistent with set inclusion, that is, pointing towards supersets. This induces a scwol structure on $\Sigma$.
\end{Example}

\noindent We can describe any polyhedral complex combinatorially as a scwol as follows. Fix any $\kappa \in \R$. To any $M_{\kappa}$-polyhedral complex $P$ we take the first barycentric subdivision $\BS(P)$ and define a scwol $X$ with vertex set the (barycentres of) cells in $P$ and edge set the edges in $\BS(P)$ oriented towards cells of smaller dimension. Conversely, for any scwol $X'$, there exists a natural geometric realisation of $X'$ as a piecewise Euclidean (i.e. $M_{0}$) polyhedral complex $P'$ with each cell isometric to a geodesic Euclidean simplex (see 1.3 of Part III.$\mathcal C$ of \cite{BH}). Observe that $P'$ is not necessarily a simplicial complex as the intersection of two distinct cells is not necessarily a single face, rather it is a union of faces. For example, interpret the $2$-torus as a polyhedral complex consisting of a square with opposite sides identified in parallel. Its barycentric subdivision is not a simplicial complex.

\begin{Definition}A \textit{morphism} of scwols is a functor between categories. An \textit{automorphism} of a scwol $X$ is a morphism $\phi:X \to X$ with an inverse, or equivalently, an oriention-preserving polyhedral isometry of the geometric realisation of $X$. An \textit{action} of a group $G$ on a scwol $X$ is a group homomorphism $G \to \Aut(X)$ such that, for any $g \in G$ and $e \in E(X)$, $g \cdot i(e) \neq t(e)$ and if $g$ fixes $i(e)$ then it must also fix $e$.
\end{Definition}

\begin{Definition} A \textit{(simple) complex of groups} $\mathcal A$ over a scwol $X$ associates a \textit{local} (a.k.a. \textit{vertex}) group $\mathcal A _v$ to each vertex $v \in V(X)$ and a monomorphism $\psi_e: \mathcal A _{i(e)} \to \mathcal A _{t(e)}$ to each edge $e \in E(X)$ such that $\psi_{e'} \circ \psi_e = \psi_{ee'}$ for all edges $e$ and $e'$ satisfying $t(e)=i(e')$.
\end{Definition}

\noindent That is, a complex of groups is a commuting diagram of monomorphisms of groups.

\begin{Definition}\label{graphofgroups} A \textit{graph of groups} is a complex of groups over a scwol with all vertices falling into one of two types, a \textit{source} or a \textit{sink}. A sink is not the initial vertex of any edge, whilst a source is the initial (resp. terminal) vertex of precisely two (resp. zero) edges. Confusingly, in the language of Bass--Serre theory, sink vertices are called "vertices" and source vertices are called "edges".
\end{Definition} 

\noindent For any group $G$ and complex of groups $\mathcal A$, we denote by $G \times \mathcal {A}$ the complex of groups obtained by taking the direct product of $G$ with each vertex group in $\mathcal {A}$ along with the corresponding canonical monomorphisms.

\begin{Definition}\label{localdevelopment} Let $\mathcal A$ be a complex of groups over a scwol $X$. Let $v$ be a vertex in $X$. The \textit{local development} or \textit{star} of $\mathcal A$ at $v$ is the polyhedral complex $\St(v)$ given by the affine realisation of the union of the posets 

\vspace{2mm}\noindent $(i)$ $\mathcal P := \Big\{ \big( g \psi_e(\mathcal A_{i(e)}), i(e)  \big)  \hspace{1mm}\Big| \hspace{1mm}  g \in \mathcal A_{v}, e \in E(X) \textnormal{ with } t(e)=v \Big\}$ with partial ordering given by $\big( g \psi_e(\mathcal A_{i(e)}), i(e)  \big) \prec \big( g' \psi_{e'}(\mathcal A_{i(e')}), i(e')  \big)$ if there exists a $f' \in E(X)$ with $i(f')=i(e')$, $t(f')=i(e)$ and $g^{-1}g' \in \psi_e(\mathcal A_{i(e)})$; 

\vspace{2mm}\noindent $(ii)$ $\mathcal Q := \Big\{ \big(\mathcal A_v, t(e) \big)  \hspace{1mm} \Big| \hspace{1mm} e \in E(X) \textnormal{ with } i(e)=v \Big\} $ with partial ordering given by $ \big(\mathcal A_v, t(e) \big) \prec  \big(\mathcal A_v, t(e') \big) $ if there exists a $f \in E(X)$ with $i(f)=t(e')$ and $t(f)=t(e)$; and 

\vspace{2mm}\noindent $(iii)$ the singleton poset $\big\{ (\mathcal A_v, v )\big\}$

\vspace{2mm}\noindent with the additional requirement that $\beta \prec \big(\mathcal A_v, v \big) \prec \alpha$ for any $\alpha \in \mathcal P$ and $\beta \in \mathcal Q$. The \textit{link} $\Lk(v)$ of $\mathcal A$ at $v$ is defined similarly except that we omit the singleton poset $\big\{ (\mathcal A_v, v )\big\}$.

\end{Definition}

\begin{Proposition} Let $\mathcal A$ be a complex of groups over a scwol $X$ and let $\kappa \in \R$. Consider a geometric realisation of $X$ as a $M_{\kappa}$-polyhedral complex with only finitely many isometry classes of cells. Then there exists an induced $M_{\kappa}$-polyhedral complex structure on each local development of $\mathcal A$.
\begin{proof} Refer to 4.14. and 4.16 of Part III.$\mathcal C$ of \cite{BH}.
\end{proof}
\end{Proposition}

\noindent In Chapter~\ref{results} we geometrically realise all of our constructed scwols as piecewise Euclidean polyhedral complexes. We do this by equipping each embedded chamber with Moussong's metric from Theorem~\ref{hjh}. So we may speak of isometry classes of links of our complex of groups constructions. We will only need to compute the link in certain special cases, which are described in the following example.

\begin{Example}\label{helpies} Let $\mathcal A$ be a complex of groups over a scwol $X$. Let $v$ be a vertex in $X$ with local group $\mathcal A_v$. Let $\Phi$ be the poset of all vertices $\{v_1,...,v_k\}$ in $X$ adjacent to $v$ where, for any $i,j \in \{1,...,k\}$, $v_i \prec v_j$ if there exists an edge in $X$ with initial vertex $v_j$ and terminal vertex $v_i$. We use Definition~\ref{localdevelopment} to (combinatorially) compute $\Lk(v)$ of $\mathcal A$ in the following cases.

\vspace{2mm}\noindent \textbf{Case I}: \textit{$\mathcal A_v$ is the trivial group.} Then $\mathcal P$ is greatly simplified as all left cosets of monomorphic images of groups in $\mathcal A_v$ must be trivial. Hence, as a poset, $\Lk(v)$ is isomorphic to $\Phi$.

\vspace{2mm}\noindent \textbf{Case II}: \textit{All edges adjacent to $v$ in $X$ are oriented away from $v$.} Then $\mathcal P$ is empty. Hence, as a poset, $\Lk(v)$ is isomorphic to $\Phi$.

\vspace{2mm}\noindent \textbf{Case III}: \textit{$\Phi$ has only one element $v'$, the edge in $X$ adjacent to both $v$ and $v'$ is oriented towards $v$, and $\mathcal A_{v'}$ has index $n$ in $\mathcal A_v$.} Then $\mathcal Q$ is empty and $\Lk(v) = \mathcal P$ is the poset consisting of all $n$ left cosets of $\mathcal A_{v'}$ in $\mathcal A_{v}$ with no ordering between any of them.

\vspace{2mm}\noindent \textbf{Case IV}: \textit{The geometric realisation of $\Phi$ has two connected components, $\Phi_1$ and $\Phi_2$.} Then all edges adjacent to $v$ in $X$ must be oriented towards $v$, so $\mathcal Q$ is empty. Let $\mathcal A_1$ (resp. $\mathcal A_2$) be the sub-complex of groups of $\mathcal A$ given by removing all vertices in $\Phi_1$ (resp. $\Phi_2$) and their adjacent edges from the underlying scwol $X$. Then $\Lk(v)$ of $\mathcal A$ is the disjoint union of $\Lk(v)$ of $\mathcal A_1$ and $\Lk(v)$ of $\mathcal A_2$.
\end{Example}

\vspace{1mm}\noindent Analogous to the construction in Bass--Serre theory \cite{B,S}, a group $G$ acting on a connected, simply connected scwol $X$ can always be written in the form of a complex of groups over the quotient scwol $G \backslash X$, with vertex groups as the stabilisers of a (choice of) preimage of the vertices under the natural projection $X \to G \backslash X$ and monomorphisms chosen accordingly. It can be shown that the resulting complex of groups is unique up to a notion of isomorphism (see Section $3$ of \cite{T2} for more details). If at least one vertex group is trivial then the action of $G$ on $X$ is faithful \cite{T2}.

\vspace{2mm}\noindent However, the converse does not hold! That is, unlike in Bass--Serre theory, it is not always the case that a complex of groups is associated with such a group action (see Example 12.17. of \cite{BH}). We call a complex of groups $\mathcal A$ over a scwol $Y$ that does correspond with such a group action \textit{strictly developable}, and call $G$ the \textit{fundamental group} $\pi_1(\mathcal A)$, and $X$ the \textit{universal cover} $\widetilde{\mathcal A}$, where $Y$ is isomorphic to $G \backslash X$. Applying Corollary~\ref{lllllllllll} tells us that if $\mathcal A$ is a finite complex of finite groups with at least one trivial vertex group, then $G$ is a uniform lattice in $\Aut(X)$.

\vspace{2mm}\noindent Any graph of groups is strictly developable \cite{B,S}. Our definition of the universal cover of a graph of groups is the barycentric subdivision of Bass and Serre's identically named definition.

\vspace{2mm}\noindent We will prove that all complexes of groups constructed in Chapter~\ref{results} are indeed strictly developable using the following theorem.

\begin{Theorem}[\textbf{Haefliger, 12.28. of \cite{BH}}]\label{Haefliger} A complex of groups is strictly developable if each of its local developments is locally $\CAT(0)$ (a.k.a. \textit{non-positively curved}).
\end{Theorem}

\section{A family of uniform lattices acting on $\mathbf{\Sigma}$ with a non-discrete set of covolumes}\label{results}

\noindent Let $(W,S)$ be a Coxeter system of type $M=\{m_{ij}\}_{i,j}$ with nerve $L$. Denote $S:=\{s_1,...,s_n\}$. Recall that a generator $s_k$ is \textit{free} if $ m_{ik}=\infty$ for all $i \neq k$. Let $S'$ denote the subset of free generators of $S$. If $S'$ is non-empty, without loss of generality write $S':=\{s_1,...,s_j\}$ for some $1 \leq j \leq n$. Take $j=0$ to mean that $S'$ is empty.

\vspace{2mm}\noindent In this chapter we prove the following theorem.

\begin{Theorem}\label{centraltheorem} The set of covolumes of uniform lattices in the automorphism group of a polyhedral complex $\Sigma$ is non-discrete, where $\Sigma$ can be either of the following.

\vspace{2mm}\noindent $(i)$ The Davis complex of $(W,S)$, where $j\geq 3$, such that the connected components of $L$ either: \begin{itemize}
\item correspond to spherical subgroups of $(W,S)$; or 
\item are cycle graphs with an even number of edges and constant edge-labels. Distinct cycles are not required to have the same labels or number of edges.
\end{itemize}

\vspace{0mm}\noindent $(ii)$ The regular right-angled building of type $(W,S)$, with $p_i \geq 2$ chambers in each $s_i$-equivalence class for all $1 \leq i \leq n$, where \begin{itemize} 
\item $j \geq 2$;
\item $(W,S)$ is right-angled; and
\item there exists at least one $ i \in \{1,...,j\}$ with $p_i>2$.
\end{itemize}
\end{Theorem}

\noindent We will also prove the following corollary, which is obtained from the constructions we use to establish Theorem~\ref{centraltheorem}.

\begin{Corollary}\label{centralcorollary} Let $\Sigma$ be as in Theorem~\ref{centraltheorem}. In case $(i)$ (resp. case $(ii)$) let $p <j$ (resp. $ p < \Max_{1 \leq i \leq j}\{p_i\} $) be prime. Then for arbitrarily large $\alpha \in \N$ there exists a uniform lattice $\Gamma$ in $\Aut(\Sigma)$ with covolume $\mu(\Gamma \backslash \Aut(\Sigma))=a/b$ (in lowest terms) such that $b$ is divisible by $p^{\alpha}$.
\end{Corollary}

\vspace{1mm}\noindent A sketch of the proof of Theorem~\ref{centraltheorem} is as follows. We explicitly construct a family of finite complexes of finite groups with at least one trivial local group. Assuming strict developability, we apply Corollary~\ref{lllllllllll} of Serre's theorem \cite{S1} to confirm that the fundamental groups of our constructions are indeed uniform lattices in the (polyhedral) automorphism groups of their respective universal covers, compute their respective covolumes as the sum of reciprocals of the orders of the vertex groups and observe that this set of covolumes is non-discrete. We check that the set of isometry classes of links of our complex of groups constructions is the same as that for $\Sigma$ (interpreted as a complex of groups with all vertex groups trivial over itself as a scwol). In particular, they are all $\CAT(0)$ under the restriction of Moussong's metric, so we can apply Theorem~\ref{Haefliger} to ensure that all of our constructions are indeed strictly developable. Finally, we use certain case-specific uniqueness arguments to show that the universal cover of each of our constructions is isomorphic to $\Sigma$ as a polyhedral complex.

\vspace{2mm}\noindent We then prove Corollary~\ref{centralcorollary} by generalising our constructions so that the associated set of covolumes converges at a rate of $1/p$ for any such prime $p$.

\vspace{2mm}\noindent A few remarks on notation. It should be assumed that every monomorphism not explicitly defined in our complex of groups constructions in Figures~\ref{gg} to ~\ref{gggg} is the natural inclusion. Vertices of the same colour in these diagrams have isometric local developments unless otherwise specified. Let $\boldsymbol{1}$ denote the trivial group. Let $\mathcal D_{2m}$ denote the dihedral group of order $2m$ and let $\mathcal C^{\alpha}_{\beta}$ denote the direct product of $\alpha$ copies of the cyclic group of order $\beta$. Assume the convention that $\mathcal C_{\beta}^0$ is the trivial group. A rational number $a/b$ is said to be in lowest terms if $a$ and $b$ are coprime integers. Let $\mu$ be the normalisation of the Haar measure on $\Aut(\Sigma)$ as in Theorem~\ref{ss}. Recall that a convergent sequence $(x_1,x_2,...) \to x$ of real numbers has a \textit{rate of convergence} of $ \lim_{l \to \infty} \frac{|x_{l+1}-x|}{|x_l-x|} $.

\subsection{$\mathbf{\Sigma}$ is the Davis complex of a free Coxeter system}\label{result1}

\noindent In this section we prove Theorem~\ref{centraltheorem} and Corollary~\ref{centralcorollary} in the particular case where $\Sigma$ is the Davis complex of a free Coxeter system. Recall from Definition~\ref{graphofgroups} that a graph of groups is a complex of groups consisting of sink and source vertices. In this section we use the term "vertex" (resp. "edge") to refer to a sink (resp. source) vertex. This is consistent with the language of Bass--Serre theory.

\vspace{2mm}\noindent Let $(W,\{s_1,...,s_n\})$ be a free Coxeter system with Davis complex $\Sigma$. The nerve $L$ consists of $n$ discrete points and the chamber $K$ is isomorphic to the complete $(n,1)$-bipartite graph. Since each $s_i$-mirror of $K$ consists of a single vertex therefore the Davis complex $\Sigma$ is isomorphic to the unique $(n,2)$-biregular tree.

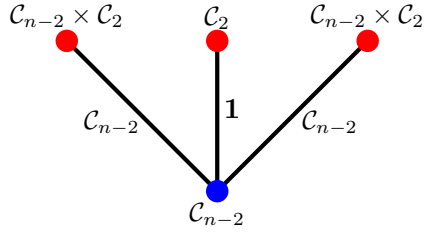
\begin{figure} [h]\
\centering
\usetikzlibrary{arrows,positioning}
\begin{tikzpicture} [%
    nd/.style = {circle,fill=red,text=black,inner sep=3pt},
    tn/.style = {circle,fill=blue,text=black,inner sep=3pt}]

 \draw [black, ultra thick,fill=none,-] (-2,0) node[above, black] {\small $\mathcal C_{n-2} \times \mathcal C_2$} -- (0,-2) node [below, black] {\small $\mathcal C_{n-2}$};
 \draw [black, ultra thick,fill=none,-] (0,0) node[above, black] {\small $\mathcal C_2$} -- (0,-2) node {};
 \draw [black, ultra thick,fill=none,-] (2,0) node[above, black] {\small $\mathcal C_{n-2} \times \mathcal C_2$} -- (0,-2) node {};

\node[tn] at (0,-2) {};
\node[nd] at (-2,0) {};
\node[nd] at (0,0) {};
\node[nd] at (2,0) {};

\node[] at (1.5,-1.1) [black]{\small $\mathcal C_{n-2}$};
\node[] at (0.2,-0.9) [black]{$\boldsymbol{1}$};
\node[] at (-1.4,-1.1) [black]{\small $\mathcal C_{n-2}$};

\end{tikzpicture}
\caption{A finite graph of finite groups $(G,A)$ } \label{gg'}
\end{figure}

\begin{Proposition}\label{|||} The set of covolumes of uniform lattices in $\Aut(\Sigma)$ is non-discrete if $n \geq 3$.

\begin{figure} [h]\
\centering
\usetikzlibrary{arrows,positioning}
\begin{tikzpicture} [%
    nd/.style = {circle,fill=red,text=black,inner sep=3pt},
    tn/.style = {circle,fill=blue,text=black,inner sep=3pt}]

 \draw [black, ultra thick,fill=none,-] (-4,0) node[above, black] {\small $\mathcal C_2$} -- (-3,-2) node [below, black] {\small $\mathcal C_{n-1}$};
\draw [black, ultra thick,fill=none,-] (-3,-2) node{} -- (-2,0) node [above, black] {\small $\mathcal C_{n-1}$};
\draw [black, ultra thick,fill=none,-] (-2,0) node{} -- (-1,-2) node [below, black] {\small $\mathcal C_{n-1}^2$};
\draw [black, ultra thick,fill=none,-] (-1,-2) node{} -- (0,0) node [above, black] {\small $\mathcal C_{n-1}^2$};

\draw [black, ultra thick,fill=none,-] (3,-2) node[below, black] {\small $\mathcal C_{n-1}^{k-1}$} -- (4,0) node [above, black] {\small $\mathcal C_{n-1}^{k-1}$};
\draw [black, ultra thick,fill=none,-] (4,0) node{} -- (5,-2) node [below, black] {\small $\mathcal C_{n-1}^{k}$};
\draw [black, ultra thick,fill=none,-] (5,-2) node{} -- (6,0) node [above, black] {\small $ \mathcal C_{n-1}^{k}$};
\draw [black, ultra thick,fill=none,-] (6,0) node{} -- (7,-2) node [below, black] {\small $\mathcal C_{n} \times \mathcal C_{n-1}^{k}$};

\draw [black, ultra thick,fill=none,->] (0,0) node{} -- (0.5,-1) node {};
\draw [black, ultra thick,fill=none,-] (2.5,-1) node{} -- (3,-2) node {};

\node[nd] at (-4,0) {};
\node[tn] at (-3,-2) {};
\node[nd] at (-2,0) {};
\node[tn] at (-1,-2) {};
\node[nd] at (0,0) {};

\node[circle,fill=black,text=black,inner sep=1.1pt] at (1.2,-1) {};
\node[circle,fill=black,text=black,inner sep=1.1pt] at (1.5,-1) {};
\node[circle,fill=black,text=black,inner sep=1.1pt] at (1.8,-1) {};

\node[tn] at (3,-2) {};
\node[nd] at (4,0) {};
\node[tn] at (5,-2) {};
\node[nd] at (6,0) {};
\node[tn] at (7,-2) {};

\node[] at (-3.75,-1.1) [black]{$\boldsymbol{1}$};
\node[] at (-2.8,-0.7) [black]{\small $\mathcal C_{n-1}$};
\node[] at (-1.75,-1.3) [black]{\small $\mathcal C_{n-1}$};
\node[] at (-0.8,-0.6) [black]{\small $\mathcal C_{n-1}^2$};

\node[] at (3.2,-0.6) [black]{\small $\mathcal C_{n-1}^{k-1}$};
\node[] at (4.2,-1.4) [black]{\small $\mathcal C_{n-1}^{k-1}$};
\node[] at (5.2,-0.6) [black]{\small $\mathcal C_{n-1}^{k}$};
\node[] at (6.25,-1.35) [black]{\small $\mathcal C_{n-1}^{k}$};
\end{tikzpicture}
\caption{A family $\{(G,A)_k\}_{k \in \N}$ of finite graphs of finite groups} \label{gg}
\end{figure}
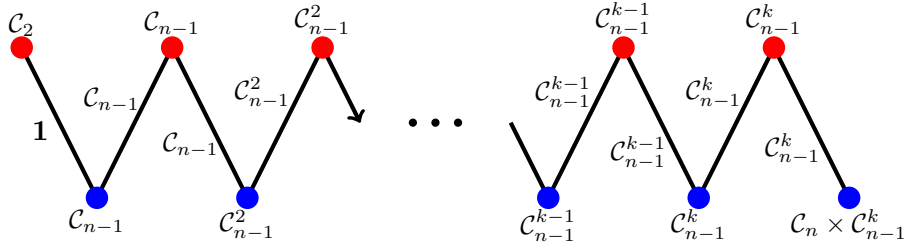

\begin{proof} Consider the family $\{(G,A)_k\}_{k \in \N}$ of graphs of groups in Figure~\ref{gg}. Fix a $k \in \N$. Since $(G,A)_k$ has a trivial edge group, the fundamental group $\pi_1 (G,A)_k$ acts faithfully on the universal cover $\widetilde{(G,A)_k}$. 

\vspace{2mm}\noindent Cases III and IV of Example~\ref{helpies} tell us that the link of each blue (resp. red) vertex of $A_k$ is $n$ (resp. $2$) disconnected edges. Moreover, the link of each edge is the disjoint union of one blue vertex and one red vertex. Since the universal cover of any graph of groups is a tree \cite{B,S}, therefore $\widetilde{(G,A)_k}$ is isomorphic to the (barycentric subdivision $\BS$ of the) unique $(n,2)$-biregular tree $\Sigma$. Observe that $\Aut(\Sigma)=\Aut(\BS(\Sigma))$ and that the orientation on the edges of $\Sigma$ naturally prohibits any edge-inversions if $n\geq 3$.

\vspace{2mm}\noindent Having satisfied all the conditions required to apply Corollary~\ref{lllllllllll} to $(G,A)_k$, we deduce that $\pi_1 (G,A)_k$ is a uniform lattice in $\Aut(\Sigma)$ with covolume \begin{align*} 
\mu \big(\pi_1 (G,A)_k \backslash & \Aut(\Sigma) \big) = \sum\limits_{v \in V(A_k)} \frac{1}{|G_{v}|}   \\
&=  \frac{1}{2}+ 2\bigg(\sum\limits_{i=1}^{k} \frac{1}{(n-1)^i}\bigg)+\frac{1}{n(n-1)^{k}}\\
&\longrightarrow \frac{1}{2}+ \frac{2}{n-2} \hspace{2mm} \textnormal{ as } \hspace{2mm} k \longrightarrow \infty \end{align*} with a rate of convergence of $\frac{1}{n-1}$.

\vspace{2mm}\noindent Observe that the limit point found above may be realised as the covolume of the uniform lattice $\pi_1 (G,A)$ in $\Aut(\Sigma)$; see Figure~\ref{gg'}.
\end{proof}
\end{Proposition}

\noindent In fact the constructions in the proof of Proposition~\ref{|||} give us a new proof of the following result of Rosenberg's (see Theorem 9.2.1 of \cite{R}).

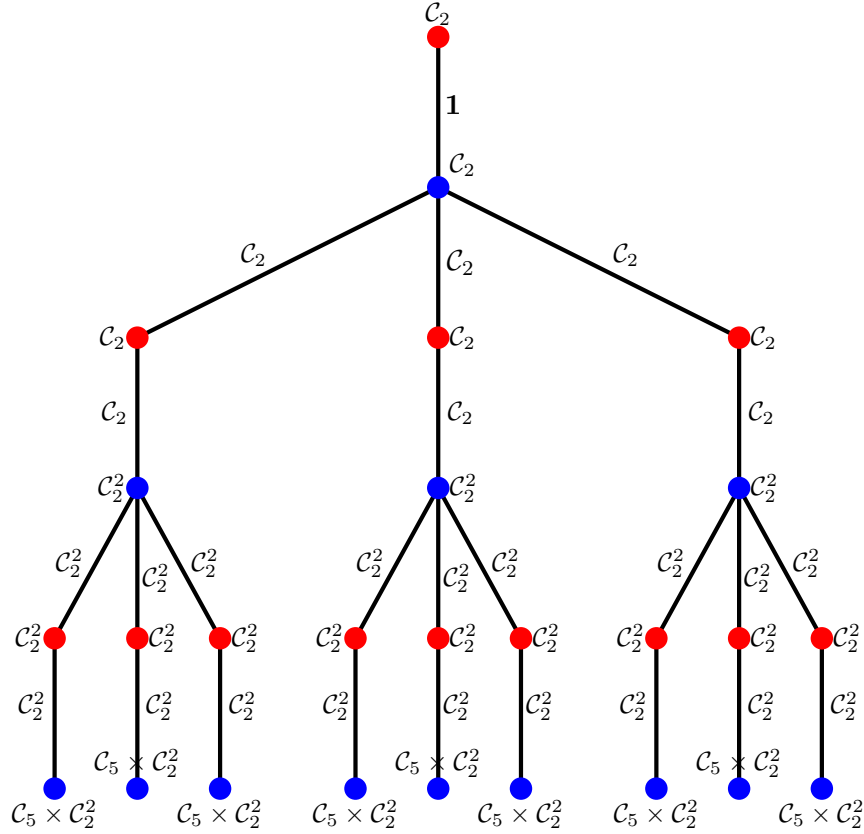
\begin{figure} [h]\
\centering
\usetikzlibrary{arrows,positioning}
\begin{tikzpicture} [%
    red/.style = {circle,fill=red,text=black,inner sep=3pt},
    blue/.style = {circle,fill=blue,text=black,inner sep=3pt}]

\draw [black, ultra thick,fill=none,-] (0,4) node[above, black] {\small $\mathcal C_2$} -- (0,2) node {};
 \draw [black, ultra thick,fill=none,-] (0,2) node[above right, black] {\small $ \mathcal C_2$} -- (0,0) node [right, black] {\small $\mathcal C_2$};
 \draw [black, ultra thick,fill=none,-] (4,0) node[right, black] {\small $\mathcal C_2$} -- (0,2) node {};
 \draw [black, ultra thick,fill=none,-] (-4,0) node[left, black] {\small $\mathcal C_2$} -- (0,2) node {};
 \draw [black, ultra thick,fill=none,-] (4,-2) node[right, black] {\small $\mathcal C_2^2$} -- (4,0) node {};
 \draw [black, ultra thick,fill=none,-] (-4,-2) node[left, black] {\small $\mathcal C_2^2$} -- (-4,0) node {};
 \draw [black, ultra thick,fill=none,-] (0,-2) node[right, black] {\small $\mathcal C_2^2$} -- (0,0) node {};
 \draw [black, ultra thick,fill=none,-] (-2.9,-4) node[right, black] {\small $\mathcal C_2^2$} -- (-4,-2) node {};
 \draw [black, ultra thick,fill=none,-] (-4,-4) node[right, black] {\small $\mathcal C_2^2$} -- (-4,-2) node {};
 \draw [black, ultra thick,fill=none,-] (-5.1,-4) node[left, black] {\small $\mathcal C_2^2$} -- (-4,-2) node {};
 \draw [black, ultra thick,fill=none,-] (1.1,-4) node[right, black] {\small $\mathcal C_2^2$} -- (0,-2) node {};
 \draw [black, ultra thick,fill=none,-] (0,-4) node[right, black] {\small $\mathcal C_2^2$} -- (0,-2) node {};
 \draw [black, ultra thick,fill=none,-] (-1.1,-4) node[left, black] {\small $\mathcal C_2^2$} -- (0,-2) node {};
 \draw [black, ultra thick,fill=none,-] (2.9,-4) node[left, black] {\small $\mathcal C_2^2$} -- (4,-2) node {};
 \draw [black, ultra thick,fill=none,-] (4,-4) node[right, black] {\small $\mathcal C_2^2$} -- (4,-2) node {};
 \draw [black, ultra thick,fill=none,-] (5.1,-4) node[right, black] {\small $\mathcal C_2^2$} -- (4,-2) node {};
 \draw [black, ultra thick,fill=none,-] (-2.9,-6) node[below, black] {\small $\mathcal C_5\times \mathcal C_2^2$} -- (-2.9,-4) node {};
 \draw [black, ultra thick,fill=none,-] (-4,-6) node[above, black] {\small $\mathcal C_5\times \mathcal C_2^2$} -- (-4,-4) node {};
 \draw [black, ultra thick,fill=none,-] (-5.1,-6) node[below, black] {\small $\mathcal C_5\times \mathcal C_2^2$} -- (-5.1,-4) node {};
 \draw [black, ultra thick,fill=none,-] (1.1,-6) node[below, black] {\small $\mathcal C_5\times \mathcal C_2^2$} -- (1.1,-4) node {};
 \draw [black, ultra thick,fill=none,-] (0,-6) node[above, black] {\small $\mathcal C_5\times \mathcal C_2^2$} -- (0,-4) node {};
 \draw [black, ultra thick,fill=none,-] (-1.1,-6) node[below, black] {\small $\mathcal C_5\times \mathcal C_2^2$} -- (-1.1,-4) node {};
 \draw [black, ultra thick,fill=none,-] (2.9,-6) node[below, black] {\small $\mathcal C_5\times \mathcal C_2^2$} -- (2.9,-4) node {};
 \draw [black, ultra thick,fill=none,-] (4,-6) node[above, black] {\small $\mathcal C_5\times \mathcal C_2^2$} -- (4,-4) node {};
 \draw [black, ultra thick,fill=none,-] (5.1,-6) node[below, black] {\small $\mathcal C_5\times \mathcal C_2^2$} -- (5.1,-4) node {};

\node[red] at (0,4) {};
\node[blue] at (0,2) {};
\node[red] at (-4,0) {};
\node[red] at (0,0) {};
\node[red] at (4,0) {};
\node[blue] at (-4,-2) {};
\node[blue] at (0,-2) {};
\node[blue] at (4,-2) {};
\node[red] at (2.9,-4) {};
\node[red] at (4,-4) {};
\node[red] at (5.1,-4) {};
\node[red] at (1.1,-4) {};
\node[red] at (0,-4) {};
\node[red] at (-1.1,-4) {};
\node[red] at (-2.9,-4) {};
\node[red] at (-4,-4) {};
\node[red] at (-5.1,-4) {};
\node[blue] at (2.9,-6) {};
\node[blue] at (4,-6) {};
\node[blue] at (5.1,-6) {};
\node[blue] at (1.1,-6) {};
\node[blue] at (0,-6) {};
\node[blue] at (-1.1,-6) {};
\node[blue] at (-2.9,-6) {};
\node[blue] at (-4,-6) {};
\node[blue] at (-5.1,-6) {};

\node[] at (0.2,3.1) [black]{$\boldsymbol{1}$};
\node[] at (2.5,1.1) [black]{\small $\mathcal C_{2}$};
\node[] at (0.3,1) [black]{$\mathcal C_{2}$};
\node[] at (-2.45,1.1) [black]{\small $\mathcal C_{2}$};
\node[] at (4.3,-1) [black]{\small $\mathcal C_{2}$};
\node[] at (0.3,-1) [black]{\small $\mathcal C_{2}$};
\node[] at (-4.3,-1) [black]{\small $\mathcal C_{2}$};
\node[] at (3.1,-3) [black]{\small $\mathcal C_{2}^2$};
\node[] at (4.9,-3) [black]{\small $\mathcal C_{2}^2$};
\node[] at (4.25,-3.2) [black]{\small $\mathcal C_{2}^2$};
\node[] at (-0.9,-3) [black]{\small $\mathcal C_{2}^2$};
\node[] at (0.9,-3) [black]{\small $\mathcal C_{2}^2$};
\node[] at (0.25,-3.2) [black]{\small $\mathcal C_{2}^2$};
\node[] at (-3.1,-3) [black]{\small $\mathcal C_{2}^2$};
\node[] at (-4.9,-3) [black]{\small $\mathcal C_{2}^2$};
\node[] at (-3.75,-3.2) [black]{\small $\mathcal C_{2}^2$};
\node[] at (2.6,-4.9) [black]{\small $\mathcal C_{2}^2$};
\node[] at (4.3,-4.9) [black]{\small $\mathcal C_{2}^2$};
\node[] at (5.4,-4.9) [black]{\small $\mathcal C_{2}^2$};
\node[] at (1.35,-4.9) [black]{\small $\mathcal C_{2}^2$};
\node[] at (0.25,-4.9) [black]{\small $\mathcal C_{2}^2$};
\node[] at (-1.35,-4.9) [black]{\small $\mathcal C_{2}^2$};
\node[] at (-2.6,-4.9) [black]{\small $\mathcal C_{2}^2$};
\node[] at (-3.7,-4.9) [black]{\small $\mathcal C_{2}^2$};
\node[] at (-5.4,-4.9) [black]{\small $\mathcal C_{2}^2$};

\end{tikzpicture}
\caption{The graph of groups $(G^2_2,A^2_2)$ when $n=5$} \label{poo}
\end{figure}

\noindent\begin{Corollary}\label{||||} For any $n \geq 3$ let $\Sigma$ be the $(n,2)$-biregular tree. Let $p < n$ be prime. For any (arbitrarily large) $\alpha \in \N$ there exists a uniform lattice $\Gamma$ in $\Aut(\Sigma)$ with covolume $\mu \big(\Gamma \backslash \Aut(\Sigma) \big) = a/b$ (in lowest terms) such that $b$ is divisible by $p^{\alpha}$.

\begin{proof} We construct a family of graphs of groups by generalising the construction in Figure~\ref{gg}, which will become the case $p=n-1$. We give an example of the case where $p=2$, $n=5$ and $k=2$ in Figure~\ref{poo}.

\vspace{2mm}\noindent Take any integer $k$. We will construct a tree $A^p_k$ with blue vertices and red vertices, and a complex of groups $(G^p_k,A^p_k)$. No vertex in $A^p_k$ is adjacent to another of the same colour. We partition $(G^p_k,A^p_k)$ into descending "levels" from $i=0$ to $i=k+1$.  Level $0$ consists of a single red vertex with local group $\mathcal C_2$ and a single edge with trivial local group connecting the lone red vertex in level $0$ to the lone blue vertex in level $1$. Level $k+1$ consists of $(n-p)^{k}$ blue vertices each with local group $\mathcal C_n \times \mathcal C_p^k$, no red vertices and no edges.

\vspace{2mm}\noindent Every remaining level $i$ for $i=1,...,k$ consists of $(n-p)^{i-1}$ blue vertices, $(n-p)^{i}$ red vertices, and $2(n-p)^{i}$ edges connecting each blue vertex on level $i$ with $n-p$ red vertices on level $i$ and each red vertex on level $i$ with one blue vertex on level $i+1$. Every vertex and edge group on level $i$, for $i=1,...,k$, is $\mathcal C_p^i$.

\vspace{2mm}\noindent The fundamental group of each graph of groups in $\{(G^p_k,A^p_k)\}_{k \in \N}$ is a uniform lattice in $\Sigma$ by the same argument as in the proof of Proposition~\ref{|||}. A simple calculation will show that the corresponding sequence of covolumes converges at a rate of $\frac{1}{p}$ as $k \to \infty$. The sum of reciprocals of orders of the vertex groups (in lowest terms) will therefore have a denominator that is divisible by ever-increasing powers of $p$ as $k$ gets larger. 
\end{proof}
\end{Corollary}

\subsection{$\mathbf{\Sigma}$ is the Davis complex of a non-free Coxeter system}\label{result2}

\noindent In this section we prove case $(i)$ of Theorem~\ref{centraltheorem} and Corollary~\ref{centralcorollary} in full generality.

\vspace{2mm}\noindent Let $(W,S)$ be any Coxeter system with nerve $L$. We can uniquely decompose $W$ into a free product $W= \ast_{i=1,...,\lambda} W_i$ of Coxeter subgroups with respective generating sets disjoint subsets of $S$. Each $W_i$ is indecomposable if $\lambda$ is a maximal choice of integer. Observe that this decomposition of $W$ corresponds with a decomposition of $L$ as a disjoint union of its connected components $L_i$, where $L_i$ is the subgraph of $L$ corresponding to $W_i$.

\begin{Theorem}\label{coolthm} Let $(W,S)$ be a Coxeter system with Davis complex $\Sigma$ and nerve $L$ containing at least four connected components.  Completely decompose $W= \ast_{i=1,...,l+4} W_i$ as a free product of Coxeter subgroups for some non-negative integer $l$. Then the set of covolumes of uniform lattices in $\Aut(\Sigma)$ is non-discrete if the following conditions hold:
\begin{itemize}
\item (at least) one connected component $L_{l+1}$ of $L$ is the $2m$-cycle graph for some $m\geq 2$ with edges labelled by some constant integer $x \geq 2$;
\item (at least) three connected components $L_{l+2}$, $L_{l+3}$ and $L_{l+4}$ of $L$ are lone vertices (that is, $S$ contains at least three free generators); and

\item if $l >0$ each $W_i$ for $1 \leq i \leq l$ is either finite or corresponds to an $L_i$ that is the $2m_i$-cycle graph for some $m_i\geq 2$ with edges labelled by some constant integer $x_i \geq 2$.
\end{itemize}

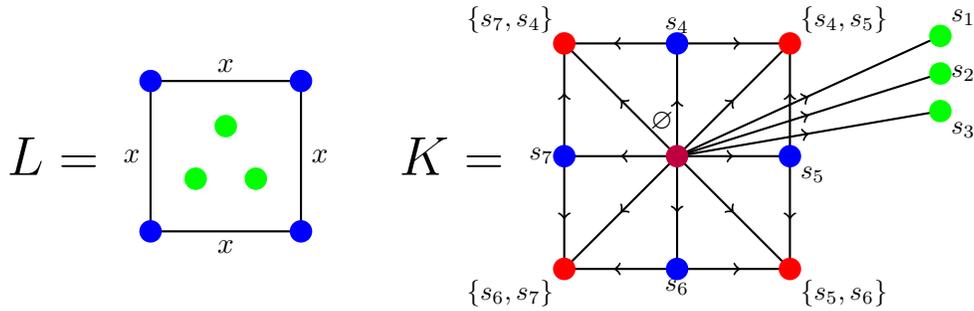
\begin{figure}[h] \
\centering
\usetikzlibrary{arrows,positioning}
\begin{tikzpicture} [%
    purple/.style = {circle,fill=purple,text=black,inner sep=3pt},
    blue/.style = {circle,fill=blue,text=black,inner sep=3pt},
    red/.style = {circle,fill=red,text=black,inner sep=3pt},
   green/.style = {circle,fill=green,text=black,inner sep=3pt}]

\draw [black, thick,fill=none,->-] (0,0) node {} -- (1.5,1.5) node [above right, black]{\small $\{s_4,s_5\}$};
\node[] at (-0.2,0.2) [above, black] {$\varnothing$};
\draw [black, thick,fill=none,->-] (0,0) node {} -- (-1.5,1.5) node [above left, black]{\small $\{s_7,s_4\}$};
\draw [black, thick,fill=none,->-] (0,0) node {} -- (-1.5,-1.5) node [below left, black] {\small $\{s_6,s_7\}$};
\draw [black, thick,fill=none,->-] (0,0) node {} -- (1.5,-1.5) node [below right, black] {\small $\{s_5,s_6\}$};
\draw [black, thick,fill=none,->-] (0,0) node {} -- (0,1.5) node [above, black] {\small $s_4$};
\draw [black, thick,fill=none,->-] (0,0) node {} -- (1.5,0) node [below right, black] {\small $s_5$};
\draw [black, thick,fill=none,->-] (0,0) node {} -- (-1.5,0) node [left, black] {\small $s_7$};
\draw [black, thick,fill=none,->-] (0,0) node {} -- (0,-1.5) node [below, black] {\small $s_6$};
\draw [black, thick,fill=none,->-] (0,0) node {} -- (3.5,1.6) node [above right, black] {\small $s_1$};
\draw [black, thick,fill=none,->-] (0,0) node {} -- (3.5,1.1) node [right, black] {\small $s_2$};
\draw [black, thick,fill=none,->-] (0,0) node {} -- (3.5,0.6) node [below right, black] {\small $s_3$};
\draw [black, thick,fill=none,-<-] (1.5,1.5) node {} -- (0,1.5) node {};
\draw [black, thick,fill=none,-<-] (1.5,1.5) node {} -- (1.5,0) node {};
\draw [black, thick,fill=none,-<-] (-1.5,1.5) node {} -- (-1.5,0) node {};
\draw [black, thick,fill=none,-<-] (-1.5,1.5) node {} -- (0,1.5) node {};
\draw [black, thick,fill=none,-<-] (1.5,-1.5) node {} -- (1.5,0) node {};
\draw [black, thick,fill=none,-<-] (1.5,-1.5) node {} -- (0,-1.5) node {};
\draw [black, thick,fill=none,-<-] (-1.5,-1.5) node {} -- (0,-1.5) node {};
\draw [black, thick,fill=none,-<-] (-1.5,-1.5) node {} -- (-1.5,0) node {};
\node[purple] at (0,0) {};
\node[red] at (1.5,1.5) {};
\node[red] at (1.5,-1.5) {};
\node[red] at (-1.5,1.5) {};
\node[red] at (-1.5,-1.5) {};
\node[blue] at (0,1.5) {};
\node[blue] at (0,-1.5) {};
\node[blue] at (1.5,0) {};
\node[blue] at (-1.5,0) {};
\node[green] at (3.5,1.6) {};
\node[green] at (3.5,1.1) {};
\node[green] at (3.5,0.6) {};
\node[] at (-3,0) [black]{\huge $K=$};

\draw [red, thick,fill=none,-] (-5,1) node {} -- (-5,-1) node{};
\draw [red, thick,fill=none,-] (-5,-1) node {} -- (-7,-1) node {};
\draw [red, thick,fill=none,-] (-7,-1) node {} -- (-7,1) node {};
\draw [red, thick,fill=none,-] (-7,1) node {} -- (-5,1) node {};
\node[blue] at (-5,1) {};
\node[blue] at (-5,-1) {};
\node[blue] at (-7,-1) {};
\node[blue] at (-7,1) {};
\node[green] at (-6.4,-0.3) {};
\node[green] at (-5.6,-0.3) {};
\node[green] at (-6,0.4) {};
\node[] at (-6,1.2) [black]{$x$};
\node[] at (-6,-1.2) [black]{$x$};
\node[] at (-4.75,0) [black]{$x$};
\node[] at (-7.25,0) [black]{$x$};
\node[] at (-8.3,0) [black]{\huge $L=$};
\end{tikzpicture}
\caption{The nerve $L$ and chamber $K$ when $l=0$, $m=2$} \label{ppppp}
\end{figure}

\begin{proof} The structure of the proof is as follows. We consider two cases, $l=0$ and $l>0$. When $l=0$ we construct a family $\{\mathcal X_k^{(0)}\}_{k \in \N} \cup \{ \mathcal X^{(0)}\}$ of complexes of groups. For any $\mathcal F$ in this family, we show that the set of isometry classes of links of $\mathcal F$ is the same as that of $\Sigma$ (Proposition~\ref{links}), that $\mathcal F$ is strictly developable (Corollary~\ref{pfpf}) and that the universal cover of $\mathcal F$ is $\Sigma$ (Proposition~\ref{lplpl}). We then apply Corollary~\ref{lllllllllll} to show that $\pi_1(\mathcal F)$ is a uniform lattice in $\Aut(\Sigma)$ and compute the corresponding covolume. When $l>0$ we inductively construct a family $\{\mathcal X_k\}_{k \in \N} \cup \{ \mathcal X\}$ of complexes of groups and modify the argument in the $l=0$ case to establish the result.

\vspace{4mm}\noindent \textbf{\underline{Case. {$ \boldsymbol1$}}:} $l=0$

\vspace{2mm}\noindent Take $\{s_1,s_2,s_3\}$ to be the subset of free generators of $S:=\{s_1,s_2,...,s_n\}$. 

\vspace{2mm}\noindent In an abuse of notation, we interpret the Davis complex $\Sigma$ as a complex of groups with all vertex groups trivial over itself as a scwol. We colour $L$, $K$ and $\Sigma$ by assigning the colour purple to all vertices of type $\varnothing$, green to all vertices of type $s_1$, $s_2$ or $s_3$, blue to all vertices of type $s_4$, $s_5$, ... or $s_n$ and red to all vertices of type $\{s_i,s_j\}$ where $4 \leq i,j \leq n$ and $i \neq j$. Equip $\Sigma$ with the $\CAT(0)$ metric introduced by Moussong in Theorem~\ref{hjh}.

\begin{Definition}\label{thingy}
\vspace{2mm}\noindent Take any integers $\alpha,\beta \geq 2$. Let $(W^{\alpha,\beta},\{t_1,...,t_{2\beta}\})$ be the Coxeter system with nerve $L^{\alpha,\beta}$ the $2\beta$-cycle graph with constant edge labels $\alpha$, and let $K^{\alpha,\beta}$ be the corresponding chamber. Let $\mathcal P^{\alpha,\beta}$ be the $2\beta$-gon of groups where each vertex group is $\mathcal D_{2\alpha}$, each edge group is $\mathcal C_2$, the face group is trivial and all monomorphisms are natural inclusions. Observe that $\mathcal P^{\alpha,\beta}$ is the complex of groups associated with the natural action of $(W^{\alpha,\beta},\{t_1,...,t_{2\beta}\})$ on its Davis complex $\Sigma^{\alpha,\beta}$ (see Section $2.3$ of \cite{T2}). That is, $\mathcal P^{\alpha,\beta}$ is strictly developable. Each polygon in $\mathcal P^{\alpha,\beta}$ has a type induced by that of the underlying scwol $K^{\alpha,\beta}$. Colour the vertices of $K^{\alpha,\beta}$, $\Sigma^{\alpha,\beta}$ and $\mathcal P^{\alpha,\beta}$ as follows. Assign the colour red to all vertices of type $\varnothing$, blue to all vertices of type $t_1$, $t_2$, ... or $t_{2\beta}$ and purple to all vertices of type $\{t_i,t_j\}$ for any $1 \leq i,j \leq 2\beta$ and $i \neq j$. As always, equip $\Sigma^{\alpha,\beta}$ with Moussong's metric.
\end{Definition}

\noindent Now consider the family $\{\mathcal X_k^{(0)}\}_{k \in \N} \cup \{ \mathcal X^{(0)}\}$ of complexes of groups shown in Figures~\ref{ggggg} and ~\ref{ggggg'}. Each complex of groups in this family consists of "multiples" of copies of $\mathcal P^{m,x}$ stacked on top of one another, glued together via additional green vertices. For each $i \in \{0,1,...,k\}$ (resp. $j\in \{0,1\}$), and for any integer $k$, we call the $i+1$'th (resp. $j+1$'th) lowest embedded "multiple" of $\mathcal P^{m,x}$ the \textit{$i$'th platform} $\mathcal P_i^{m,x} = \mathcal C_2^i \times \mathcal P^{m,x}$ (resp. \textit{$j$'th platform} $\mathcal P_j^{m,x} = \mathcal P^{m,x}$) of the structure.

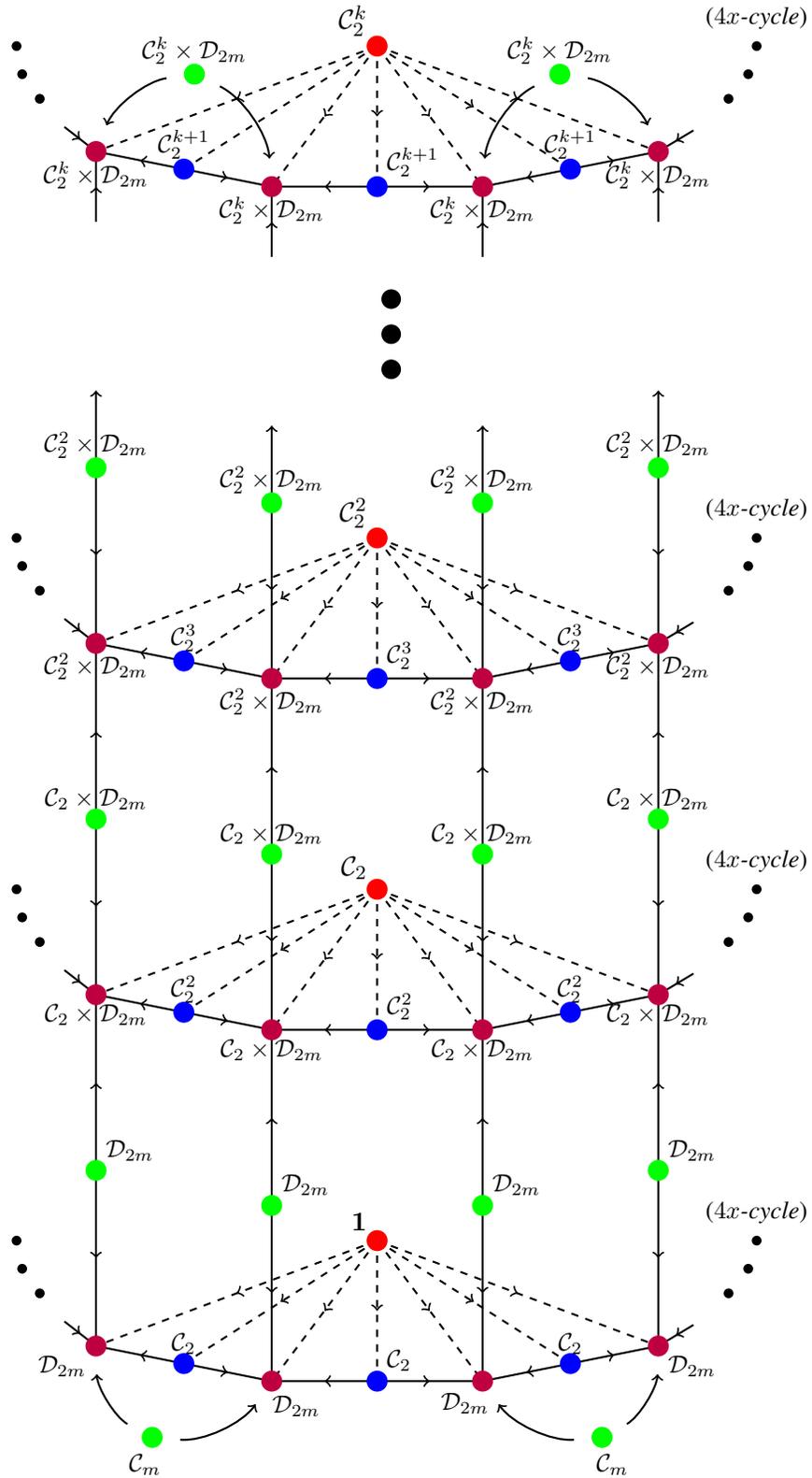
\begin{figure} \
\centering
\usetikzlibrary{arrows,positioning}
\begin{tikzpicture} [%
    purple/.style = {circle,fill=purple,text=black,inner sep=3pt},
    blue/.style = {circle,fill=blue,text=black,inner sep=3pt},
    red/.style = {circle,fill=red,text=black,inner sep=3pt},
   green/.style = {circle,fill=green,text=black,inner sep=3pt}, 
    line/.style={->,shorten >=0.4cm,shorten <=0.4cm},thick]

\draw [dashed, black, thick,fill=none,->-] (0,0) node [above left, black] {$\boldsymbol{1}$} -- (-1.5,-2) node {};
\node[] at (-1.15,-2.06) [below, black]{\small $\mathcal D_{2m}$};
\draw [dashed, black, thick,fill=none,->-] (0,0) node {} -- (1.5,-2) node {};
\node[] at (1.25,-2.03) [below, black]{\small $\mathcal D_{2m}$};
\draw [dashed, black, thick,fill=none,->-] (0,0) node {} -- (4,-1.5) node [below right, black] {\small $\mathcal D_{2m}$};
\draw [dashed, black, thick,fill=none,->-] (0,0) node {} -- (-4,-1.5) node [below left, black] {\small $\mathcal D_{2m}$};
\draw [dashed, black, thick,fill=none,->-] (0,0) node {} -- (0,-2) node [above right, black] {\small $\mathcal C_{2}$};
\draw [dashed, black, thick,fill=none,->-] (0,0) node {} -- (2.75,-1.75) node [above, black] {\small $\mathcal C_{2}$};
\draw [dashed, black, thick,fill=none,->-] (0,0) node {} -- (-2.75,-1.75) node [above, black] {\small $\mathcal C_{2}$};
\draw [black, thick,fill=none,->-] (0,-2) node {} -- (-1.5,-2) node {};
\draw [black, thick,fill=none,->-] (0,-2) node {} -- (1.5,-2) node {};
\draw [black, thick,fill=none,->-] (2.75,-1.75) node {} -- (1.5,-2) node {};
\draw [black, thick,fill=none,->-] (2.75,-1.75) node {} -- (4,-1.5) node {};
\draw [black, thick,fill=none,->-] (-2.75,-1.75) node {} -- (-1.5,-2) node {};
\draw [black, thick,fill=none,->-] (-2.75,-1.75) node {} -- (-4,-1.5) node {};
\draw [black, thick,fill=none,->-]  (4.5,-1.2) node {} -- (4,-1.5) node {};
\draw [black, thick,fill=none,->-]  (-4.45,-1.15) node {} -- (-4,-1.5) node {};
\node[circle,fill=black,text=black,inner sep=1.4pt] at (5,-0.75) {};
\node[circle,fill=black,text=black,inner sep=1.4pt] at (5.28,-0.4) {};
\node[circle,fill=black,text=black,inner sep=1.4pt] at (5.4,0) {};
\node[circle,fill=black,text=black,inner sep=1.4pt] at (-4.8,-0.75) {};
\node[circle,fill=black,text=black,inner sep=1.4pt] at (-5.05,-0.4) {};
\node[circle,fill=black,text=black,inner sep=1.4pt] at (-5.13,0) {};
\draw [black, thick,fill=none,->-] (-1.5,0.5) node [above right, black] {\small $\mathcal D_{2m}$} -- (-1.5,-2) node {};
\draw [black, thick,fill=none,->-] (1.5,0.5) node [above right, black] {\small $\mathcal D_{2m}$} -- (1.5,-2) node {};
\draw [black, thick,fill=none,->-] (4,1) node [above right, black] {\small $\mathcal D_{2m}$} -- (4,-1.5) node {};
\draw [black, thick,fill=none,->-] (-4,1) node [above right, black] {\small $\mathcal D_{2 m}$} -- (-4,-1.5) node {};
\draw [black, thick,fill=none,->-] (-1.5,0.5) node {} -- (-1.5,3) node {};
\draw [black, thick,fill=none,->-] (1.5,0.5) node {} -- (1.5,3) node {};
\draw [black, thick,fill=none,->-] (4,1) node {} -- (4,3.5) node {};
\draw [black, thick,fill=none,->-] (-4,1) node {} -- (-4,3.5) node {};
\node[red] at (0,0) {};
\node[purple] at (-1.5,-2) {};
\node[purple] at (1.5,-2) {};
\node[purple] at (4,-1.5) {};
\node[purple] at (-4,-1.5) {};
\node[blue] at (0,-2) {};
\node[blue] at (2.75,-1.75) {};
\node[blue] at (-2.75,-1.75) {};
\node[green] at (-1.5,0.5) {};
\node[green] at (1.5,0.5) {};
\node[green] at (4,1) {};
\node[green] at (-4,1) {};

\draw [dashed, black, thick,fill=none,->-] (0,5) node [above left, black] {$\mathcal C _2$} -- (-1.5,3) node [below, black] {\small $\mathcal C_2 \times \mathcal D_{2 m}$};
\draw [dashed, black, thick,fill=none,->-] (0,5) node {} -- (1.5,3) node [below, black] {\small $\mathcal C_2 \times \mathcal D_{2 m}$};
\draw [dashed, black, thick,fill=none,->-] (0,5) node {} -- (4,3.5) node [below, black] {\small $\mathcal C_2 \times \mathcal D_{2 m}$};
\draw [dashed, black, thick,fill=none,->-] (0,5) node {} -- (-4,3.5) node [below, black] {\small $\mathcal C_2 \times \mathcal D_{2m}$};
\draw [dashed, black, thick,fill=none,->-] (0,5) node {} -- (0,3) node [above right, black] {\small $\mathcal C_{2}^2$};
\draw [dashed, black, thick,fill=none,->-] (0,5) node {} -- (2.75,3.25) node [above, black] {\small $\mathcal C_{2}^2$};
\draw [dashed, black, thick,fill=none,->-] (0,5) node {} -- (-2.75,3.25) node [above, black] {\small $\mathcal C_{2}^2$};
\draw [black, thick,fill=none,->-] (0,3) node {} -- (-1.5,3) node {};
\draw [black, thick,fill=none,->-] (0,3) node {} -- (1.5,3) node {};
\draw [black, thick,fill=none,->-] (2.75,3.25) node {} -- (1.5,3) node {};
\draw [black, thick,fill=none,->-] (2.75,3.25) node {} -- (4,3.5) node {};
\draw [black, thick,fill=none,->-] (-2.75,3.25) node {} -- (-1.5,3) node {};
\draw [black, thick,fill=none,->-] (-2.75,3.25) node {} -- (-4,3.5) node {};
\draw [black, thick,fill=none,->-]  (4.5,3.8) node {} -- (4,3.5) node {};
\draw [black, thick,fill=none,->-]  (-4.45,3.85) node {} -- (-4,3.5) node {};
\node[circle,fill=black,text=black,inner sep=1.4pt] at (5,4.25) {};
\node[circle,fill=black,text=black,inner sep=1.4pt] at (5.28,4.6) {};
\node[circle,fill=black,text=black,inner sep=1.4pt] at (5.4,5) {};
\node[circle,fill=black,text=black,inner sep=1.4pt] at (-4.8,4.25) {};
\node[circle,fill=black,text=black,inner sep=1.4pt] at (-5.05,4.6) {};
\node[circle,fill=black,text=black,inner sep=1.4pt] at (-5.13,5) {};
\draw [black, thick,fill=none,->-] (-1.5,5.5) node [above, black] {\small $\mathcal C_2 \times \mathcal D_{2 m}$} -- (-1.5,3) node {};
\draw [black, thick,fill=none,->-] (1.5,5.5) node [above, black] {\small $\mathcal C_2 \times \mathcal D_{2 m}$} -- (1.5,3) node {};
\draw [black, thick,fill=none,->-] (4,6) node [above, black] {\small $\mathcal C_2 \times \mathcal D_{2 m}$} -- (4,3.5) node {};
\draw [black, thick,fill=none,->-] (-4,6) node [above, black] {\small $\mathcal C_2 \times \mathcal D_{2 m}$} -- (-4,3.5) node {};
\draw [black, thick,fill=none,->-] (-1.5,5.5) node {} -- (-1.5,8) node {};
\draw [black, thick,fill=none,->-] (1.5,5.5) node {} -- (1.5,8) node {};
\draw [black, thick,fill=none,->-] (4,6) node {} -- (4,8.5) node {};
\draw [black, thick,fill=none,->-] (-4,6) node {} -- (-4,8.5) node {};
\node[red] at (0,5) {};
\node[purple] at (-1.5,3) {};
\node[purple] at (1.5,3) {};
\node[purple] at (4,3.5) {};
\node[purple] at (-4,3.5) {};
\node[blue] at (0,3) {};
\node[blue] at (2.75,3.25) {};
\node[blue] at (-2.75,3.25) {};
\node[green] at (-1.5,5.5) {};
\node[green] at (1.5,5.5) {};
\node[green] at (4,6) {};
\node[green] at (-4,6) {};

\draw [dashed, black, thick,fill=none,->-] (0,10) node [above left, black] {$\mathcal C^2 _2$} -- (-1.5,8) node [below, black] {\small $\mathcal C_2^2 \times \mathcal D_{2 m}$};
\draw [dashed, black, thick,fill=none,->-] (0,10) node {} -- (1.5,8) node [below, black] {\small $\mathcal C_2^2 \times \mathcal D_{2 m}$};
\draw [dashed, black, thick,fill=none,->-] (0,10) node {} -- (4,8.5) node [below, black] {\small $\mathcal C_2^2 \times \mathcal D_{2m}$};
\draw [dashed, black, thick,fill=none,->-] (0,10) node {} -- (-4,8.5) node [below, black] {\small $\mathcal C_2^2 \times \mathcal D_{2 m}$};
\draw [dashed, black, thick,fill=none,->-] (0,10) node {} -- (0,8) node [above right, black] {\small $\mathcal C_{2}^3$};
\draw [dashed, black, thick,fill=none,->-] (0,10) node {} -- (2.75,8.25) node [above, black] {\small $\mathcal C_{2}^3$};
\draw [dashed, black, thick,fill=none,->-] (0,10) node {} -- (-2.75,8.25) node [above, black] {\small $\mathcal C_{2}^3$};
\draw [black, thick,fill=none,->-] (0,8) node {} -- (-1.5,8) node {};
\draw [black, thick,fill=none,->-] (0,8) node {} -- (1.5,8) node {};
\draw [black, thick,fill=none,->-] (2.75,8.25) node {} -- (1.5,8) node {};
\draw [black, thick,fill=none,->-] (2.75,8.25) node {} -- (4,8.5) node {};
\draw [black, thick,fill=none,->-] (-2.75,8.25) node {} -- (-1.5,8) node {};
\draw [black, thick,fill=none,->-] (-2.75,8.25) node {} -- (-4,8.5) node {};
\draw [black, thick,fill=none,->-]  (4.5,8.8) node {} -- (4,8.5) node {};
\draw [black, thick,fill=none,->-]  (-4.45,8.85) node {} -- (-4,8.5) node {};
\node[circle,fill=black,text=black,inner sep=1.4pt] at (5,9.25) {};
\node[circle,fill=black,text=black,inner sep=1.4pt] at (5.28,9.6) {};
\node[circle,fill=black,text=black,inner sep=1.4pt] at (5.4,10) {};
\node[circle,fill=black,text=black,inner sep=1.4pt] at (-4.8,9.25) {};
\node[circle,fill=black,text=black,inner sep=1.4pt] at (-5.05,9.6) {};
\node[circle,fill=black,text=black,inner sep=1.4pt] at (-5.13,10) {};
\draw [black, thick,fill=none,->-] (-1.5,10.5) node [above, black] {\small $\mathcal C_2^2 \times \mathcal D_{2 m}$} -- (-1.5,8) node {};
\draw [black, thick,fill=none,->-] (1.5,10.5) node [above, black] {\small $\mathcal C_2^2 \times \mathcal D_{2  m}$} -- (1.5,8) node {};
\draw [black, thick,fill=none,->-] (4,11) node [above, black] {\small $\mathcal C_2^2 \times \mathcal D_{2 m}$} -- (4,8.5) node {};
\draw [black, thick,fill=none,->-] (-4,11) node [above, black] {\small $\mathcal C_2^2 \times \mathcal D_{2 m}$} -- (-4,8.5) node {};
\draw [black, thick,fill=none,->] (-1.5,10.5) node {} -- (-1.5,11.6) node {};
\draw [black, thick,fill=none,->] (1.5,10.5) node {} -- (1.5,11.6) node {};
\draw [black, thick,fill=none,->] (4,11) node {} -- (4,12.1) node {};
\draw [black, thick,fill=none,->] (-4,11) node {} -- (-4,12.1) node {};
\node[red] at (0,10) {};
\node[purple] at (-1.5,8) {};
\node[purple] at (1.5,8) {};
\node[purple] at (4,8.5) {};
\node[purple] at (-4,8.5) {};
\node[blue] at (0,8) {};
\node[blue] at (2.75,8.25) {};
\node[blue] at (-2.75,8.25) {};
\node[green] at (-1.5,10.5) {};
\node[green] at (1.5,10.5) {};
\node[green] at (4,11) {};
\node[green] at (-4,11) {};

\draw [black, thick,fill=none,->-] (-1.5,14) node {} -- (-1.5,15) node {};
\draw [black, thick,fill=none,->-] (1.5,14) node {} -- (1.5,15) node {};
\draw [black, thick,fill=none,->-] (-4,14.5) node {} -- (-4,15.5) node {};
\draw [black, thick,fill=none,->-] (4,14.5) node {} -- (4,15.5) node {};

\draw [dashed, black, thick,fill=none,->-] (0,17) node [above left, black] {$\mathcal C^k _2$} -- (-1.5,15) node [below, black] {\small $\mathcal C_2^k \times \mathcal D_{2 m}$};
\draw [dashed, black, thick,fill=none,->-] (0,17) node {} -- (1.5,15) node [below, black] {\small $\mathcal C_2^k \times \mathcal D_{2  m}$};
\draw [dashed, black, thick,fill=none,->-] (0,17) node {} -- (4,15.5) node [below, black] {\small $\mathcal C_2^k \times \mathcal D_{2  m}$};
\draw [dashed, black, thick,fill=none,->-] (0,17) node {} -- (-4,15.5) node [below, black] {\small $\mathcal C_2^k \times \mathcal D_{2 m}$};
\draw [dashed, black, thick,fill=none,->-] (0,17) node {} -- (0,15) node [above right, black] {\small $\mathcal C_{2}^{k+1}$};
\draw [dashed, black, thick,fill=none,->-] (0,17) node {} -- (2.75,15.25) node [above, black] {\small $\mathcal C_{2}^{k+1}$};
\draw [dashed, black, thick,fill=none,->] (0,17) node {} -- (-2.75,15.25) node [above, black] {\small $\mathcal C_{2}^{k+1}$};
\draw [black, thick,fill=none,->-] (0,15) node {} -- (-1.5,15) node {};
\draw [black, thick,fill=none,->-] (0,15) node {} -- (1.5,15) node {};
\draw [black, thick,fill=none,->-] (2.75,15.25) node {} -- (1.5,15) node {};
\draw [black, thick,fill=none,->-] (2.75,15.25) node {} -- (4,15.5) node {};
\draw [black, thick,fill=none,->-] (-2.75,15.25) node {} -- (-1.5,15) node {};
\draw [black, thick,fill=none,->-] (-2.75,15.25) node {} -- (-4,15.5) node {};
\draw [black, thick,fill=none,->-]  (4.5,15.8) node {} -- (4,15.5) node {};
\draw [black, thick,fill=none,->-]  (-4.45,15.85) node {} -- (-4,15.5) node {};
\node[circle,fill=black,text=black,inner sep=1.4pt] at (5,16.25) {};
\node[circle,fill=black,text=black,inner sep=1.4pt] at (5.28,16.6) {};
\node[circle,fill=black,text=black,inner sep=1.4pt] at (5.4,17) {};
\node[circle,fill=black,text=black,inner sep=1.4pt] at (-4.8,16.25) {};
\node[circle,fill=black,text=black,inner sep=1.4pt] at (-5.05,16.6) {};
\node[circle,fill=black,text=black,inner sep=1.4pt] at (-5.13,17) {};
\node[red] at (0,17) {};
\node[purple] at (-1.5,15) {};
\node[purple] at (1.5,15) {};
\node[purple] at (4,15.5) {};
\node[purple] at (-4,15.5) {};
\node[blue] at (0,15) {};
\node[blue] at (2.75,15.25) {};
\node[blue] at (-2.75,15.25) {};
\node[circle,fill=black,text=black,inner sep=2.8pt] at (0.2,12.4) {};
\node[circle,fill=black,text=black,inner sep=2.8pt] at (0.2,12.9) {};
\node[circle,fill=black,text=black,inner sep=2.8pt] at (0.2,13.4) {};

\path [black,thick,line,bend right] (-3.2,-2.8) edge (-1.5,-2);
\path [black,thick,line,bend left] (3.2,-2.8) edge (1.5,-2);
\path [black,thick,line,bend left] (-3.2,-2.8) edge (-4,-1.5);
\path [black,thick,line,bend right] (3.2,-2.8) edge (4,-1.5);
\node[green] at (-3.2,-2.8) {};
\node[green] at (3.2,-2.8) {};
\node[] at (-3.3,-2.9) [below, black]{\small $\mathcal C_m$};
\node[] at (3.35,-2.9) [below, black]{\small $\mathcal C_m$};

\path [black,thick,line,bend left] (-2.6,16.6) edge (-1.5,15);
\path [black,thick,line,bend right] (2.6,16.6) edge (1.5,15);
\path [black,thick,line,bend right] (-2.6,16.6) edge (-4,15.5);
\path [black,thick,line,bend left] (2.6,16.6) edge (4,15.5);
\node[green] at (-2.6,16.6) {};
\node[green] at (2.6,16.6) {};
\node[] at (-2.6,16.6) [above, black]{\small $\mathcal C_2^k \times \mathcal D_{2  m}$};
\node[] at (2.6,16.6) [above, black]{\small $\mathcal C_2^k \times \mathcal D_{2  m}$};

\node[] at (5.4,0.4) [black]{\small ($4x$\textit{-cycle})};
\node[] at (5.4,5.4) [black]{\small ($4x$\textit{-cycle})};
\node[] at (5.4,10.4) [black]{\small ($4x$\textit{-cycle})};
\node[] at (5.4,17.4) [black]{\small ($4x$\textit{-cycle})};
\end{tikzpicture}
\caption{A family $\{\mathcal X_k^{(0)}\}_{k \in \N}$ of finite complexes of finite groups} (with embedded "platforms" $\mathcal P_i^{m,x} = \mathcal C_2^i \times \mathcal P^{m,x}$ for $i=0,1,...,k$)\label{ggggg}
\end{figure}

\begin{figure}[h] \
\centering
\usetikzlibrary{arrows,positioning}
\begin{tikzpicture} [%
    purple/.style = {circle,fill=purple,text=black,inner sep=3pt},
    blue/.style = {circle,fill=blue,text=black,inner sep=3pt},
    red/.style = {circle,fill=red,text=black,inner sep=3pt},
   green/.style = {circle,fill=green,text=black,inner sep=3pt}, 
    line/.style={->,shorten >=0.4cm,shorten <=0.4cm},thick]

\draw [dashed, black, thick,fill=none,->-] (0,0) node [above left, black] {$\boldsymbol{1}$} -- (-1.5,-2) node [below, black] {\small $\mathcal D_{2 m}$};
\draw [dashed, black, thick,fill=none,->-] (0,0) node {} -- (1.5,-2) node [below, black] {\small $\mathcal D_{2 m}$};
\draw [dashed, black, thick,fill=none,->-] (0,0) node {} -- (4,-1.5) node [below, black] {\small $\mathcal D_{2  m}$};
\draw [dashed, black, thick,fill=none,->-] (0,0) node {} -- (-4,-1.5) node [below, black] {\small $\mathcal D_{2  m}$};
\draw [dashed, black, thick,fill=none,->-] (0,0) node {} -- (0,-2) node [above right, black] {\small $\mathcal C_{2}$};
\draw [dashed, black, thick,fill=none,->-] (0,0) node {} -- (2.75,-1.75) node [above, black] {\small $\mathcal C_{2}$};
\draw [dashed, black, thick,fill=none,->-] (0,0) node {} -- (-2.75,-1.75) node [above, black] {\small $\mathcal C_{2}$};
\draw [black, thick,fill=none,->-] (0,-2) node {} -- (-1.5,-2) node {};
\draw [black, thick,fill=none,->-] (0,-2) node {} -- (1.5,-2) node {};
\draw [black, thick,fill=none,->-] (2.75,-1.75) node {} -- (1.5,-2) node {};
\draw [black, thick,fill=none,->-] (2.75,-1.75) node {} -- (4,-1.5) node {};
\draw [black, thick,fill=none,->-] (-2.75,-1.75) node {} -- (-1.5,-2) node {};
\draw [black, thick,fill=none,->-] (-2.75,-1.75) node {} -- (-4,-1.5) node {};
\draw [black, thick,fill=none,->-]  (4.5,-1.2) node {} -- (4,-1.5) node {};
\draw [black, thick,fill=none,->-]  (-4.45,-1.15) node {} -- (-4,-1.5) node {};
\node[circle,fill=black,text=black,inner sep=1.4pt] at (5,-0.75) {};
\node[circle,fill=black,text=black,inner sep=1.4pt] at (5.28,-0.4) {};
\node[circle,fill=black,text=black,inner sep=1.4pt] at (5.4,0) {};
\node[circle,fill=black,text=black,inner sep=1.4pt] at (-4.8,-0.75) {};
\node[circle,fill=black,text=black,inner sep=1.4pt] at (-5.05,-0.4) {};
\node[circle,fill=black,text=black,inner sep=1.4pt] at (-5.13,0) {};
\draw [black, thick,fill=none,->-] (-1.8,0.5) node [above left, black] {\small $\mathcal D_{2  m}$} -- (-1.5,-2) node {};
\draw [black, thick,fill=none,->-] (1.8,0.5) node [above right, black] {\small $\mathcal C_{m}$} -- (1.5,-2) node {};
\draw [black, thick,fill=none,->-] (4.4,1) node [above right, black] {\small $\mathcal C_{m}$} -- (4,-1.5) node {};
\draw [black, thick,fill=none,->-] (-4.4,1) node [above left, black] {\small $\mathcal D_{2 m}$} -- (-4,-1.5) node {};
\draw [black, thick,fill=none,->-] (-1.8,0.5) node {} -- (-1.5,3) node {};
\draw [black, thick,fill=none,->-] (1.8,0.5) node {} -- (1.5,3) node {};
\draw [black, thick,fill=none,->-] (4.4,1) node {} -- (4,3.5) node {};
\draw [black, thick,fill=none,->-] (-4.4,1) node {} -- (-4,3.5) node {};
\draw [black, thick,fill=none,->-] (-1.2,0.5) node [above right, black] {\small $\mathcal C_{m}$} -- (-1.5,-2) node {};
\draw [black, thick,fill=none,->-] (1.2,0.5) node [above left, black] {\small $\mathcal D_{2 m}$} -- (1.5,-2) node {};
\draw [black, thick,fill=none,->-] (3.6,1) node [above left, black] {\small $\mathcal D_{2 m}$} -- (4,-1.5) node {};
\draw [black, thick,fill=none,->-] (-3.6,1) node [above right, black] {\small $\mathcal C_{m}$} -- (-4,-1.5) node {};
\draw [black, thick,fill=none,->-] (-1.2,0.5) node {} -- (-1.5,3) node {};
\draw [black, thick,fill=none,->-] (1.2,0.5) node {} -- (1.5,3) node {};
\draw [black, thick,fill=none,->-] (3.6,1) node {} -- (4,3.5) node {};
\draw [black, thick,fill=none,->-] (-3.6,1) node {} -- (-4,3.5) node {};
\node[red] at (0,0) {};
\node[purple] at (-1.5,-2) {};
\node[purple] at (1.5,-2) {};
\node[purple] at (4,-1.5) {};
\node[purple] at (-4,-1.5) {};
\node[blue] at (0,-2) {};
\node[blue] at (2.75,-1.75) {};
\node[blue] at (-2.75,-1.75) {};
\node[green] at (-1.8,0.5) {};
\node[green] at (1.8,0.5) {};
\node[green] at (4.4,1) {};
\node[green] at (-4.4,1) {};
\node[green] at (-1.2,0.5) {};
\node[green] at (1.2,0.5) {};
\node[green] at (3.6,1) {};
\node[green] at (-3.6,1) {};

\draw [dashed, black, thick,fill=none,->-] (0,5) node [above left, black] {$\boldsymbol{1}$} -- (-1.5,3) node [below left, black] {\small $\mathcal D_{2 m}$};
\draw [dashed, black, thick,fill=none,->-] (0,5) node {} -- (1.5,3) node [below right, black] {\small $\mathcal D_{2 m}$};
\draw [dashed, black, thick,fill=none,->-] (0,5) node {} -- (4,3.5) node [below right, black] {\small $\mathcal D_{2 m}$};
\draw [dashed, black, thick,fill=none,->-] (0,5) node {} -- (-4,3.5) node [below left, black] {\small $\mathcal D_{2 m}$};
\draw [dashed, black, thick,fill=none,->-] (0,5) node {} -- (0,3) node [above right, black] {\small $\mathcal C_2$};
\draw [dashed, black, thick,fill=none,->-] (0,5) node {} -- (2.75,3.25) node [above, black] {\small $\mathcal C_2$};
\draw [dashed, black, thick,fill=none,->-] (0,5) node {} -- (-2.75,3.25) node [above, black] {\small $\mathcal C_2$};
\draw [black, thick,fill=none,->-] (0,3) node {} -- (-1.5,3) node {};
\draw [black, thick,fill=none,->-] (0,3) node {} -- (1.5,3) node {};
\draw [black, thick,fill=none,->-] (2.75,3.25) node {} -- (1.5,3) node {};
\draw [black, thick,fill=none,->-] (2.75,3.25) node {} -- (4,3.5) node {};
\draw [black, thick,fill=none,->-] (-2.75,3.25) node {} -- (-1.5,3) node {};
\draw [black, thick,fill=none,->-] (-2.75,3.25) node {} -- (-4,3.5) node {};
\draw [black, thick,fill=none,->-]  (4.5,3.8) node {} -- (4,3.5) node {};
\draw [black, thick,fill=none,->-]  (-4.45,3.85) node {} -- (-4,3.5) node {};
\node[circle,fill=black,text=black,inner sep=1.4pt] at (5,4.25) {};
\node[circle,fill=black,text=black,inner sep=1.4pt] at (5.28,4.6) {};
\node[circle,fill=black,text=black,inner sep=1.4pt] at (5.4,5) {};
\node[circle,fill=black,text=black,inner sep=1.4pt] at (-4.8,4.25) {};
\node[circle,fill=black,text=black,inner sep=1.4pt] at (-5.05,4.6) {};
\node[circle,fill=black,text=black,inner sep=1.4pt] at (-5.13,5) {};
\node[red] at (0,5) {};
\node[purple] at (-1.5,3) {};
\node[purple] at (1.5,3) {};
\node[purple] at (4,3.5) {};
\node[purple] at (-4,3.5) {};
\node[blue] at (0,3) {};
\node[blue] at (2.75,3.25) {};
\node[blue] at (-2.75,3.25) {};

\node[] at (5.4,0.4) [black]{\small ($4x$\textit{-cycle})};
\node[] at (5.4,5.4) [black]{\small ($4x$\textit{-cycle})};
\end{tikzpicture}
\caption{A finite complex of finite groups $\mathcal X^{(0)}$} \label{ggggg'}
\end{figure}
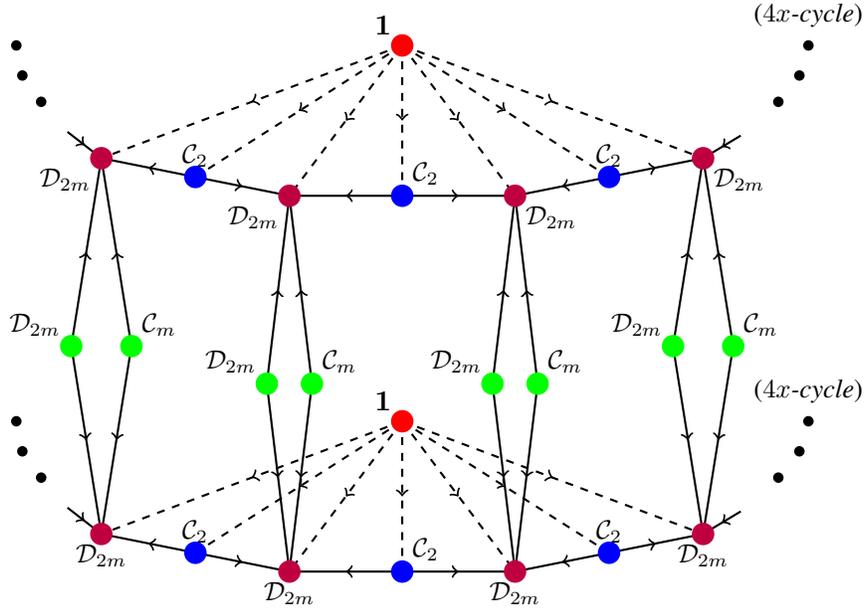

\begin{Proposition}\label{links} Let $\mathcal F$ be any complex of groups in the family $\{\mathcal X_k^{(0)}\}_{k \in \N} \cup \{ \mathcal X^{(0)}\}$. Then the set of isometry classes of links of vertices in $\mathcal F$ is the same as that of $\Sigma$, as illustrated in Figure~\ref{hrr}, and in particular is consistent with the colour schemes of $\mathcal F$ and $\Sigma$. That is, two vertices in either $\mathcal F$ or $\Sigma$ have isometric links iff they have the same colour.

\begin{proof} Consider the Coxeter system $(W^{m,x},\{t_1,...,t_{2x}\})$, its Davis complex $\Sigma^{m,x}$ and the associated complex of groups $\mathcal P^{m,x}$ (as in Definition~\ref{thingy}). Observe that the isometry class of link of each polygon in $\mathcal P^{m,x}$ is the same as that of any vertex of the same type in $\Sigma^{m,x}$ since the universal covering map is a type-preserving local isometry.

\vspace{2mm}\noindent For any distinct $i,j \in \{1,...,2x\}$ let $\phi_{i,j}$ be the automorphism of $(W^{m,x},\{t_1,...,t_{2x}\})$ that swaps the generators $t_i$ and $t_j$. Observe that $\phi_{i,j}$ induces an automorphism on $\Sigma^{m,x}$ for all $i,j \in \{1,...,2x\}$. Exploiting this symmetry means that there are only three inequivalent isometry classes of links of vertices in $\Sigma^{m,x}$; that of vertices of type $\varnothing$ (red vertices), of type $t_i$ (blue vertices) and of type $\{t_i,t_j\}$ for $i \neq j$ (purple vertices).

\vspace{2mm}\noindent Interpret $\Sigma^{m,x}$ as a complex of groups with all vertex groups trivial over itself as a scwol. Case I of Example~\ref{helpies} tells us that the link of any vertex $v$ in $\Sigma^{m,x}$ is the intersection of $\Sigma^{m,x}$ (as a polyhedral complex) with a $2$-sphere of sufficiently small radius centred at $v$. In particular, the (combinatorial) link of any red vertex in $\Sigma^{m,x}$ is the $4x$-cycle graph with vertices alternating between purple and blue, the link of any blue vertex in $\Sigma^{m,x}$ is the $4$-cycle graph with vertices alternating between purple and red, and the link of any purple vertex in $\Sigma^{m,x}$ is the $4m$-cycle graph with vertices alternating between blue and red.

\vspace{2mm}\noindent Now consider the complex of groups $\mathcal F$. The universal cover of each of its embedded platforms $\mathcal P_i^{m,x}$ is $\Sigma^{m,x}$ by the same argument as for $\mathcal P^{m,x}$. The link of any red (resp. blue) vertex in $\mathcal F$ is the same as that of any red (resp. blue) vertex in $\mathcal P^{m,x}$. It remains to show that the link of each green and purple vertex in $\mathcal F$ is as in Figure~\ref{hrr}. Let $v_g$ (resp. $v_p$) be (without loss of generality) any green (resp. purple) vertex in $\mathcal F$.

\vspace{2mm}\noindent Denote the scwol underlying $\mathcal F$ by $X_{\mathcal F}$. Since both of the edges adjacent to $v_g$ in $X_{\mathcal F}$ are oriented towards a purple vertex, Case II of Example~\ref{helpies} tells us that $\Lk(v_g)$ is the graph with two purple vertices and no edges.

\vspace{2mm}\noindent Let $\Phi_p$ be the poset $\{v_{red}, v_{blue},v'_{blue},v_{green},v'_{green}\}$ of all vertices in $X_{\mathcal F}$ adjacent to $v_p$ with ordering $v_{blue} \prec v_{red}$ and $v'_{blue} \prec v_{red}$. Observe that the geometric realisation of $\Phi_p$ has three connected components, one of which corresponds to the poset of vertices adjacent to any purple vertex in $\mathcal P^{m,x}$, and the other two of which are lone green vertices with local groups of index $1$ and $2$ respectively in $\mathcal F_{v_p}$. Cases III and IV of Example~\ref{helpies} then tell us that $\Lk(v_p)$ of $\mathcal F$ is the disjoint union of the link of any purple vertex in $\mathcal P^{m,x}$ with three disconnected green vertices.
\end{proof} 
\end{Proposition}

\begin{figure}[h] \
\centering
\usetikzlibrary{arrows,positioning}
\begin{tikzpicture} [%
    purple/.style = {circle,fill=purple,text=black,inner sep=3pt},
    blue/.style = {circle,fill=blue,text=black,inner sep=3pt},
    red/.style = {circle,fill=red,text=black,inner sep=3pt},
   green/.style = {circle,fill=green,text=black,inner sep=3pt}, 
    line/.style={->,shorten >=0.4cm,shorten <=0.4cm},thick]

\draw [black, thick,fill=none,-] (-2.25,0) node {} -- (-2.25,-0.8) node {};
\draw [black, thick,fill=none,-] (-2.25,0) node {} -- (-2.25,0.8) node {};
\draw [black, thick,fill=none,-] (-1.85,1.45) node {} -- (-2.25,0.8) node {};
\draw [black, thick,fill=none,-] (-1.85,1.45) node {} -- (-1.45,2.1) node {};
\draw [black, thick,fill=none,-] (-1.85,-1.45) node {} -- (-2.25,-0.8) node {};
\draw [black, thick,fill=none,-] (-1.85,-1.45) node {} -- (-1.45,-2.1) node {};
\draw [black, thick,fill=none,-]  (-1.45,2.1) node {} -- (-1.05,2.3) node {};
\draw [black, thick,fill=none,-] (-1.45,-2.1) node {} -- (-1.05,-2.3) node {};
\node[circle,fill=black,text=black,inner sep=1pt] at (0,2.5) {};
\node[circle,fill=black,text=black,inner sep=1pt] at (0.25,2.48) {};
\node[circle,fill=black,text=black,inner sep=1pt] at (-0.25,2.48) {};
\node[circle,fill=black,text=black,inner sep=1pt] at (0,-2.5) {};
\node[circle,fill=black,text=black,inner sep=1pt] at (0.25,-2.48) {};
\node[circle,fill=black,text=black,inner sep=1pt] at (-0.25,-2.48) {};
\draw [black, thick,fill=none,-] (2.25,0) node {} -- (2.25,-0.8) node {};
\draw [black, thick,fill=none,-] (2.25,0) node {} -- (2.25,0.8) node {};
\draw [black, thick,fill=none,-] (1.85,1.45) node {} -- (2.25,0.8) node {};
\draw [black, thick,fill=none,-] (1.85,1.45) node {} -- (1.45,2.1) node {};
\draw [black, thick,fill=none,-] (1.85,-1.45) node {} -- (2.25,-0.8) node {};
\draw [black, thick,fill=none,-] (1.85,-1.45) node {} -- (1.45,-2.1) node {};
\draw [black, thick,fill=none,-]  (1.45,2.1) node {} -- (1.05,2.3) node {};
\draw [black, thick,fill=none,-] (1.45,-2.1) node {} -- (1.05,-2.3) node {};
\node[green] at (-0.45,-0.333) {};
\node[green] at (0,0.5) {};
\node[green] at (0.45,-0.333) {};
\node[red] at (-2.25,-0.8) {};
\node[red] at (-2.25,0.8) {};
\node[red] at (-1.45,2.1) {};
\node[red] at (-1.45,-2.1) {};
\node[blue] at (-2.25,0) {};
\node[blue] at (-1.85,1.45) {};
\node[blue] at (-1.85,-1.45) {};
\node[red] at (2.25,-0.8) {};
\node[red] at (2.25,0.8) {};
\node[red] at (1.45,2.1) {};
\node[red] at (1.45,-2.1) {};
\node[blue] at (2.25,0) {};
\node[blue] at (1.85,1.45) {};
\node[blue] at (1.85,-1.45) {};
\node[] at (0,3) [black]{$\Lk($\textcolor{purple}{purple}$)=$};
\node[] at (0,2) [black]{\small ($4m$\textit{-cycle})};

\draw [black, thick,fill=none,-] (4.75,0) node {} -- (4.75,-0.8) node {};
\draw [black, thick,fill=none,-] (4.75,0) node {} -- (4.75,0.8) node {};
\draw [black, thick,fill=none,-] (5.15,1.45) node {} -- (4.75,0.8) node {};
\draw [black, thick,fill=none,-] (5.15,1.45) node {} -- (5.55,2.1) node {};
\draw [black, thick,fill=none,-] (5.15,-1.45) node {} -- (4.75,-0.8) node {};
\draw [black, thick,fill=none,-] (5.15,-1.45) node {} -- (5.55,-2.1) node {};
\draw [black, thick,fill=none,-]  (5.55,2.1) node {} -- (5.95,2.3) node {};
\draw [black, thick,fill=none,-] (5.55,-2.1) node {} -- (5.95,-2.3) node {};
\node[circle,fill=black,text=black,inner sep=1pt] at (7,2.5) {};
\node[circle,fill=black,text=black,inner sep=1pt] at (7.25,2.48) {};
\node[circle,fill=black,text=black,inner sep=1pt] at (6.75,2.48) {};
\node[circle,fill=black,text=black,inner sep=1pt] at (7,-2.5) {};
\node[circle,fill=black,text=black,inner sep=1pt] at (7.25,-2.48) {};
\node[circle,fill=black,text=black,inner sep=1pt] at (6.75,-2.48) {};
\draw [black, thick,fill=none,-] (9.25,0) node {} -- (9.25,-0.8) node {};
\draw [black, thick,fill=none,-] (9.25,0) node {} -- (9.25,0.8) node {};
\draw [black, thick,fill=none,-] (8.85,1.45) node {} -- (9.25,0.8) node {};
\draw [black, thick,fill=none,-] (8.85,1.45) node {} -- (8.45,2.1) node {};
\draw [black, thick,fill=none,-] (8.85,-1.45) node {} -- (9.25,-0.8) node {};
\draw [black, thick,fill=none,-] (8.85,-1.45) node {} -- (8.45,-2.1) node {};
\draw [black, thick,fill=none,-]  (8.45,2.1) node {} -- (8.05,2.3) node {};
\draw [black, thick,fill=none,-] (8.45,-2.1) node {} -- (8.05,-2.3) node {};
\node[purple] at (4.75,-0.8) {};
\node[purple] at (4.75,0.8) {};
\node[purple] at (5.55,2.1) {};
\node[purple] at (5.55,-2.1) {};
\node[blue] at (4.75,0) {};
\node[blue] at (5.15,1.45) {};
\node[blue] at (5.15,-1.45) {};
\node[purple] at (9.25,-0.8) {};
\node[purple] at (9.25,0.8) {};
\node[purple] at (8.45,2.1) {};
\node[purple] at (8.45,-2.1) {};
\node[blue] at (9.25,0) {};
\node[blue] at (8.85,1.45) {};
\node[blue] at (8.85,-1.45) {};
\node[] at (7,3) [black]{$\Lk($\textcolor{red}{red}$)=$};
\node[] at (7,2) [black]{\small ($4x$\textit{-cycle})};

\draw [black, thick,fill=none,-] (1.5,-3.5) node {} -- (3,-5) node {};
\draw [black, thick,fill=none,-] (1.5,-3.5) node {} -- (0,-5) node {};
\draw [black, thick,fill=none,-] (1.5,-6.5) node {} -- (3,-5) node {};
\draw [black, thick,fill=none,-] (1.5,-6.5) node {} -- (0,-5) node {};
\node[purple] at (1.5,-3.5) {};
\node[purple] at (1.5,-6.5) {};
\node[red] at (3,-5) {};
\node[red] at (0,-5) {};
\node[] at (-1.5,-5) [black]{$\Lk($\textcolor{blue}{blue}$)=$};

\node[purple] at (8.5,-4.5) {};
\node[purple] at (8.5,-5.5) {};
\node[] at (6.5,-5) [black]{$\Lk($\textcolor{green}{green}$)=$};
\end{tikzpicture}
\caption{The set of isometry classes of links of $\Sigma$} \label{hrr}
\end{figure}
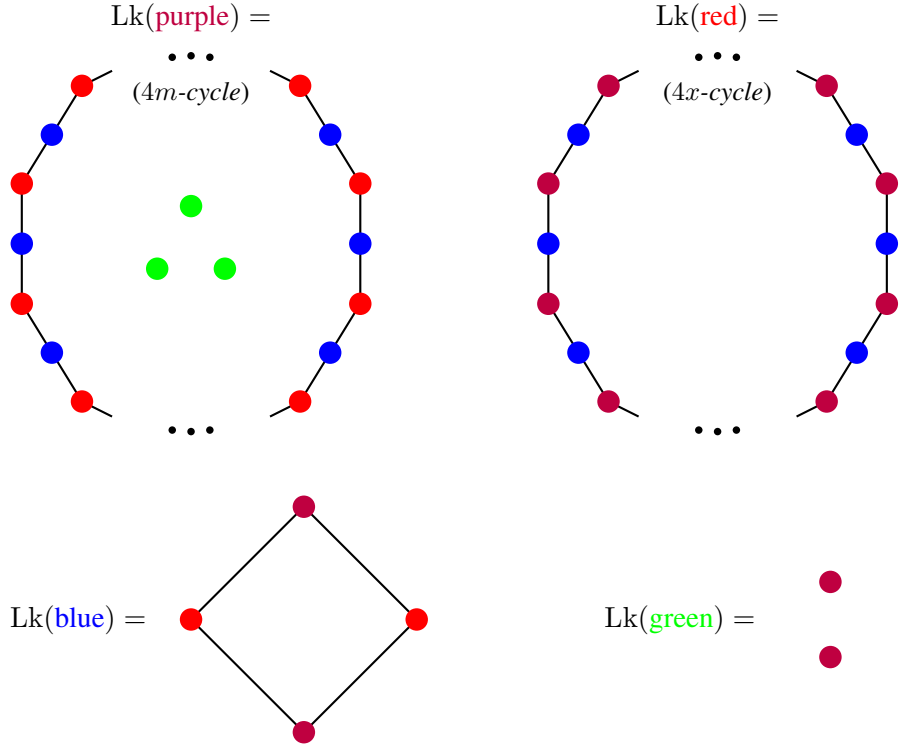

\begin{Corollary}\label{pfpf} Each complex of groups in the family $\{\mathcal X_k^{(0)}\}_{k \in \N} \cup \{ \mathcal X^{(0)}\}$ is strictly developable.
\begin{proof} Proposition~\ref{links} says that each complex of groups in the family $\{\mathcal X_k^{(0)}\}_{k \in \N} \cup \{ \mathcal X^{(0)}\}$ has a set of isometry classes of local developments identical to that of the Davis complex $\Sigma$. Hence they are $\CAT(0)$ under the restriction of Moussong's metric on $\Sigma$. The result then follows from Theorem~\ref{Haefliger}.
\end{proof} 
\end{Corollary}

\vspace{1mm}\noindent It is in fact not true, for an arbitrary polyhedral complex, that the local structure (the set of isometry classes of links) uniquely determines the global structure \cite{La}. This is an area of current research. So we need the following argument to compute the universal covers.

\begin{Proposition}\label{lplpl} Each complex of groups in the family $\{\mathcal X_k^{(0)}\}_{k \in \N} \cup \{ \mathcal X^{(0)}\}$ has a universal cover that is isomorphic to $\Sigma$ as a polyhedral complex. 
\begin{proof} It is known \cite{La} that, up to isomorphism, there exists a unique polygonal complex $\Xi$ such that the link of each vertex in $\Xi$ is the $2m$-cycle graph and each $2$-cell is a regular $2x$-gon for the integers $m, x \geq 2$  (again abusing notation and interpreting $\Xi$ as a complex of groups with all vertex groups trivial over itself as a scwol, its orientation irrelevant).

\vspace{2mm}\noindent Observe that the Davis complex $\Sigma_1$ of the Coxeter system $(W_1,\{s_4,...,s_n\})$ with nerve $L_1$ the $2m$-cycle graph with constant edge labels $x$ is isomorphic to the barycentric subdivision of $\Xi$ under the map sending the purple vertices of type $\varnothing$ to the pre-subdivision vertices of $\Xi$. So in this particular case the set of isometry classes of links does indeed uniquely determine the polyhedral complex $\Sigma_1$.

\vspace{2mm}\noindent We may think of $\Sigma$ as a "$(3,2)$-biregular tree of copies of $\Sigma_1$" glued at the purple vertices. If we wish to construct a simply connected polyhedral complex with the same set of isometry classes of links as $\Sigma$, we in fact have no choice in how to glue together the copies of $\Sigma_1$, for precisely the same reason that there exists a unique $(3,2)$-biregular tree. So $\Sigma$ is uniquely determined by its set of isometry classes of links, and hence the universal cover of each complex of groups in the family $\{\mathcal X_k^{(0)}\}_{k \in \N} \cup \{ \mathcal X^{(0)}\}$ is isomorphic to $\Sigma$ as a polyhedral complex. 
\end{proof} 
\end{Proposition}

\vspace{2mm}\noindent Having satisfied all the conditions required to apply Corollary~\ref{lllllllllll} to any $\mathcal X_k^{(0)}$, we deduce that $\pi_1 (\mathcal X_k^{(0)})$ is a uniform lattice in $\Aut(\Sigma)$ with covolume \begin{align*} 
\mu &\big(\pi_1 (\mathcal X_k^{(0)}) \backslash \Aut(\Sigma) \big) = \sum\limits_{v \in V(\mathcal X_k^{(0)})} \frac{1}{|(\mathcal X_{k})_{v}^{(0)}|}   \\
&=  \bigg(1+\frac{2x}{2}+\frac{4x}{2m}\bigg)\bigg(\sum\limits_{i=0}^{k} \frac{1}{2^i}\bigg)  + \frac{x}{m} - \frac{x}{2^k \cdot (2m)}\\
&\longrightarrow 2+2x+5x/m \hspace{2mm} \textnormal{ as } \hspace{2mm} k \longrightarrow \infty \end{align*} with a rate of convergence of $\frac{1}{2}$.

\vspace{3mm}\noindent Observe that the limit point found above may be realised as the covolume of the uniform lattice $\pi_1 (\mathcal X^{(0)})$ in $\Aut(\Sigma)$; see Figure~\ref{ggggg'}.

\vspace{4mm}\noindent \textbf{\underline{Case. {$ \boldsymbol2$}}:} $l>0$

\vspace{2mm}\noindent Recall that the nerve $L$ contains strictly more than four connected components in this case. We construct the relevant family of complexes of groups inductively, beginning with the family $\{\mathcal X_k^{(0)}\}_{k \in \N} \cup \{\mathcal X^{(0)}\}$ in Figures~\ref{ggggg} and ~\ref{ggggg'}. Let $L^{(0)}$ and $\Sigma^{(0)}$ be the nerve and Davis complex respectively of $(W,S)$ in the $l = 0$ case. 

\vspace{2mm}\noindent Consider the subgroup $W_1$ of $(W,S)$ with corresponding nerve $L_1$ satisfying either Condition~\ref{(i)} or ~\ref{(ii)} below. Denote the Davis complex of $W_1$ by $\Sigma_1$. As previously, we abuse notation and interpret $\Sigma_1$ as a complex of groups with all vertex groups trivial over itself as a scwol. The \textit{opposite complex of groups} $\Sigma_1^{op}$ is obtained by reversing the direction of all the arrows in $\Sigma_1$. Assign the colour purple to all vertices of type $\varnothing$ in $\Sigma_1$ (and $\Sigma_1^{op}$). Recall from Definition~\ref{thingy} the colour scheme of the complex of groups $\mathcal P^{m_1,x_1}$.

\vspace{2mm}\noindent Let $L^{(1)}:=L^{(0)} \sqcup L_1$ have corresponding Davis complex $\Sigma^{(1)}$, where $\sqcup$ denotes disjoint union. We show how to construct a family $\{\mathcal X_k^{(1)}\}_{k \in \N}$ of complexes of groups with sums of reciprocals of the orders of the vertex groups converging to that of a complex of groups $\mathcal X^{(1)}$, all with universal cover isomorphic to $\Sigma^{(1)}$ as a polyhedral complex.

\begin{Condition}\label{(i)} $W_1$ is a spherical subgroup of $(W,S)$.
\end{Condition}

\noindent Let $c_1$ be the lowest common multiple of the integers $|W_1|$ and $2x$. For each $i \in \{0,1,...,k\}$ (resp. $j\in \{0,1\}$) let $A_i$ (resp. $A_j$) be the set of all $c_1$ purple vertices on the $i$'th (resp. $j$'th) platforms of $c_1/(2x)$ copies of $\mathcal X_k^{(0)}$ (resp. $\mathcal X^{(0)}$), let $B_i$ (resp. $B_j$) be the set of all $c_1$ purple vertices of $c_1/ |W_1|$ copies of $\mathcal C_2^i \times \Sigma_1^{op}$ (resp. $\Sigma_1^{op}$) and let $\phi_i :A_i \to B_i$ (resp. $\phi_j :A_j \to B_j$) be any choice of bijection. Fix an integer $k$. We construct a complex of groups $\mathcal X_k^{(1)}$ (resp. $\mathcal X^{(1)}$) by first taking the disjoint union of $c_1/(2x)$ copies of $\mathcal X_k^{(0)}$ (resp. $\mathcal X^{(0)}$) with $c_1/ |W_1|$ copies of $\sqcup_{i \in \{0,1,...,k\}} \mathcal C_2^i \times \Sigma_1^{op}$ (resp. $\sqcup_{j\in \{0,1\}}\Sigma_1^{op}$), and then identifying pairs of purple vertices under each bijection in the family $\big\{\phi_i \mid i \in \{0,1,...,k\}\big\}$ (resp. $\big\{\phi_j \mid j\in \{0,1\}\big\}$).

\begin{Condition}\label{(ii)} $L_1$ is the $2m_1$-cycle graph with edges labelled by a constant integer $x_1 \geq 2$. 
\end{Condition}

\noindent Let $c_1$ be the lowest common multiple of the integers $2x_1$ and $2x$. For each $i \in \{0,1,...,k\}$ (resp. $j\in \{0,1\}$) recall that $A_i$ (resp. $A_j$) is the set of all $c_1$ purple vertices on the $i$'th (resp. $j$'th) platforms of $c_1/(2x)$ copies of $\mathcal X_k^{(0)}$ (resp. $\mathcal X^{(0)}$), let $C_i$ (resp. $C_j$) be the set of all $c_1$ purple vertices of $c_1/(2x_1)$ copies of $\mathcal P_i^{m_1,x_1} = \mathcal C_2^i \times \mathcal P^{m_1,x_1}$ (resp. $\mathcal P_j^{m_1,x_1} = \mathcal P^{m_1,x_1}$) and let $\psi_i :A_i \to C_i$ (resp. $\psi_j :A_j \to C_j$) be any choice of bijection. Fix an integer $k$. We construct a complex of groups $\mathcal X_k^{(1)}$ (resp. $\mathcal X^{(1)}$) by first taking the disjoint union of $c_1/(2x)$ copies of $\mathcal X_k^{(0)}$ (resp. $\mathcal X^{(0)}$) with $c_1/(2x_1)$ copies of $\sqcup_{i \in \{0,1,...,k\}} \mathcal P_i^{m_1,x_1}$ (resp. $\sqcup_{j\in \{0,1\}}\mathcal P_j^{m_1,x_1}$), and then identifying pairs of purple vertices under each bijection in the family $\big\{\psi_i \mid i \in \{0,1,...,k\}\big\}$ (resp. $\big\{\psi_j \mid j\in \{0,1\}\big\}$).

\begin{Proposition}\label{links2} Each complex of groups in the family $\{\mathcal X_k^{(1)}\}_{k \in \N} \cup \{ \mathcal X^{(1)}\}$ has a set of isometry classes of links identical to that of the Davis complex $\Sigma^{(1)}$.

\begin{proof} We use a similar argument as in Proposition~\ref{links}. Let $\mathcal F$ be any complex of groups in the family $\{\mathcal X_k^{(0)}\}_{k \in \N} \cup \{ \mathcal X^{(0)}\}$ and let $\mathcal G$ be any complex of groups in the family $\{\mathcal X_k^{(1)}\}_{k \in \N} \cup \{ \mathcal X^{(1)}\}$.

\vspace{2mm}\noindent Assume Condition~\ref{(i)} (resp. Condition~\ref{(ii)}) holds. It is not difficult to see that the links of the green, blue and red vertices in $\mathcal G$ are isometric to those of the corresponding coloured vertices in $\mathcal F$, and it follows from Case IV of Example~\ref{helpies} that the link of any purple vertex in $\mathcal G$ is the disjoint union of the link of any purple vertex in $\mathcal F$ with the link of any purple vertex in $\Sigma_1$ (resp. $\mathcal P^{m_1,x_1}$). Finally, the links of the remaining vertices in $\mathcal G$ are isometric to the links of the corresponding (non-specified but non-purple) coloured vertices in $\Sigma_1$ (resp. $\mathcal P^{m_1,x_1}$).
\end{proof}
\end{Proposition}

\begin{Corollary}\label{pfpfp} Each complex of groups in the family $\{\mathcal X_k^{(1)}\}_{k \in \N} \cup \{ \mathcal X^{(1)}\}$ is strictly developable.
\begin{proof} The result follows from Proposition~\ref{links2} and Theorem~\ref{Haefliger}, since each link of $\Sigma^{(1)}$ is $\CAT(0)$ under Moussong's metric.
\end{proof} 
\end{Corollary}

\begin{Proposition}\label{lplplp} Each complex of groups in the family $\{\mathcal X_k^{(1)}\}_{k \in \N} \cup \{ \mathcal X^{(1)}\}$ has a universal cover that is isomorphic to $\Sigma^{(1)}$ as a polyhedral complex. 
\begin{proof} Recall from Definition~\ref{thingy} that the complex of groups $\mathcal P^{m_1,x_1}$ is strictly developable. Interpret the universal cover $\Sigma^{m_1,x_1}$ of $\mathcal P^{m_1,x_1}$ as a complex of groups with all vertex groups trivial over itself as a scwol.

\vspace{2mm}\noindent Recall from Proposition~\ref{lplpl} that $\Sigma^{(0)}$ is uniquely determined by its set of isometry classes of links. In fact, this is also true of $\Sigma^{m_1,x_1}$ and the finite Davis complex $\Sigma_1$. For $\Sigma^{m_1,x_1}$ it follows from the proof of Proposition~\ref{lplpl} and for $\Sigma_1$ it follows from the fact that any finite Coxeter group naturally acts on a sphere of appropriate dimension, inducing a tesselation by spherical simplicies, with corresponding Davis complex the cone on this tesselation.

\vspace{2mm}\noindent Assume that Condition~\ref{(i)} (resp. Condition~\ref{(ii)}) holds. Observe that $\Sigma^{(1)}$ can be thought of as a "$(\Sigma^{(0)},\Sigma_1)$-biregular tree" (resp. "$(\Sigma^{(0)},\Sigma^{m_1,x_1})$-biregular tree") glued together at the purple vertices. By a similar argument as in the proof of Proposition~\ref{lplpl}, in the construction of $\Sigma^{(1)}$ we have no choice in how to glue together the copies of $\Sigma^{(0)}$ and $\Sigma_1$ (resp. $\Sigma^{m_1,x_1}$) since $\Sigma^{(1)}$ is simply connected.
\end{proof}
\end{Proposition}

\noindent We complete the proof of Theorem~\ref{coolthm} by iterating the above construction for the sequence of subgroups $(W_1,W_2, ..., W_l)$ of $(W,S)$ in order to obtain a family $\{\mathcal X_k:=\mathcal X_k^{(l)}\}_{k \in \N}$ of complexes of groups with sums of reciprocals of the orders of the vertex groups converging to that of a complex of groups $\mathcal X:=\mathcal X^{(l)}$, all with universal cover isomorphic to $\Sigma=\Sigma^{(l)}$ as a polyhedral complex. The construction of $\mathcal X_k^{(r)}$ (resp. $\mathcal X^{(r)}$) from $\mathcal X_k^{(r-1)}$ (resp. $\mathcal X^{(r-1)}$) for any integer $1 < r \leq  l$ is the same as the construction of $\mathcal X_k^{(1)}$ (resp. $\mathcal X^{(1)}$) from $\mathcal X_k^{(0)}$ (resp. $\mathcal X^{(0)}$) except that we take $c_r$ to be the lowest common multiple of the integers $|W_{r}|$ and $c_{r-1}$ if Condition~\ref{(i)} holds, or alternatively the lowest common multiple of the integers $2x_r$ and $c_{r-1}$ if Condition~\ref{(ii)} holds.
\end{proof}
\end{Theorem}

\begin{Remark}\label{importantremark} In fact Theorem~\ref{coolthm} still holds even if we relax the assumption that there exists at least one connected component of the nerve that is a cycle graph with an even number of sides and constant edge labels. To see this we modify the proof of Theorem~\ref{coolthm} as follows. Let $(\ast_{i=1,...,l+4} W_i,S)$ be as in the statement of Theorem~\ref{coolthm}, except that now we take $W_{l+1}$ to be finite. The case where each $W_i$ is isomorphic to $\mathcal C_2$ is Proposition~\ref{|||}. Interpret the Davis complex $\Sigma_{l+1}$ corresponding to the Coxeter group $W_{l+1}$ as a complex of groups with all vertex groups trivial over itself as a scwol. Assign the colour purple to all vertices of type $\varnothing$ of $\Sigma_{l+1}$ (and $\Sigma_{l+1}^{op}$). For any $k \in \N$, we construct a complex of groups $\mathcal Y_k^{(0)}$ (resp. $\mathcal Y^{(0)}$) from $\mathcal X_k^{(0)}$ (resp. $\mathcal X^{(0)}$) by replacing each embedded platform $\mathcal C_2^i \times \mathcal P^{m,x}$ with the complex of groups $\mathcal C_2^{i+1} \times \Sigma_{l+1}^{op}$ for every $i \in \{0,1,...,k\}$. As in Figure~\ref{ggggg} (resp. Figure~\ref{ggggg'}) we glue together the sets of purple vertices of any two adjacent platforms via green vertices in a $1-1$ correspondence. For every $i \in \{0,1,...,k\}$ we replace each green vertex group $\mathcal C_2^i \times \mathcal D_{2m}$ with $\mathcal C_2^{i+1}$ and replace each green vertex group $\mathcal C_{m}$ with the trivial group. It is then not hard to follow the arguments in the proof of Theorem~\ref{coolthm} and deduce that each complex of groups in the family $\{\mathcal Y_k^{(0)}\}_{k \in \N} \cup \{ \mathcal Y^{(0)}\}$ has a universal cover that is isomorphic to $\Sigma$ in the case where $l=0$. The remainder of the proof is identical to that of Theorem~\ref{coolthm}.
\end{Remark}

\vspace{1mm}\noindent By modifying the family $\{\mathcal X_k^{(0)}\}_{k \in \N}$ (resp. $\{\mathcal Y_k^{(0)}\}_{k \in \N}$) of complexes of groups described in the proof of Theorem~\ref{coolthm} (resp. Remark~\ref{importantremark}), we arrive at the following result.

\noindent\begin{Corollary}\label{|||||} Let $(W,S)$ be a Coxeter system with $j \geq 3$ free generators as in the statement of Theorem~\ref{coolthm} (resp. Remark~\ref{importantremark}). Let $\Sigma$ be the Davis complex of $(W,S)$ and let $p < j$ be prime. For any (arbitrarily large) $\alpha \in \N$ there exists a uniform lattice $\Gamma$ in $\Aut(\Sigma)$ with covolume $\mu \big(\Gamma \backslash \Aut(\Sigma) \big) = a/b$ (in lowest terms) such that $b$ is divisible by $p^{\alpha}$.

\begin{proof} Take any $k \in \N$. Construct a complex of groups ${\mathcal X'}_k^{(0)}$ from ${\mathcal X}_k^{(0)}$ by replacing the $i$'th platform in $\mathcal X_k^{(0)}$ with $\mathcal C_p^i \times \mathcal P^{m,x}$  for each $i \in \{0,1,...,k\}$ and adjusting each green vertex group immediately above platform $i$ from $\mathcal C_2^i \times \mathcal D_{2 m}$ to $\mathcal C_p^i \times \mathcal D_{2 m}$ accordingly. We also replace each green source vertex with local group $\mathcal C_m$ below platform $0$ with $p$ green source vertices each with local group $\mathcal D_{2m}$ and each adjacent to the same pair of purple vertices. All else is left the same. Similarly we construct a complex of groups ${\mathcal Y'}_k^{(0)}$ from ${\mathcal Y}_k^{(0)}$ by replacing each embedded copy of $\mathcal C_2$ with $\mathcal C_p$.

\vspace{2mm}\noindent We then construct a complex of groups ${\mathcal X'}_k$ from ${\mathcal X'}_k^{(0)}$ (resp. ${\mathcal Y'}_k$ from ${\mathcal Y'}_k^{(0)}$) by exactly the same iterative process as we obtained $\mathcal X_k$ from $\mathcal X_k^{(0)}$ in the proof of Theorem~\ref{coolthm} (with the only difference being that there are now $l-p+2$ iterations rather than $l$). 

\vspace{2mm}\noindent The fundamental group of ${\mathcal X'}_k$ (resp. ${\mathcal Y'}_k$) is a uniform lattice in $\Sigma$ by the same argument as in the proof of Theorem~\ref{coolthm}. A simple calculation will show that the corresponding sequence of covolumes converges at a rate of $\frac{1}{p}$ as $k \to \infty$. The sum of reciprocals of orders of the vertex groups (in lowest terms) will therefore have a denominator that is divisible by ever-increasing powers of $p$ as $k$ gets larger.
\end{proof}
\end{Corollary}

\vspace{2mm}\noindent We conclude this section with some observations and open questions.

\vspace{2mm}\noindent We could define an infinite complex of finite groups $\mathcal X_{\infty}$ by taking the obvious limit of the family $\{\mathcal X_k\}_{k \in \N}$ as $k \to \infty$. This would give us a family $\{\pi_1(\mathcal X_k)\}_{k \in \N}$ of uniform lattices in $\Aut(\Sigma)$ with a set of covolumes converging to that of the (necessarily) non-uniform lattice $\pi_1(\mathcal X_{\infty})$. In fact, this is what Thomas proved for a different subclass of Davis complexes in \cite{T2}. 

\vspace{2mm}\noindent In the proof of Theorem~\ref{coolthm}, observe that the limit point \[ \lim_{k \to \infty} \mu \big(\pi_1 (\mathcal X_k) \backslash \Aut(\Sigma) \big) =  \mu \big(\pi_1 (\mathcal X_{\infty}) \backslash \Aut(\Sigma)  \big)\] is a rational number. One wonders if it would be possible to find a similarly convergent family of complexes of groups $\{\mathcal Z_k\}_{k \in \N} \cup \{\mathcal Z_{\infty}\}$, sharing all the same properties of $\{\mathcal X_k\}_{k \in \N} \cup \{\mathcal X_{\infty}\}$, except that the limit point $\mu \big(\pi_1 (\mathcal Z_{\infty}) \backslash \Aut(\Sigma) \big)$ is irrational? Non-uniform lattices with irrational covolumes have been studied by Farb and Hruska \cite{FH} in the case of trees.

\vspace{2mm}\noindent Now recall the following definition.

\begin{Definition} An \textit{end} of a graph is an equivalence class of embedded rays under the equivalence relation that two rays $r_1$ and $r_2$ are equivalent iff there exists a ray $r_3$ that intersects with each of $r_1$ and $r_2$ at infinitely many vertices. An \textit{end} of a finitely generated group with a given generating set is an end of its Cayley graph.
\end{Definition} 

\noindent It is well-known (Theorem $8.32$ of \cite{BH}) that the number of ends of an infinite finitely generated group is a quasi-isometric invariant, hence independent of the choice of generating set, which can take only the values $1$, $2$ or $\infty$.

\vspace{2mm}\noindent Theorem $8.7.4.$ of \cite{D} characterises the number of ends of a Coxeter group. Precisely, an infinite Coxeter system $(W,S)$ has exactly $1$ end iff the \textit{punctured nerve} $L\smallsetminus \sigma_{T}$ is connected for each spherical subset $T$ of $S$, where $\sigma_{T}$ denotes the simplex associated to $T$ in the nerve. $(W,S)$ has exactly $2$ ends iff there exists a decomposition $(W,S)=(W_0 \times W_1,S_0 \cup S_1)$ where $W_0$ is the infinite dihedral group and $W_1$ is a finite Coxeter subgroup of $W$. 

\vspace{2mm}\noindent Now let $(W,S)$ be a Coxeter system as in either Proposition~\ref{|||} or Theorem~\ref{coolthm}. Recall that $S$ has $j \geq 3$ free generators $\{s_1,...,s_j\}$. Assume that we can decompose $(W,S)=(W_0 \times W_1,S_0 \cup S_1)$, where $W_1$ is finite and $W_0$ is (without loss of generality) the infinite dihedral group $\langle s_1,s_2 \rangle$. Then $s_3$ must commute with both $s_1$ and $s_2$ in $(W,S)$, a contradiction. Moreover observe that the punctured nerve $L \smallsetminus \sigma_{\{s_1\}}$ is disconnected. Hence $(W,S)$ has infinitely many ends. This raises the question:
\noindent\textit{does there exist a Coxeter system with $1$ end that has a non-discrete set of covolumes of uniform lattices acting on its Davis complex?}

\subsection{$\mathbf{\Sigma}$ is a regular right-angled building}\label{result3}

\noindent In this section we prove case $(ii)$ of Theorem~\ref{centraltheorem} and Corollary~\ref{centralcorollary}. Recall that $S:=\{s_1,...,s_n\}$.

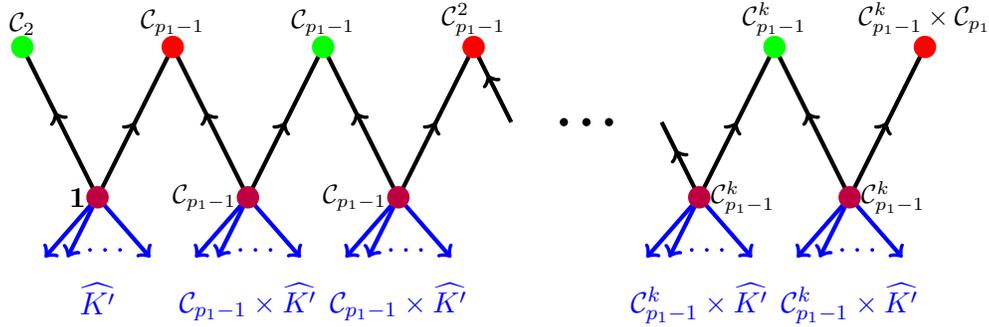
\begin{figure} [h]\
\centering
\usetikzlibrary{arrows,positioning}
\begin{tikzpicture} [%
    nd/.style = {circle,fill=red,text=black,inner sep=3pt},
    tn/.style = {circle,fill=blue,text=black,inner sep=3pt},
    aa/.style = {circle,fill=green,text=black,inner sep=3pt},
    purple/.style = {circle,fill=purple,text=black,inner sep=3pt}]

 \draw [black, ultra thick,fill=none,-<-] (-6,0) node[above, black] {\small $\mathcal C_2$} -- (-5,-2) node [left, black] {$\boldsymbol{1}$};

\draw [blue, ultra thick,fill=none,->] (-5,-2) node {} -- (-5.7,-2.8) node {};
\draw [blue, ultra thick,fill=none,->] (-5,-2) node {} -- (-5.4,-2.8) node {};
\draw [blue, ultra thick,fill=none,->] (-5,-2) node {} -- (-4.3,-2.8) node {};
\node[circle,fill=blue,text=black,inner sep=0.5pt] at (-5.1,-2.7) {};
\node[circle,fill=blue,text=black,inner sep=0.5pt] at (-4.9,-2.7) {};
\node[circle,fill=blue,text=black,inner sep=0.5pt] at (-4.7,-2.7) {};
\node[] at (-5,-3) [below, blue]{\large $\widehat{K'}$};

\draw [black, ultra thick,fill=none,->-] (-5,-2) node{} -- (-4,0) node [above, black] {\small $\mathcal C_{p_1-1}$};
\draw [black, ultra thick,fill=none,-<-] (-4,0) node{} -- (-3,-2) node [left, black] {\small $\mathcal C_{p_1-1}$};

\draw [blue, ultra thick,fill=none,->] (-3,-2) node {} -- (-3.7,-2.8) node {};
\draw [blue, ultra thick,fill=none,->] (-3,-2) node {} -- (-3.4,-2.8) node {};
\draw [blue, ultra thick,fill=none,->] (-3,-2) node {} -- (-2.3,-2.8) node {};
\node[circle,fill=blue,text=black,inner sep=0.5pt] at (-3.1,-2.7) {};
\node[circle,fill=blue,text=black,inner sep=0.5pt] at (-2.9,-2.7) {};
\node[circle,fill=blue,text=black,inner sep=0.5pt] at (-2.7,-2.7) {};
\node[] at (-3,-3) [below, blue]{$\mathcal C_{p_1-1} \times \widehat{K'}$};

\draw [black, ultra thick,fill=none,->-] (-3,-2) node{} -- (-2,0) node [above, black] {\small $\mathcal C_{p_1-1}$};
\draw [black, ultra thick,fill=none,-<-] (-2,0) node{} -- (-1,-2) node [left, black] {\small $\mathcal C_{p_1-1}$};

\draw [blue, ultra thick,fill=none,->] (-1,-2) node {} -- (-1.7,-2.8) node {};
\draw [blue, ultra thick,fill=none,->] (-1,-2) node {} -- (-1.4,-2.8) node {};
\draw [blue, ultra thick,fill=none,->] (-1,-2) node {} -- (-0.3,-2.8) node {};
\node[circle,fill=blue,text=black,inner sep=0.5pt] at (-1.1,-2.7) {};
\node[circle,fill=blue,text=black,inner sep=0.5pt] at (-0.9,-2.7) {};
\node[circle,fill=blue,text=black,inner sep=0.5pt] at (-0.7,-2.7) {};
\node[] at (-1,-3) [below, blue]{$\mathcal C_{p_1-1} \times \widehat{K'}$};

\draw [black, ultra thick,fill=none,->-] (-1,-2) node{} -- (0,0) node [above, black] {\small $\mathcal C_{p_1-1}^2$};

\draw [black, ultra thick,fill=none,->-] (3,-2) node[right, black] {\small $\mathcal C_{p_1-1}^{k}$} -- (4,0) node [above, black] {\small $\mathcal C_{p_1-1}^{k}$};

\draw [blue, ultra thick,fill=none,->] (3,-2) node {} -- (2.3,-2.8) node {};
\draw [blue, ultra thick,fill=none,->] (3,-2) node {} -- (2.6,-2.8) node {};
\draw [blue, ultra thick,fill=none,->] (3,-2) node {} -- (3.7,-2.8) node {};
\node[circle,fill=blue,text=black,inner sep=0.5pt] at (2.9,-2.7) {};
\node[circle,fill=blue,text=black,inner sep=0.5pt] at (3.1,-2.7) {};
\node[circle,fill=blue,text=black,inner sep=0.5pt] at (3.3,-2.7) {};
\node[] at (3,-3) [below, blue]{$\mathcal C_{p_1-1}^{k} \times \widehat{K'}$};

\draw [black, ultra thick,fill=none,-<-] (4,0) node{} -- (5,-2) node [right, black] {\small $\mathcal C_{p_1-1}^{k}$};

\draw [blue, ultra thick,fill=none,->] (5,-2) node {} -- (4.3,-2.8) node {};
\draw [blue, ultra thick,fill=none,->] (5,-2) node {} -- (4.6,-2.8) node {};
\draw [blue, ultra thick,fill=none,->] (5,-2) node {} -- (5.7,-2.8) node {};
\node[circle,fill=blue,text=black,inner sep=0.5pt] at (4.9,-2.7) {};
\node[circle,fill=blue,text=black,inner sep=0.5pt] at (5.1,-2.7) {};
\node[circle,fill=blue,text=black,inner sep=0.5pt] at (5.3,-2.7) {};
\node[] at (5,-3) [below, blue]{$\mathcal C_{p_1-1}^k \times \widehat{K'}$};

\draw [black, ultra thick,fill=none,->-] (5,-2) node{} -- (6,0) node [above, black] {\small $\mathcal C_{p_1-1}^{k} \times \mathcal C_{p_1}$};

\draw [black, ultra thick,fill=none,-<-] (0,0) node{} -- (0.5,-1) node {};
\draw [black, ultra thick,fill=none,-<-] (2.5,-1) node{} -- (3,-2) node {};

\node[aa] at (-6,0) {};
\node[purple] at (-5,-2) {};
\node[nd] at (-4,0) {};
\node[purple] at (-3,-2) {};
\node[aa] at (-2,0) {};
\node[purple] at (-1,-2) {};
\node[nd] at (0,0) {};

\node[circle,fill=black,text=black,inner sep=1.1pt] at (1.2,-1) {};
\node[circle,fill=black,text=black,inner sep=1.1pt] at (1.5,-1) {};
\node[circle,fill=black,text=black,inner sep=1.1pt] at (1.8,-1) {};

\node[purple] at (3,-2) {};
\node[aa] at (4,0) {};
\node[purple] at (5,-2) {};
\node[nd] at (6,0) {};

\end{tikzpicture}
\caption{A family $\{\mathcal A_k\}_{k \in \N}$ of finite complexes of finite groups} (constructed in the proof of Corollary~\ref{bldg} as $\{\mathcal A^{p_1-1}_k\}_{k \in \N}$ when $p_2=2$)\label{ggg}
\end{figure}

\begin{Theorem}\label{buildingthm} Let $(W,S)$ be a right-angled Coxeter system with (at least) two free generators, which we call $s_1$ and $s_2$. Let $\Sigma$ be the unique regular right-angled building of type $(W,S)$ with $p_i \geq 2$ chambers in each $s_i$-equivalence class, for all $i\in \{1,...,n\}$. Then the set of covolumes of uniform lattices in $\Aut(\Sigma)$ is non-discrete if either $p_1$ or $p_2$ is strictly greater than $2$.

\begin{proof} The building $\Sigma$ has a geometric realisation which is a union of embedded copies of the Davis complex of $(W,S)$, and which can be equipped with a piecewise Euclidean metric such that it is a complete $\CAT(0)$ space (Theorems $18.3.1$ and $18.3.9$ of \cite{D}). 

\vspace{2mm}\noindent Without loss of generality assume $p_1 \geq p_2$.

\vspace{2mm}\noindent The structure of the proof is as follows. We consider two cases, $p_1 > p_2 =2$ and $p_1 \geq p_2 >2$. If $p_1 > p_2 =2$ (resp. $p_1 \geq p_2 >2$) we construct a family $\{\mathcal A_k\}_{k \in \N}  \cup \{ \mathcal A\}$ (resp. $\{\mathcal B_k\}_{k \in \N}  \cup \{ \mathcal B\}$) of complexes of groups. For any $\mathcal J$ in either family, we show that the set of isometry classes of links of $\mathcal J$ is the same as that of $\Sigma$ (Proposition~\ref{links3}), that $\mathcal J$ is strictly developable (Corollary~\ref{pfpfpf}) and that the universal cover of $\mathcal J$ is $\Sigma$ (Proposition~\ref{lplplpl}). We then apply Corollary~\ref{lllllllllll} to show that $\pi_1(\mathcal J)$ is a uniform lattice in $\Aut(\Sigma)$ and compute the corresponding set of covolumes to establish the result.

\vspace{2mm}\noindent As previously, we abuse notation and interpret $\Sigma$ as a complex of groups with all vertex groups trivial over itself as a scwol. We colour $\Sigma$ by assigning the colour red to all vertices of type $s_1$, green to all vertices of type $s_2$, purple to all vertices of type $\varnothing$ and blue to all remaining vertices.

\vspace{2mm}\noindent Assume $n>2$. Let $S'$ be the subset $\{s_3,...,s_n\}$ of $S$. Let $K'$ denote the chamber of $(W_{S'},S')$ and let $\Sigma'$ be the unique regular right-angled building of type $(W_{S'},S')$ with $p_i \geq 2$ chambers in each $s_i$-equivalence class, for all $i \in \{3,...,n\}$. Recall that the vertices of $K'$ are labelled by the spherical subsets of $S'$. Construct a complex of groups $\widehat{K'}$ over the scwol $K'$ by associating to each vertex $\{s_{i_1},...,s_{i_N}\}$ in $K'$ the group $\mathcal C _{p_{i_1}} \times ... \times \mathcal C _{p_{i_N}}$ with natural inclusions as monomorphisms. One can check that $\widehat{K'}$ is strictly developable, with universal cover isomorphic to $\Sigma'$ as a polyhedral complex (refer to Section $2.3$ of \cite{T2}). If $n=2$ we take $\widehat{K'}$ to be the complex of groups consisting only of one vertex with trivial local group. We colour $\widehat{K'}$ completely blue except for its unique vertex of type $\varnothing$ which we colour purple. This is consistent with the colouring of $\Sigma$.

\begin{figure} [h!]\
\centering
\usetikzlibrary{arrows,positioning}
\begin{tikzpicture} [%
    nd/.style = {circle,fill=red,text=black,inner sep=3pt},
    tn/.style = {circle,fill=blue,text=black,inner sep=3pt},
    aa/.style = {circle,fill=green,text=black,inner sep=3pt},
    purple/.style = {circle,fill=purple,text=black,inner sep=3pt}]

 \draw [black, ultra thick,fill=none,-<-] (-7,0) node[above, black] {\small $\mathcal C_2 \times C_{p_1-2}$} -- (-6,-2) node [left, black] {\small $\mathcal C_{p_1-2}$};

\draw [blue, ultra thick,fill=none,->] (-6,-2) node {} -- (-6.7,-2.8) node {};
\draw [blue, ultra thick,fill=none,->] (-6,-2) node {} -- (-6.4,-2.8) node {};
\draw [blue, ultra thick,fill=none,->] (-6,-2) node {} -- (-5.3,-2.8) node {};
\node[circle,fill=blue,text=black,inner sep=0.5pt] at (-6.1,-2.7) {};
\node[circle,fill=blue,text=black,inner sep=0.5pt] at (-5.9,-2.7) {};
\node[circle,fill=blue,text=black,inner sep=0.5pt] at (-5.7,-2.7) {};
\node[] at (-6,-3) [below, blue]{$\mathcal C_{p_1-2} \times \widehat{K'}$};

\draw [black, ultra thick,fill=none,->-] (-6,-2) node{} -- (-3,1) node [above, black] {\small $\mathcal C_{p_1-2}$};
\draw [black, ultra thick,fill=none,-<-] (-3,1) node{} -- (-3,-2) node [left, black] {\small $\boldsymbol{1}$};

\draw [blue, ultra thick,fill=none,->] (-3,-2) node {} -- (-3.7,-2.8) node {};
\draw [blue, ultra thick,fill=none,->] (-3,-2) node {} -- (-3.4,-2.8) node {};
\draw [blue, ultra thick,fill=none,->] (-3,-2) node {} -- (-2.3,-2.8) node {};
\node[circle,fill=blue,text=black,inner sep=0.5pt] at (-3.1,-2.7) {};
\node[circle,fill=blue,text=black,inner sep=0.5pt] at (-2.9,-2.7) {};
\node[circle,fill=blue,text=black,inner sep=0.5pt] at (-2.7,-2.7) {};
\node[] at (-3,-3) [below, blue]{\large $\widehat{K'}$};

\draw [black, ultra thick,fill=none,->-] (-3,-2) node{} -- (-2,-1) node [above, black] {\small $\mathcal C_{2}$};
\draw [black, ultra thick,fill=none,-<-] (-3,1) node{} -- (0,-2) node [right, black] {\small $\mathcal C_{p_1-2}$};

\draw [blue, ultra thick,fill=none,->] (0,-2) node {} -- (-0.7,-2.8) node {};
\draw [blue, ultra thick,fill=none,->] (0,-2) node {} -- (-0.4,-2.8) node {};
\draw [blue, ultra thick,fill=none,->] (0,-2) node {} -- (0.7,-2.8) node {};
\node[circle,fill=blue,text=black,inner sep=0.5pt] at (-0.1,-2.7) {};
\node[circle,fill=blue,text=black,inner sep=0.5pt] at (0.1,-2.7) {};
\node[circle,fill=blue,text=black,inner sep=0.5pt] at (0.3,-2.7) {};
\node[] at (0,-3) [below, blue]{$\mathcal C_{p_1-2} \times \widehat{K'}$};

\draw [black, ultra thick,fill=none,->-] (0,-2) node{} -- (1,0) node [above, black] {\small $\mathcal C_2 \times \mathcal C_{p_1-2}$};

\node[aa] at (-7,0) {};
\node[purple] at (-6,-2) {};
\node[nd] at (-3,1) {};
\node[purple] at (-3,-2) {};
\node[aa] at (-2,-1) {};
\node[purple] at (0,-2) {};
\node[aa] at (1,0) {};

\end{tikzpicture}
\caption{A finite complex of finite groups $\mathcal A$} \label{ggg'}
\end{figure}
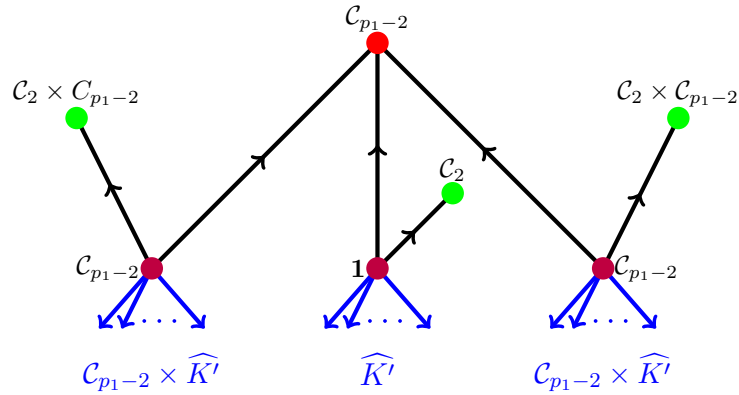

\vspace{2mm}\noindent In Figures~\ref{ggg}, ~\ref{ggg'}, ~\ref{gggg'} and ~\ref{gggg} we construct a family $\cup_{k \in \N} \{\mathcal A_k, \mathcal B_k\} \cup \{ \mathcal A, \mathcal B\}$ of complexes of groups that each include multiple embedded copies of $\mathcal C_{p-1}^i \times \widehat{K'}$ for some $i$ in $\{0,1,...,k\}$. We colour such embedded copies of "multiples" of $\widehat{K'}$ accordingly, that is, all vertices are blue other than those of type $\varnothing$ which we colour purple. We also use the colour red (resp. green) to denote vertices of type $s_1$ (resp. $s_2$), that is, with pre-image in the universal cover of degree $p_1$ (resp. $p_2$).

\vspace{2mm}\noindent A quick remark on how to interpret Figure~\ref{gggg}. For any integer $k$, the complex of groups $\mathcal B_k$ consists of "multiples" of copies of an embedded "block" that are stacked below one another, glued together via additional green vertices with local groups that are even powers of $\mathcal C_{p_2-1}$. Each embedded block is a $\mathcal C^2_{p_2-1}$-multiple of the block immediately above it. We draw two such blocks in Figure~\ref{gggg} - the "base" case and the $\mathcal C^{2k}_{p_2-1}$-multiple.

\vspace{2mm}\noindent Let $\Omega$ denote the sum of reciprocals of orders of the vertex groups in $\widehat{K'}$.

\begin{Proposition}\label{links3} If $p_1 > p_2 =2$ (resp. $p_1 \geq p_2 >2$) then each complex of groups in the family $\{\mathcal A_k\}_{k \in \N}  \cup \{ \mathcal A\}$ (resp. $\{\mathcal B_k\}_{k \in \N}  \cup \{ \mathcal B\}$) has a set of isometry classes of links identical to that of $\Sigma$.

\begin{proof} We use a similar argument as in Proposition~\ref{links}. Let $\mathcal J$ be any complex of groups in the family $\cup_{k \in \N} \{\mathcal A_k, \mathcal B_k\} \cup \{ \mathcal A, \mathcal B\} $.

\vspace{2mm}\noindent It follows from Cases III and IV of Example~\ref{helpies} that the link of each red (resp. green) vertex in $\mathcal J$ is $p_1$ (resp. $p_2$) disconnected purple vertices and that the link of each purple vertex in $\mathcal J$ is the disjoint union of the link of the unique purple vertex of type $\varnothing$ in $\widehat{K'}$ with one lone red vertex and one lone green vertex. The link of each of the remaining (blue) vertices in $\mathcal J$ is isometric to the link of any vertex of the same type in $\widehat{K'}$.

\vspace{2mm}\noindent It is clear from the construction of the building $\Sigma$, and the fact that $\widehat{K'}$ has universal cover $\Sigma'$, that the link of each vertex in $\mathcal J$ is isometric to the link of any vertex of the same type in $\Sigma$.
\end{proof}
\end{Proposition}

\begin{Corollary}\label{pfpfpf}Each complex of groups in the family $\cup_{k \in \N} \{\mathcal A_k, \mathcal B_k\} \cup \{ \mathcal A, \mathcal B\} $ is strictly developable.
\begin{proof} The result follows from Proposition~\ref{links3} and Theorem~\ref{Haefliger}, since each link of $\Sigma$ is $\CAT(0)$ under the restriction of the metric in Theorem $18.3.9$ of \cite{D}.
\end{proof} 
\end{Corollary}

\begin{Proposition}\label{lplplpl} If $p_1 > p_2 =2$ (resp. $p_1 \geq p_2 >2$) then each complex of groups in the family $\{\mathcal A_k\}_{k \in \N}  \cup \{ \mathcal A\}$ (resp. $\{\mathcal B_k\}_{k \in \N}  \cup \{ \mathcal B\}$) has a universal cover that is isomorphic to $\Sigma$ as a polyhedral complex.
\begin{proof} Recall that there exists a unique regular right-angled building $\Sigma$ such that for all $i\in \{1,...,n\}$ each $s_i$-equivalence class contains $p_i$ chambers (Theorem $8$ of \cite{T4}). Hence $\Sigma$ is uniquely determined by its set of isometry classes of links. The result then follows from Proposition~\ref{links3} and Corollary~\ref{pfpfpf}.
\end{proof}
\end{Proposition}

\noindent \textbf{\underline{Case. {$ \boldsymbol1$}}:} $p_1 > p_2 =2$

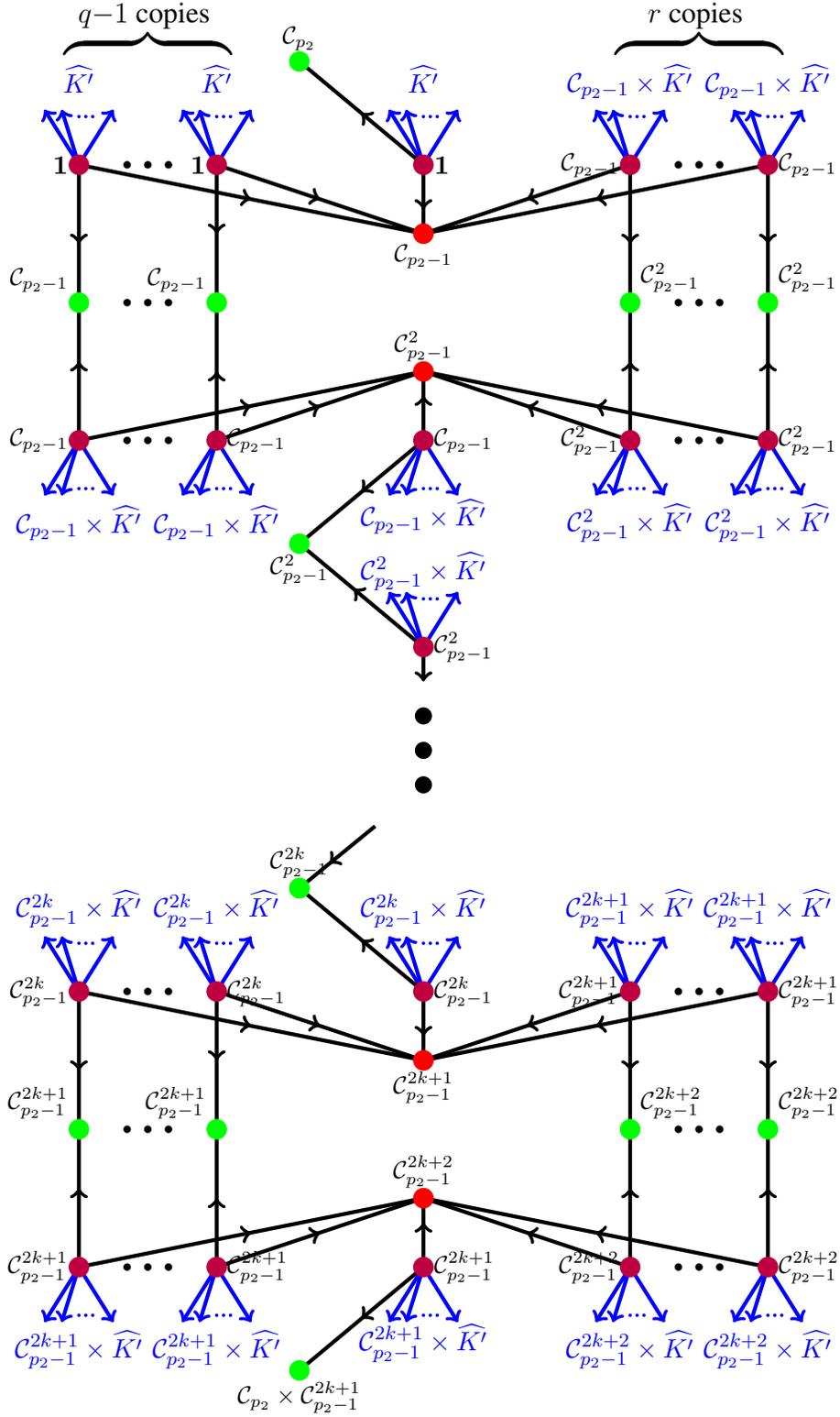
\begin{figure} \
\centering
\usetikzlibrary{arrows,positioning}
\begin{tikzpicture} [%
    nd/.style = {circle,fill=green,text=black,inner sep=3pt},
    tn/.style = {circle,fill=blue,text=black,inner sep=3pt},
    aa/.style = {circle,fill=red,text=black,inner sep=3pt},
    purple/.style = {circle,fill=purple,text=black,inner sep=3pt}]

\node[] at (7,3.2) [above, black]{\LARGE $\overbrace{\textnormal{\textcolor{white}{aaaaaaaaa}}}^{r \textnormal{ copies }}$};
\node[] at (-1,3.2) [above, black]{\LARGE $\overbrace{\textnormal{\textcolor{white}{aaaaaaaaa}}}^{q-1 \textnormal{ copies }}$};

\draw [black, ultra thick,fill=none,-<-] (-2,0) node [above left, black] {\small $\mathcal C_{p_2-1}$} -- (-2,-2) node [left, black] {\small $\mathcal C_{p_2-1}$};
\draw [blue, ultra thick,fill=none,->] (-2,-2) node {} -- (-2.5,-2.8) node {};
\draw [blue, ultra thick,fill=none,->] (-2,-2) node {} -- (-2.25,-2.8) node {};
\draw [blue, ultra thick,fill=none,->] (-2,-2) node {} -- (-1.5,-2.8) node {};
\node[circle,fill=blue,text=black,inner sep=0.5pt] at (-2,-2.7) {};
\node[circle,fill=blue,text=black,inner sep=0.5pt] at (-1.9,-2.7) {};
\node[circle,fill=blue,text=black,inner sep=0.5pt] at (-1.8,-2.7) {};
\node[circle,text=blue,inner sep=0.5pt] at (-2,-3.2) {$\mathcal C_{p_2-1} \times \widehat{K'}$};

\draw [black, ultra thick,fill=none,-<-] (-2,0) node{} -- (-2,2) node [left, black] {$\boldsymbol{1}$};
\draw [blue, ultra thick,fill=none,->] (-2,2) node {} -- (-2.5,2.8) node {};
\draw [blue, ultra thick,fill=none,->] (-2,2) node {} -- (-2.25,2.8) node {};
\draw [blue, ultra thick,fill=none,->] (-2,2) node {} -- (-1.5,2.8) node {};
\node[circle,fill=blue,text=black,inner sep=0.5pt] at (-2,2.7) {};
\node[circle,fill=blue,text=black,inner sep=0.5pt] at (-1.9,2.7) {};
\node[circle,fill=blue,text=black,inner sep=0.5pt] at (-1.8,2.7) {};
\node[circle,text=blue,inner sep=0.5pt] at (-2,3.25) {$ \widehat{K'}$};

\draw [black, ultra thick,fill=none,->-] (0,-2) node [right, black] {\small $\mathcal C_{p_2-1}$} -- (0,0) node [above left, black] {\small $\mathcal C_{p_2-1}$};
\draw [blue, ultra thick,fill=none,->] (0,-2) node {} -- (-0.5,-2.8) node {};
\draw [blue, ultra thick,fill=none,->] (0,-2) node {} -- (-0.25,-2.8) node {};
\draw [blue, ultra thick,fill=none,->] (0,-2) node {} -- (0.5,-2.8) node {};
\node[circle,fill=blue,text=black,inner sep=0.5pt] at (0,-2.7) {};
\node[circle,fill=blue,text=black,inner sep=0.5pt] at (0.1,-2.7) {};
\node[circle,fill=blue,text=black,inner sep=0.5pt] at (0.2,-2.7) {};
\node[circle,text=blue,inner sep=0.5pt] at (0,-3.2) {$\mathcal C_{p_2-1} \times \widehat{K'}$};

\draw [black, ultra thick,fill=none,->-] (0,2) node [left, black] {$\boldsymbol{1}$} -- (0,0) node {};
\draw [blue, ultra thick,fill=none,->] (0,2) node {} -- (-0.5,2.8) node {};
\draw [blue, ultra thick,fill=none,->] (0,2) node {} -- (-0.25,2.8) node {};
\draw [blue, ultra thick,fill=none,->] (0,2) node {} -- (0.5,2.8) node {};
\node[circle,fill=blue,text=black,inner sep=0.5pt] at (0,2.7) {};
\node[circle,fill=blue,text=black,inner sep=0.5pt] at (0.1,2.7) {};
\node[circle,fill=blue,text=black,inner sep=0.5pt] at (0.2,2.7) {};
\node[circle,text=blue,inner sep=0.5pt] at (0,3.25) {$ \widehat{K'}$};

\draw [black, ultra thick,fill=none,-<-] (6,0) node{} -- (6,-2) node [ left, black] {\small $\mathcal C_{p_2-1}^2$};
\draw [blue, ultra thick,fill=none,->] (6,-2) node {} -- (5.5,-2.8) node {};
\draw [blue, ultra thick,fill=none,->] (6,-2) node {} -- (5.75,-2.8) node {};
\draw [blue, ultra thick,fill=none,->] (6,-2) node {} -- (6.5,-2.8) node {};
\node[circle,fill=blue,text=black,inner sep=0.5pt] at (6,-2.7) {};
\node[circle,fill=blue,text=black,inner sep=0.5pt] at (6.1,-2.7) {};
\node[circle,fill=blue,text=black,inner sep=0.5pt] at (6.2,-2.7) {};
\node[circle,text=blue,inner sep=0.5pt] at (6,-3.2) {$\mathcal C_{p_2-1}^2 \times \widehat{K'}$};

\draw [black, ultra thick,fill=none,-<-] (6,0) node [above right, black] {\small $\mathcal C_{p_2-1}^2$} -- (6,2) node [left, black] {\small $\mathcal C_{p_2-1}$};
\draw [blue, ultra thick,fill=none,->] (6,2) node {} -- (5.5,2.8) node {};
\draw [blue, ultra thick,fill=none,->] (6,2) node {} -- (5.75,2.8) node {};
\draw [blue, ultra thick,fill=none,->] (6,2) node {} -- (6.5,2.8) node {};
\node[circle,fill=blue,text=black,inner sep=0.5pt] at (6,2.7) {};
\node[circle,fill=blue,text=black,inner sep=0.5pt] at (6.1,2.7) {};
\node[circle,fill=blue,text=black,inner sep=0.5pt] at (6.2,2.7) {};
\node[circle,text=blue,inner sep=0.5pt] at (6,3.2) {$\mathcal C_{p_2-1} \times \widehat{K'}$};

\draw [black, ultra thick,fill=none,-<-] (8,0) node [above right, black] {\small $\mathcal C_{p_2-1}^2$} -- (8,-2) node [right, black] {\small $\mathcal C_{p_2-1}^2$};
\draw [blue, ultra thick,fill=none,->] (8,-2) node {} -- (7.5,-2.8) node {};
\draw [blue, ultra thick,fill=none,->] (8,-2) node {} -- (7.75,-2.8) node {};
\draw [blue, ultra thick,fill=none,->] (8,-2) node {} -- (8.5,-2.8) node {};
\node[circle,fill=blue,text=black,inner sep=0.5pt] at (8,-2.7) {};
\node[circle,fill=blue,text=black,inner sep=0.5pt] at (8.1,-2.7) {};
\node[circle,fill=blue,text=black,inner sep=0.5pt] at (8.2,-2.7) {};
\node[circle,text=blue,inner sep=0.5pt] at (8,-3.2) {$\mathcal C_{p_2-1}^2 \times \widehat{K'}$};

\draw [black, ultra thick,fill=none,-<-] (8,0) node{} -- (8,2) node [right, black] {\small $\mathcal C_{p_2-1}$};
\draw [blue, ultra thick,fill=none,->] (8,2) node {} -- (7.5,2.8) node {};
\draw [blue, ultra thick,fill=none,->] (8,2) node {} -- (7.75,2.8) node {};
\draw [blue, ultra thick,fill=none,->] (8,2) node {} -- (8.5,2.8) node {};
\node[circle,fill=blue,text=black,inner sep=0.5pt] at (8,2.7) {};
\node[circle,fill=blue,text=black,inner sep=0.5pt] at (8.1,2.7) {};
\node[circle,fill=blue,text=black,inner sep=0.5pt] at (8.2,2.7) {};
\node[circle,text=blue,inner sep=0.5pt] at (8,3.2) {$\mathcal C_{p_2-1} \times \widehat{K'}$};

\draw [black, ultra thick,fill=none,->-] (-2,2) node{} -- (3,1) node [below, black] {\small $\mathcal C_{p_2-1}$};
\draw [black, ultra thick,fill=none,->-] (0,2) node{} -- (3,1) node {};
\draw [black, ultra thick,fill=none,->-] (6,2) node{} -- (3,1) node {};
\draw [black, ultra thick,fill=none,->-] (8,2) node{} -- (3,1) node {};
\draw [black, ultra thick,fill=none,->-] (-2,-2) node{} -- (3,-1) node [above, black] {\small $\mathcal C_{p_2-1}^2$};
\draw [black, ultra thick,fill=none,->-] (0,-2) node{} -- (3,-1) node {};
\draw [black, ultra thick,fill=none,->-] (6,-2) node{} -- (3,-1) node {};
\draw [black, ultra thick,fill=none,->-] (8,-2) node{} -- (3,-1) node{};

\draw [black, ultra thick,fill=none,-<-] (3,1) node{} -- (3,2) node [right, black] {$\boldsymbol{1}$};
\draw [blue, ultra thick,fill=none,->] (3,2) node {} -- (2.5,2.8) node {};
\draw [blue, ultra thick,fill=none,->] (3,2) node {} -- (2.75,2.8) node {};
\draw [blue, ultra thick,fill=none,->] (3,2) node {} -- (3.5,2.8) node {};
\node[circle,fill=blue,text=black,inner sep=0.5pt] at (3,2.7) {};
\node[circle,fill=blue,text=black,inner sep=0.5pt] at (3.1,2.7) {};
\node[circle,fill=blue,text=black,inner sep=0.5pt] at (3.2,2.7) {};
\node[circle,text=blue,inner sep=0.5pt] at (3,3.25) {$ \widehat{K'}$};
\draw [black, ultra thick,fill=none,->-] (3,2) node{} -- (1.2,3.5) node [above, black] {\small $\mathcal C_{p_2}$};

\draw [black, ultra thick,fill=none,-<-] (3,-1) node{} -- (3,-2) node [right, black] {\small $\mathcal C_{p_2-1}$};
\draw [blue, ultra thick,fill=none,->] (3,-2) node {} -- (2.5,-2.8) node {};
\draw [blue, ultra thick,fill=none,->] (3,-2) node {} -- (2.75,-2.8) node {};
\draw [blue, ultra thick,fill=none,->] (3,-2) node {} -- (3.5,-2.8) node {};
\node[circle,fill=blue,text=black,inner sep=0.5pt] at (3,-2.7) {};
\node[circle,fill=blue,text=black,inner sep=0.5pt] at (3.1,-2.7) {};
\node[circle,fill=blue,text=black,inner sep=0.5pt] at (3.2,-2.7) {};
\node[circle,text=blue,inner sep=0.5pt] at (3,-3.1) {$\mathcal C_{p_2-1} \times \widehat{K'}$};

\draw [black, ultra thick,fill=none,->-] (3,-2) node{} -- (1.2,-3.5) node [below, black] {\small $\mathcal C_{p_2-1}^2$};
\draw [black, ultra thick,fill=none,-<-] (1.2,-3.5) node{} -- (3,-5) node [right, black] {\small $\mathcal C_{p_2-1}^2$};
\draw [blue, ultra thick,fill=none,->] (3,-5) node {} -- (2.5,-4.2) node {};
\draw [blue, ultra thick,fill=none,->] (3,-5) node {} -- (2.75,-4.2) node {};
\draw [blue, ultra thick,fill=none,->] (3,-5) node {} -- (3.5,-4.2) node {};
\node[circle,fill=blue,text=black,inner sep=0.5pt] at (3,-4.3) {};
\node[circle,fill=blue,text=black,inner sep=0.5pt] at (3.1,-4.3) {};
\node[circle,fill=blue,text=black,inner sep=0.5pt] at (3.2,-4.3) {};
\node[circle,text=blue,inner sep=0.5pt] at (3,-3.9) {$\mathcal C_{p_2-1}^2 \times \widehat{K'}$};

\draw [black, ultra thick,fill=none,->] (3,-5) node{} -- (3,-5.5) node {};

\node[circle,fill=black,text=black,inner sep=2.5pt] at (3,-6) {};
\node[circle,fill=black,text=black,inner sep=2.5pt] at (3,-6.5) {};
\node[circle,fill=black,text=black,inner sep=2.5pt] at (3,-7) {};

\node[nd] at (0,0) {};
\node[nd] at (-2,0) {};
\node[purple] at (-2,-2) {};
\node[purple] at (-2,2) {};
\node[purple] at (0,-2) {};
\node[purple] at (0,2) {};
\node[nd] at (6,0) {};
\node[purple] at (6,-2) {};
\node[purple] at (6,2) {};
\node[nd] at (8,0) {};
\node[purple] at (8,-2) {};
\node[purple] at (8,2) {};
\node[aa] at (3,1) {};
\node[aa] at (3,-1) {};
\node[purple] at (3,2) {};
\node[purple] at (3,-2) {};
\node[nd] at (1.2,3.5) {};
\node[nd] at (1.2,-3.5) {};
\node[purple] at (3,-5) {};

\node[circle,fill=black,text=black,inner sep=1.1pt] at (-1.3,-2) {};
\node[circle,fill=black,text=black,inner sep=1.1pt] at (-1,-2) {};
\node[circle,fill=black,text=black,inner sep=1.1pt] at (-0.7,-2) {};
\node[circle,fill=black,text=black,inner sep=1.1pt] at (-1.3,2) {};
\node[circle,fill=black,text=black,inner sep=1.1pt] at (-1,2) {};
\node[circle,fill=black,text=black,inner sep=1.1pt] at (-0.7,2) {};

\node[circle,fill=black,text=black,inner sep=1.1pt] at (6.7,-2) {};
\node[circle,fill=black,text=black,inner sep=1.1pt] at (7,-2) {};
\node[circle,fill=black,text=black,inner sep=1.1pt] at (7.3,-2) {};
\node[circle,fill=black,text=black,inner sep=1.1pt] at (6.7,2) {};
\node[circle,fill=black,text=black,inner sep=1.1pt] at (7,2) {};
\node[circle,fill=black,text=black,inner sep=1.1pt] at (7.3,2) {};

\node[circle,fill=black,text=black,inner sep=1.1pt] at (-1.3,0) {};
\node[circle,fill=black,text=black,inner sep=1.1pt] at (-1,0) {};
\node[circle,fill=black,text=black,inner sep=1.1pt] at (-0.7,0) {};
\node[circle,fill=black,text=black,inner sep=1.1pt] at (6.7,0) {};
\node[circle,fill=black,text=black,inner sep=1.1pt] at (7,0) {};
\node[circle,fill=black,text=black,inner sep=1.1pt] at (7.3,0) {};

\draw [black, ultra thick,fill=none,-<-] (-2,-12) node[above left, black] {\small $\mathcal C_{p_2-1}^{2k+1}$} -- (-2,-14) node [left, black] {\small $\mathcal C_{p_2-1}^{2k+1}$};
\draw [blue, ultra thick,fill=none,->] (-2,-14) node {} -- (-2.5,-14.8) node {};
\draw [blue, ultra thick,fill=none,->] (-2,-14) node {} -- (-2.25,-14.8) node {};
\draw [blue, ultra thick,fill=none,->] (-2,-14) node {} -- (-1.5,-14.8) node {};
\node[circle,fill=blue,text=black,inner sep=0.5pt] at (-2,-14.7) {};
\node[circle,fill=blue,text=black,inner sep=0.5pt] at (-1.9,-14.7) {};
\node[circle,fill=blue,text=black,inner sep=0.5pt] at (-1.8,-14.7) {};
\node[circle,text=blue,inner sep=0.5pt] at (-2,-15.2) {$\mathcal C_{p_2-1}^{2k+1} \times \widehat{K'}$};

\draw [black, ultra thick,fill=none,-<-] (-2,-12) node{} -- (-2,-10) node [left, black] {\small $\mathcal C_{p_2-1}^{2k}$};
\draw [blue, ultra thick,fill=none,->] (-2,-10) node {} -- (-2.5,-9.2) node {};
\draw [blue, ultra thick,fill=none,->] (-2,-10) node {} -- (-2.25,-9.2) node {};
\draw [blue, ultra thick,fill=none,->] (-2,-10) node {} -- (-1.5,-9.2) node {};
\node[circle,fill=blue,text=black,inner sep=0.5pt] at (-2,-9.3) {};
\node[circle,fill=blue,text=black,inner sep=0.5pt] at (-1.9,-9.3) {};
\node[circle,fill=blue,text=black,inner sep=0.5pt] at (-1.8,-9.3) {};
\node[circle,text=blue,inner sep=0.5pt] at (-2,-8.8) {$\mathcal C_{p_2-1}^{2k} \times \widehat{K'}$};

\draw [black, ultra thick,fill=none,->-] (0,-14) node [right, black] {\small $\mathcal C_{p_2-1}^{2k+1}$} -- (0,-12) node [above left, black] {\small $\mathcal C_{p_2-1}^{2k+1}$};
\draw [blue, ultra thick,fill=none,->] (0,-14) node {} -- (-0.5,-14.8) node {};
\draw [blue, ultra thick,fill=none,->] (0,-14) node {} -- (-0.25,-14.8) node {};
\draw [blue, ultra thick,fill=none,->] (0,-14) node {} -- (0.5,-14.8) node {};
\node[circle,fill=blue,text=black,inner sep=0.5pt] at (0,-14.7) {};
\node[circle,fill=blue,text=black,inner sep=0.5pt] at (0.1,-14.7) {};
\node[circle,fill=blue,text=black,inner sep=0.5pt] at (0.2,-14.7) {};
\node[circle,text=blue,inner sep=0.5pt] at (0,-15.2) {$\mathcal C_{p_2-1}^{2k+1} \times \widehat{K'}$};

\draw [black, ultra thick,fill=none,->-] (0,-10) node [right, black] {\small $\mathcal C_{p_2-1}^{2k}$} -- (0,-12) node {};
\draw [blue, ultra thick,fill=none,->] (0,-10) node {} -- (-0.5,-9.2) node {};
\draw [blue, ultra thick,fill=none,->] (0,-10) node {} -- (-0.25,-9.2) node {};
\draw [blue, ultra thick,fill=none,->] (0,-10) node {} -- (0.5,-9.2) node {};
\node[circle,fill=blue,text=black,inner sep=0.5pt] at (0,-9.3) {};
\node[circle,fill=blue,text=black,inner sep=0.5pt] at (0.1,-9.3) {};
\node[circle,fill=blue,text=black,inner sep=0.5pt] at (0.2,-9.3) {};
\node[circle,text=blue,inner sep=0.5pt] at (0,-8.8) {$\mathcal C_{p_2-1}^{2k} \times \widehat{K'}$};

\draw [black, ultra thick,fill=none,-<-] (6,-12) node{} -- (6,-14) node [ left, black] {\small $\mathcal C_{p_2-1}^{2k+2}$};
\draw [blue, ultra thick,fill=none,->] (6,-14) node {} -- (5.5,-14.8) node {};
\draw [blue, ultra thick,fill=none,->] (6,-14) node {} -- (5.75,-14.8) node {};
\draw [blue, ultra thick,fill=none,->] (6,-14) node {} -- (6.5,-14.8) node {};
\node[circle,fill=blue,text=black,inner sep=0.5pt] at (6,-14.7) {};
\node[circle,fill=blue,text=black,inner sep=0.5pt] at (6.1,-14.7) {};
\node[circle,fill=blue,text=black,inner sep=0.5pt] at (6.2,-14.7) {};
\node[circle,text=blue,inner sep=0.5pt] at (6,-15.2) {$\mathcal C_{p_2-1}^{2k+2} \times \widehat{K'}$};

\draw [black, ultra thick,fill=none,-<-] (6,-12) node [above right, black] {\small $\mathcal C_{p_2-1}^{2k+2}$} -- (6,-10) node [left, black] {\small $\mathcal C_{p_2-1}^{2k+1}$};
\draw [blue, ultra thick,fill=none,->] (6,-10) node {} -- (5.5,-9.2) node {};
\draw [blue, ultra thick,fill=none,->] (6,-10) node {} -- (5.75,-9.2) node {};
\draw [blue, ultra thick,fill=none,->] (6,-10) node {} -- (6.5,-9.2) node {};
\node[circle,fill=blue,text=black,inner sep=0.5pt] at (6,-9.3) {};
\node[circle,fill=blue,text=black,inner sep=0.5pt] at (6.1,-9.3) {};
\node[circle,fill=blue,text=black,inner sep=0.5pt] at (6.2,-9.3) {};
\node[circle,text=blue,inner sep=0.5pt] at (6,-8.8) {$\mathcal C_{p_2-1}^{2k+1} \times \widehat{K'}$};

\draw [black, ultra thick,fill=none,-<-] (8,-12) node [above right, black] {\small $\mathcal C_{p_2-1}^{2k+2}$} -- (8,-14) node [right, black] {\small $\mathcal C_{p_2-1}^{2k+2}$};
\draw [blue, ultra thick,fill=none,->] (8,-14) node {} -- (7.5,-14.8) node {};
\draw [blue, ultra thick,fill=none,->] (8,-14) node {} -- (7.75,-14.8) node {};
\draw [blue, ultra thick,fill=none,->] (8,-14) node {} -- (8.5,-14.8) node {};
\node[circle,fill=blue,text=black,inner sep=0.5pt] at (8,-14.7) {};
\node[circle,fill=blue,text=black,inner sep=0.5pt] at (8.1,-14.7) {};
\node[circle,fill=blue,text=black,inner sep=0.5pt] at (8.2,-14.7) {};
\node[circle,text=blue,inner sep=0.5pt] at (8,-15.2) {$\mathcal C_{p_2-1}^{2k+2} \times \widehat{K'}$};

\draw [black, ultra thick,fill=none,-<-] (8,-12) node{} -- (8,-10) node [right, black] {\small $\mathcal C_{p_2-1}^{2k+1}$};
\draw [blue, ultra thick,fill=none,->] (8,-10) node {} -- (7.5,-9.2) node {};
\draw [blue, ultra thick,fill=none,->] (8,-10) node {} -- (7.75,-9.2) node {};
\draw [blue, ultra thick,fill=none,->] (8,-10) node {} -- (8.5,-9.2) node {};
\node[circle,fill=blue,text=black,inner sep=0.5pt] at (8,-9.3) {};
\node[circle,fill=blue,text=black,inner sep=0.5pt] at (8.1,-9.3) {};
\node[circle,fill=blue,text=black,inner sep=0.5pt] at (8.2,-9.3) {};
\node[circle,text=blue,inner sep=0.5pt] at (8,-8.8) {$\mathcal C_{p_2-1}^{2k+1} \times \widehat{K'}$};

\draw [black, ultra thick,fill=none,->-] (-2,-10) node{} -- (3,-11) node [below, black] {\small $\mathcal C_{p_2-1}^{2k+1}$};
\draw [black, ultra thick,fill=none,->-] (0,-10) node{} -- (3,-11) node {};
\draw [black, ultra thick,fill=none,->-] (6,-10) node{} -- (3,-11) node {};
\draw [black, ultra thick,fill=none,->-] (8,-10) node{} -- (3,-11) node {};
\draw [black, ultra thick,fill=none,->-] (-2,-14) node{} -- (3,-13) node [above, black] {\small $\mathcal C_{p_2-1}^{2k+2}$};
\draw [black, ultra thick,fill=none,->-] (0,-14) node{} -- (3,-13) node {};
\draw [black, ultra thick,fill=none,->-] (6,-14) node{} -- (3,-13) node {};
\draw [black, ultra thick,fill=none,->-] (8,-14) node{} -- (3,-13) node{};

\draw [black, ultra thick,fill=none,-<-] (3,-11) node{} -- (3,-10) node [right, black] {\small $\mathcal C_{p_2-1}^{2k}$};
\draw [blue, ultra thick,fill=none,->] (3,-10) node {} -- (2.5,-9.2) node {};
\draw [blue, ultra thick,fill=none,->] (3,-10) node {} -- (2.75,-9.2) node {};
\draw [blue, ultra thick,fill=none,->] (3,-10) node {} -- (3.5,-9.2) node {};
\node[circle,fill=blue,text=black,inner sep=0.5pt] at (3,-9.3) {};
\node[circle,fill=blue,text=black,inner sep=0.5pt] at (3.1,-9.3) {};
\node[circle,fill=blue,text=black,inner sep=0.5pt] at (3.2,-9.3) {};
\node[circle,text=blue,inner sep=0.5pt] at (3,-8.8) {$\mathcal C_{p_2-1}^{2k} \times \widehat{K'}$};
\draw [black, ultra thick,fill=none,->-] (3,-10) node{} -- (1.2,-8.5) node {};

\draw [black, ultra thick,fill=none,-<-] (3,-13) node{} -- (3,-14) node [right, black] {\small $\mathcal C_{p_2-1}^{2k+1}$};
\draw [blue, ultra thick,fill=none,->] (3,-14) node {} -- (2.5,-14.8) node {};
\draw [blue, ultra thick,fill=none,->] (3,-14) node {} -- (2.75,-14.8) node {};
\draw [blue, ultra thick,fill=none,->] (3,-14) node {} -- (3.5,-14.8) node {};
\node[circle,fill=blue,text=black,inner sep=0.5pt] at (3,-14.7) {};
\node[circle,fill=blue,text=black,inner sep=0.5pt] at (3.1,-14.7) {};
\node[circle,fill=blue,text=black,inner sep=0.5pt] at (3.2,-14.7) {};
\node[circle,text=blue,inner sep=0.5pt] at (3,-15.1) {$\mathcal C_{p_2-1}^{2k+1} \times \widehat{K'}$};

\draw [black, ultra thick,fill=none,->-] (3,-14) node{} -- (1.2,-15.5) node [below, black] {\small $\mathcal C_{p_2} \times \mathcal C_{p_2-1}^{2k+1}$};

\node[nd] at (0,-12) {};
\node[nd] at (-2,-12) {};
\node[purple] at (-2,-14) {};
\node[purple] at (-2,-10) {};
\node[purple] at (0,-14) {};
\node[purple] at (0,-10) {};
\node[nd] at (6,-12) {};
\node[purple] at (6,-14) {};
\node[purple] at (6,-10) {};
\node[nd] at (8,-12) {};
\node[purple] at (8,-14) {};
\node[purple] at (8,-10) {};
\node[aa] at (3,-11) {};
\node[aa] at (3,-13) {};
\node[purple] at (3,-10) {};
\node[purple] at (3,-14) {};
\node[nd] at (1.2,-8.5) {};
\node[nd] at (1.2,-15.5) {};

\node[circle,fill=black,text=black,inner sep=1.1pt] at (-1.3,-14) {};
\node[circle,fill=black,text=black,inner sep=1.1pt] at (-1,-14) {};
\node[circle,fill=black,text=black,inner sep=1.1pt] at (-0.7,-14) {};
\node[circle,fill=black,text=black,inner sep=1.1pt] at (-1.3,-10) {};
\node[circle,fill=black,text=black,inner sep=1.1pt] at (-1,-10) {};
\node[circle,fill=black,text=black,inner sep=1.1pt] at (-0.7,-10) {};

\node[circle,fill=black,text=black,inner sep=1.1pt] at (6.7,-14) {};
\node[circle,fill=black,text=black,inner sep=1.1pt] at (7,-14) {};
\node[circle,fill=black,text=black,inner sep=1.1pt] at (7.3,-14) {};
\node[circle,fill=black,text=black,inner sep=1.1pt] at (6.7,-10) {};
\node[circle,fill=black,text=black,inner sep=1.1pt] at (7,-10) {};
\node[circle,fill=black,text=black,inner sep=1.1pt] at (7.3,-10) {};

\node[circle,fill=black,text=black,inner sep=1.1pt] at (-1.3,-12) {};
\node[circle,fill=black,text=black,inner sep=1.1pt] at (-1,-12) {};
\node[circle,fill=black,text=black,inner sep=1.1pt] at (-0.7,-12) {};
\node[circle,fill=black,text=black,inner sep=1.1pt] at (6.7,-12) {};
\node[circle,fill=black,text=black,inner sep=1.1pt] at (7,-12) {};
\node[circle,fill=black,text=black,inner sep=1.1pt] at (7.3,-12) {};

\draw [black, ultra thick,fill=none,-<-] (1.2,-8.5) node [above, black] {\small $\mathcal C_{p_2-1}^{2k}$} -- (2.3,-7.6) node {};
\node[nd] at (1.2,-8.5) {};

\end{tikzpicture}
\caption{A family $\{\mathcal B_k\}_{k \in \N}$ of finite complexes of finite groups}\label{gggg}
\end{figure}

\vspace{1.5mm}\noindent Having satisfied all the conditions required to apply Corollary~\ref{lllllllllll} to $\mathcal A_k$ for any $k \in \N$, therefore $\pi_1 (\mathcal A_k)$ is a uniform lattice in $\Aut (\Sigma)$ with covolume \begin{align*} 
\mu &\big(\pi_1 (\mathcal A_k) \backslash \Aut(\Sigma) \big) = \sum\limits_{v \in V(\mathcal A_k)} \frac{1}{|({\mathcal A}_{k})_{v}|}   \\
&=  \bigg(1+ \sum\limits_{i=1}^{k} \frac{2}{(p_1-1)^i}\bigg) \Omega + \frac{1}{2}+ \sum\limits_{i=1}^{k} \frac{2}{(p_1-1)^i} +\frac{1}{p_1(p_1-1)^{k}} \\
&\longrightarrow \bigg(1+ \frac{2}{p_1-2}\bigg)\Omega+ \frac{1}{2}+ \frac{2}{p_1-2} \hspace{2mm} \textnormal{ as } \hspace{2mm} k \longrightarrow \infty \end{align*} with a rate of convergence of $\frac{1}{p_1-1}$.

\vspace{2mm}\noindent Observe that the limit point found above may be realised as the covolume of the uniform lattice $\pi_1 (\mathcal A)$ in $\Aut(\Sigma)$; see Figure~\ref{ggg'}.

\vspace{4mm}\noindent \textbf{\underline{Case. {$ \boldsymbol2$}}:} $p_1 \geq p_2 >2$

\vspace{2mm}\noindent Use the Euclidean algorithm to write $p_1=(p_2-1)q+r$ for some non-negative integers $q$ and $r$.

\vspace{2mm}\noindent Consider the family ${\{ \mathcal B_k\}}_{k \in \N} \cup \{\mathcal B\}$ of complexes of groups shown in Figures~\ref{gggg'} and ~\ref{gggg}. Having satisfied all the conditions required to apply Corollary~\ref{lllllllllll} to $\mathcal B_k$ for any $k \in \N$, we deduce that $\pi_1 (\mathcal B_k)$ is a uniform lattice in $\Aut (\Sigma)$ with covolume \begin{align*} 
\mu &\big(\pi_1 (\mathcal B_k) \backslash \Aut(\Sigma) \big) = \sum\limits_{v \in V(\mathcal B_k)} \frac{1}{|(\mathcal B_{k})_{v}|}   \\
&=  \bigg[\Big(q+ \frac{r+q}{p_2 -1} +\frac{r}{{(p_2 -1)}^2} \Big) \Omega + 1+\frac{q}{p_2-1}+\frac{r+1}{(p_2-1)^2} \bigg]\bigg(\sum\limits_{i=0}^{k} \frac{1}{(p_2-1)^{2i}} \bigg)\\
&\hspace{30mm}+\frac{1}{p_2}- 1 - \frac{1}{p_2 (p_2-1)^{2k+1}} \\
&\longrightarrow \frac{p_1(p_2-1+1)\Omega+ (p_2-1)^2 +p_1+1}{(p_2-1)^2-1}+\frac{1}{p_2}- 1 \hspace{2mm} \textnormal{ as } \hspace{2mm} k \longrightarrow \infty\\
&=  \frac{p_1p_2\Omega +p_1+p_2}{p_2(p_2-2)}
\end{align*} with a rate of convergence of $\frac{1}{p_2-1}$. 

\vspace{2mm}\noindent Observe that the limit point found above may be realised as the covolume of the uniform lattice $\pi_1 (\mathcal B)$ in $\Aut(\Sigma)$; see Figure~\ref{gggg'}.
\end{proof}
\end{Theorem}

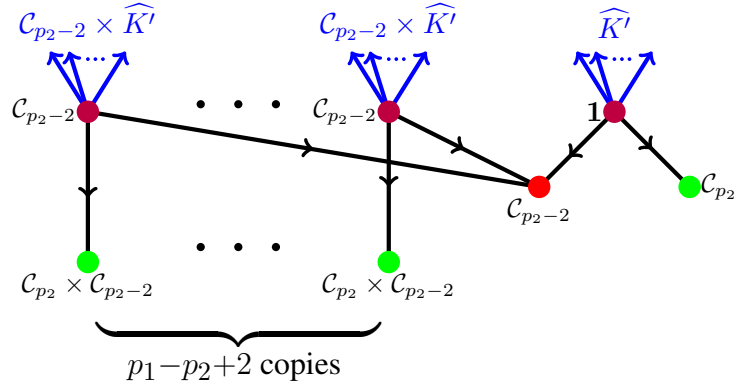
\begin{figure}[h] \
\centering
\usetikzlibrary{arrows,positioning}
\begin{tikzpicture} [%
    nd/.style = {circle,fill=green,text=black,inner sep=3pt},
    tn/.style = {circle,fill=blue,text=black,inner sep=3pt},
    aa/.style = {circle,fill=red,text=black,inner sep=3pt},
    purple/.style = {circle,fill=purple,text=black,inner sep=3pt}]

\node[] at (0,-0.4) [below, black]{\LARGE $\underbrace{\textnormal{\textcolor{white}{aaaaaaaaaaaaaa}}}_{p_1-p_2+2 \textnormal{ copies }}$};

\draw [black, ultra thick,fill=none,-<-] (-2,0) node [below, black] {\small $\mathcal C_{p_2} \times \mathcal C_{p_2-2}$} -- (-2,2) node [left, black] {\small $\mathcal C_{p_2-2}$};
\draw [blue, ultra thick,fill=none,->] (-2,2) node {} -- (-2.5,2.8) node {};
\draw [blue, ultra thick,fill=none,->] (-2,2) node {} -- (-2.25,2.8) node {};
\draw [blue, ultra thick,fill=none,->] (-2,2) node {} -- (-1.5,2.8) node {};
\node[circle,fill=blue,text=black,inner sep=0.5pt] at (-2,2.7) {};
\node[circle,fill=blue,text=black,inner sep=0.5pt] at (-1.9,2.7) {};
\node[circle,fill=blue,text=black,inner sep=0.5pt] at (-1.8,2.7) {};
\node[] at (-2,2.8) [above, blue]{$\mathcal C_{p_2-2} \times \widehat{K'}$};

\draw [black, ultra thick,fill=none,->-] (2,2) node [left, black] {\small $\mathcal C_{p_2-2}$} -- (2,0) node [below, black] {\small $\mathcal C_{p_2} \times \mathcal C_{p_2-2}$};
\draw [blue, ultra thick,fill=none,->] (2,2) node {} -- (1.5,2.8) node {};
\draw [blue, ultra thick,fill=none,->] (2,2) node {} -- (1.75,2.8) node {};
\draw [blue, ultra thick,fill=none,->] (2,2) node {} -- (2.5,2.8) node {};
\node[circle,fill=blue,text=black,inner sep=0.5pt] at (2,2.7) {};
\node[circle,fill=blue,text=black,inner sep=0.5pt] at (2.1,2.7) {};
\node[circle,fill=blue,text=black,inner sep=0.5pt] at (2.2,2.7) {};
\node[] at (2,2.8) [above, blue]{$\mathcal C_{p_2-2} \times \widehat{K'}$};

\draw [black, ultra thick,fill=none,->-] (-2,2) node{} -- (4,1) node [below, black] {\small $\mathcal C_{p_2-2}$};
\draw [black, ultra thick,fill=none,->-] (2,2) node{} -- (4,1) node {};

\draw [black, ultra thick,fill=none,-<-] (4,1) node{} -- (5,2) node [left, black] {$\boldsymbol{1}$};
\draw [blue, ultra thick,fill=none,->] (5,2) node {} -- (4.5,2.8) node {};
\draw [blue, ultra thick,fill=none,->] (5,2) node {} -- (4.75,2.8) node {};
\draw [blue, ultra thick,fill=none,->] (5,2) node {} -- (5.5,2.8) node {};
\node[circle,fill=blue,text=black,inner sep=0.5pt] at (5,2.7) {};
\node[circle,fill=blue,text=black,inner sep=0.5pt] at (5.1,2.7) {};
\node[circle,fill=blue,text=black,inner sep=0.5pt] at (5.2,2.7) {};
\node[] at (5,2.85) [above, blue]{$\widehat{K'}$};
\draw [black, ultra thick,fill=none,->-] (5,2) node{} -- (6,1) node [right, black] {\small $\mathcal C_{p_2}$};

\node[nd] at (2,0) {};
\node[nd] at (-2,0) {};
\node[purple] at (-2,2) {};
\node[purple] at (2,2) {};
\node[aa] at (4,1) {};
\node[purple] at (5,2) {};
\node[nd] at (6,1) {};

\node[circle,fill=black,text=black,inner sep=1.1pt] at (-0.5,2.1) {};
\node[circle,fill=black,text=black,inner sep=1.1pt] at (0,2.1) {};
\node[circle,fill=black,text=black,inner sep=1.1pt] at (0.5,2.1) {};
\node[circle,fill=black,text=black,inner sep=1.1pt] at (-0.5,0.2) {};
\node[circle,fill=black,text=black,inner sep=1.1pt] at (0,0.2) {};
\node[circle,fill=black,text=black,inner sep=1.1pt] at (0.5,0.2) {};

\end{tikzpicture}
\caption{A finite complex of finite groups $\mathcal B$} \label{gggg'}
\end{figure}

\noindent Generalising the family $\{\mathcal A_k\}_{k \in \N} $ of complexes of groups in Figure~\ref{ggg} gives us a new proof of the following (already known) result. Our argument is similar to the proof of Corollary~\ref{||||}. 

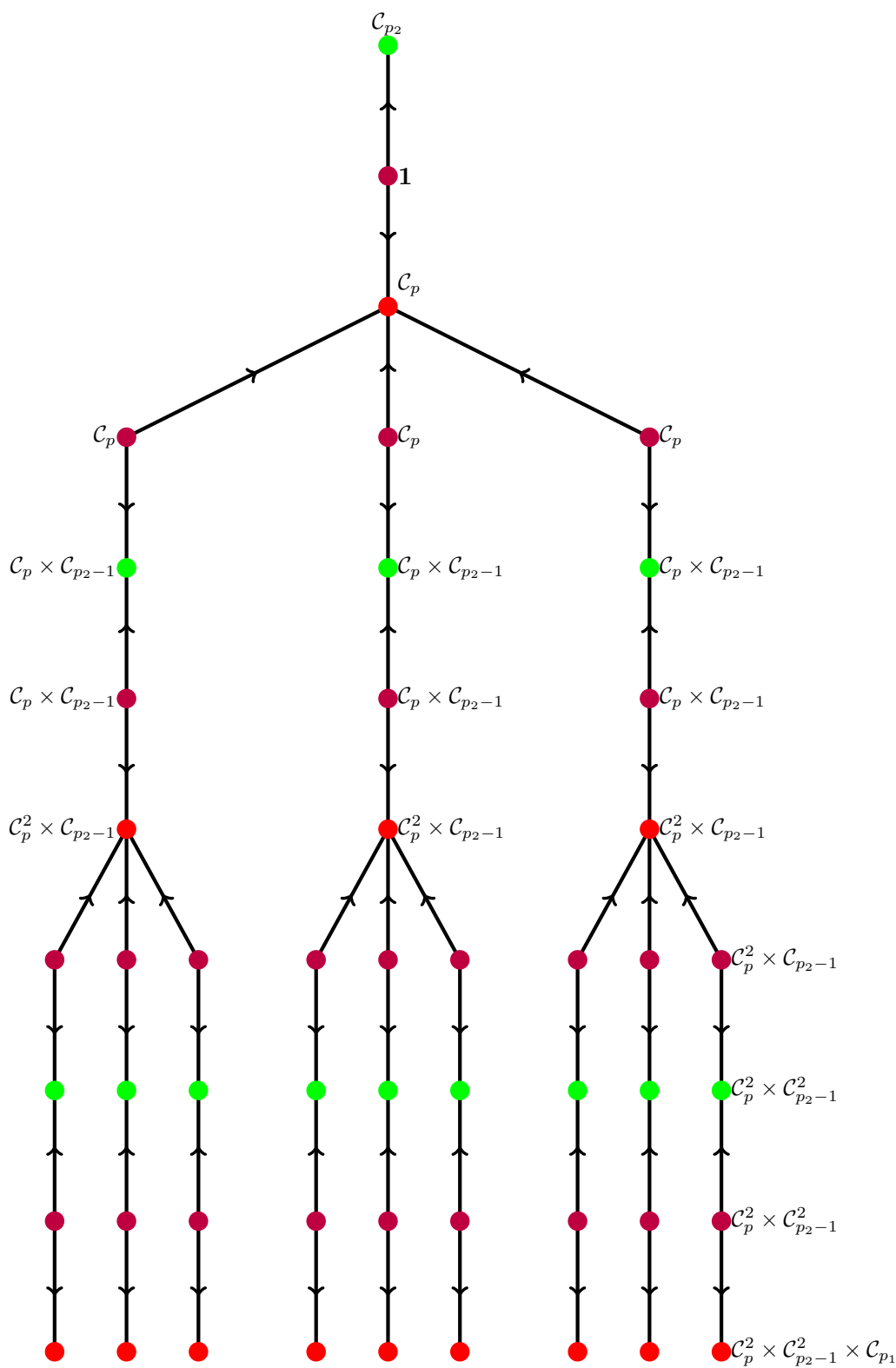
\begin{figure}\
\centering
\usetikzlibrary{arrows,positioning}
\begin{tikzpicture} [%
    red/.style = {circle,fill=red,text=black,inner sep=3pt},
    blue/.style = {circle,fill=blue,text=black,inner sep=3pt},
    green/.style = {circle,fill=green,text=black,inner sep=3pt},
   purple/.style = {circle,fill=purple,text=black,inner sep=3pt}]

\draw [black, ultra thick,fill=none,-<-] (0,10) node[above, black] {\small $\mathcal C_{p_2}$} -- (0,8) node {};
\draw [black, ultra thick,fill=none,->-] (0,8) node[right, black] {$\boldsymbol 1$} -- (0,6) node {};
 \draw [black, ultra thick,fill=none,-<-] (0,6) node[above right, black] {\small $\mathcal C_p $} -- (0,4) node [right, black] {\small $\mathcal C_p$};
 \draw [black, ultra thick,fill=none,->-] (4,4) node[right, black] {\small $\mathcal C_p$} -- (0,6) node {};
 \draw [black, ultra thick,fill=none,->-] (-4,4) node[left, black] {\small $\mathcal C_p$} -- (0,6) node {};
 \draw [black, ultra thick,fill=none,-<-] (4,2) node{} -- (4,0) node {};
 \draw [black, ultra thick,fill=none,-<-] (-4,2) node{} -- (-4,0) node {};
 \draw [black, ultra thick,fill=none,-<-] (0,2) node{} -- (0,0) node {};

 \draw [black, ultra thick,fill=none,-<-] (4,2) node[right, black] {\small $\mathcal C_p \times \mathcal C_{p_2-1}$} -- (4,4) node {};
 \draw [black, ultra thick,fill=none,-<-] (-4,2) node[left, black] {\small $\mathcal C_p\times \mathcal C_{p_2-1}$} -- (-4,4) node {};
 \draw [black, ultra thick,fill=none,-<-] (0,2) node[right, black] {\small $\mathcal C_p\times \mathcal C_{p_2-1}$} -- (0,4) node {};

 \draw [black, ultra thick,fill=none,-<-] (4,-2) node[right, black] {\small $\mathcal C_p^2 \times \mathcal C_{p_2-1}$} -- (4,0) node [right, black] {\small $\mathcal C_p\times \mathcal C_{p_2-1}$};
 \draw [black, ultra thick,fill=none,-<-] (-4,-2) node[left, black] {\small $\mathcal C_p^2\times \mathcal C_{p_2-1}$} -- (-4,0) node [left, black] {\small $\mathcal C_p\times \mathcal C_{p_2-1}$};
 \draw [black, ultra thick,fill=none,-<-] (0,-2) node[right, black] {\small $\mathcal C_p^2\times \mathcal C_{p_2-1}$} -- (0,0) node [right, black] {\small $\mathcal C_p\times \mathcal C_{p_2-1}$};
 \draw [black, ultra thick,fill=none,->-] (-2.9,-4) node{} -- (-4,-2) node {};
 \draw [black, ultra thick,fill=none,->-] (-4,-4) node{} -- (-4,-2) node {};
 \draw [black, ultra thick,fill=none,->-] (-5.1,-4) node{} -- (-4,-2) node {};
 \draw [black, ultra thick,fill=none,->-] (1.1,-4) node{} -- (0,-2) node {};
 \draw [black, ultra thick,fill=none,->-] (0,-4) node{} -- (0,-2) node {};
 \draw [black, ultra thick,fill=none,->-] (-1.1,-4) node{} -- (0,-2) node {};
 \draw [black, ultra thick,fill=none,->-] (2.9,-4) node{} -- (4,-2) node {};
 \draw [black, ultra thick,fill=none,->-] (4,-4) node{} -- (4,-2) node {};
 \draw [black, ultra thick,fill=none,->-] (5.1,-4) node{} -- (4,-2) node {};
 \draw [black, ultra thick,fill=none,-<-] (-2.9,-6) node{} -- (-2.9,-4) node {};
 \draw [black, ultra thick,fill=none,-<-] (-4,-6) node{} -- (-4,-4) node {};
 \draw [black, ultra thick,fill=none,-<-] (-5.1,-6) node{} -- (-5.1,-4) node {};
 \draw [black, ultra thick,fill=none,-<-] (1.1,-6) node{} -- (1.1,-4) node {};
 \draw [black, ultra thick,fill=none,-<-] (0,-6) node{} -- (0,-4) node {};
 \draw [black, ultra thick,fill=none,-<-] (-1.1,-6) node{} -- (-1.1,-4) node {};
 \draw [black, ultra thick,fill=none,-<-] (2.9,-6) node{} -- (2.9,-4) node {};
 \draw [black, ultra thick,fill=none,-<-] (4,-6) node{} -- (4,-4) node {};
 \draw [black, ultra thick,fill=none,-<-] (5.1,-6) node[right, black] {\small $\mathcal C_p^2\times \mathcal C_{p_2-1}^2$} -- (5.1,-4) node [right, black] {\small $\mathcal C_p^2\times \mathcal C_{p_2-1}$};

 \draw [black, ultra thick,fill=none,-<-] (-2.9,-6) node{} -- (-2.9,-8) node {};
 \draw [black, ultra thick,fill=none,-<-] (-4,-6) node{} -- (-4,-8) node {};
 \draw [black, ultra thick,fill=none,-<-] (-5.1,-6) node{} -- (-5.1,-8) node {};
 \draw [black, ultra thick,fill=none,-<-] (1.1,-6) node{} -- (1.1,-8) node {};
 \draw [black, ultra thick,fill=none,-<-] (0,-6) node{} -- (0,-8) node {};
 \draw [black, ultra thick,fill=none,-<-] (-1.1,-6) node{} -- (-1.1,-8) node {};
 \draw [black, ultra thick,fill=none,-<-] (2.9,-6) node{} -- (2.9,-8) node {};
 \draw [black, ultra thick,fill=none,-<-] (4,-6) node{} -- (4,-8) node {};
 \draw [black, ultra thick,fill=none,-<-] (5.1,-6) node{} -- (5.1,-8) node [right, black] {\small $\mathcal C_p^2\times \mathcal C_{p_2-1}^2$};

 \draw [black, ultra thick,fill=none,-<-] (-2.9,-10) node{} -- (-2.9,-8) node {};
 \draw [black, ultra thick,fill=none,-<-] (-4,-10) node{} -- (-4,-8) node{};
 \draw [black, ultra thick,fill=none,-<-] (-5.1,-10) node {} -- (-5.1,-8) node {};
 \draw [black, ultra thick,fill=none,-<-] (1.1,-10) node {} -- (1.1,-8) node {};
 \draw [black, ultra thick,fill=none,-<-] (0,-10) node{} -- (0,-8) node {};
 \draw [black, ultra thick,fill=none,-<-] (-1.1,-10) node {} -- (-1.1,-8) node {};
 \draw [black, ultra thick,fill=none,-<-] (2.9,-10) node {} -- (2.9,-8) node {};
 \draw [black, ultra thick,fill=none,-<-] (4,-10) node{} -- (4,-8) node{};
 \draw [black, ultra thick,fill=none,-<-] (5.1,-10) node [right, black] {\small $\mathcal C_p^2\times \mathcal C_{p_2-1}^2 \times \mathcal C_{p_1}$} -- (5.1,-8) node {};

\node[green] at (0,10) {};
\node[purple] at (0,8) {};
\node[red] at (0,6) {};
\node[purple] at (-4,4) {};
\node[purple] at (0,4) {};
\node[purple] at (4,4) {};
\node[green] at (-4,2) {};
\node[green] at (0,2) {};
\node[green] at (4,2) {};
\node[purple] at (-4,0) {};
\node[purple] at (0,0) {};
\node[purple] at (4,0) {};
\node[red] at (-4,-2) {};
\node[red] at (0,-2) {};
\node[red] at (4,-2) {};
\node[purple] at (2.9,-4) {};
\node[purple] at (4,-4) {};
\node[purple] at (5.1,-4) {};
\node[purple] at (1.1,-4) {};
\node[purple] at (0,-4) {};
\node[purple] at (-1.1,-4) {};
\node[purple] at (-2.9,-4) {};
\node[purple] at (-4,-4) {};
\node[purple] at (-5.1,-4) {};
\node[green] at (2.9,-6) {};
\node[green] at (4,-6) {};
\node[green] at (5.1,-6) {};
\node[green] at (1.1,-6) {};
\node[green] at (0,-6) {};
\node[green] at (-1.1,-6) {};
\node[green] at (-2.9,-6) {};
\node[green] at (-4,-6) {};
\node[green] at (-5.1,-6) {};
\node[purple] at (2.9,-8) {};
\node[purple] at (4,-8) {};
\node[purple] at (5.1,-8) {};
\node[purple] at (1.1,-8) {};
\node[purple] at (0,-8) {};
\node[purple] at (-1.1,-8) {};
\node[purple] at (-2.9,-8) {};
\node[purple] at (-4,-8) {};
\node[purple] at (-5.1,-8) {};
\node[red] at (2.9,-10) {};
\node[red] at (4,-10) {};
\node[red] at (5.1,-10) {};
\node[red] at (1.1,-10) {};
\node[red] at (0,-10) {};
\node[red] at (-1.1,-10) {};
\node[red] at (-2.9,-10) {};
\node[red] at (-4,-10) {};
\node[red] at (-5.1,-10) {};

\end{tikzpicture}
\caption{The complex of groups $\mathcal H^p_2$ with $p_1-p=3$ (vertices in any horizontal line have the same local groups, we do not label them all due to space constraints)} \label{pooo}
\end{figure}

\begin{Corollary}\label{bldg} Let $(W,S)$ be a right-angled Coxeter system with free generators $\{s_1,...,s_j\}$ for $j \geq 2$. Let $\Sigma$ be the regular right-angled building of type $(W,S)$ as in Theorem~\ref{buildingthm} and let $ p < \Max_{1 \leq i \leq j}\{p_i\} $ be prime. For any $\alpha \in \N$ there exists a uniform lattice $\Gamma$ in $\Aut(\Sigma)$ with covolume $\mu\big(\Gamma \backslash \Aut(\Sigma) \big) = a/b$ (in lowest terms) such that $b$ is divisible by $p^{\alpha}$.

\begin{proof} Without loss of generality let $p_1=\Max_{1 \leq i \leq j}\{p_i\}$. We call a set of edges $1-1$ if they have no pairwise shared vertices and they form a bijection between two given sets of vertices.

\vspace{2mm}\noindent Take any $k \in \N$. We first construct a complex of groups $\mathcal H^p_k$ with underlying scwol a tree with blue, red and purple vertices. Figure~\ref{pooo} depicts the case where $p_1-p=3$ and $k=2$. No vertex is adjacent to another of the same colour. Edges always point away from purple vertices. Partition $\mathcal H^p_k$ into descending "levels" from $i=0$ to $i=k+1$. 

\vspace{2mm}\noindent Level $0$ consists of a single green vertex with local group $\mathcal C_{p_2}$ and a single edge from the lone purple vertex with trivial vertex group in level $1$ to the lone green vertex in level $0$. Level $k+1$ consists of $1-1$ edges from $(p_1-p)^{k}$ purple vertices each with local group $\mathcal C^{k}_{p} \times \mathcal C^{k}_{p_2-1}$ to the same number of red vertices each with local group $\mathcal C^{k}_{p} \times \mathcal C^{k}_{p_2-1} \times \mathcal C_{p_1}$.

\vspace{2mm}\noindent Every remaining level $i$ for $i\in \{1,...,k\}$ consists of: \begin{itemize} 
\item $(p_1-p)^{i-1}$ purple vertices each with local group $\mathcal C^{i-1}_{p} \times \mathcal C^{i-1}_{p_2-1}$, which we say are of \textit{type} I;
\item $(p_1-p)^{i}$ purple vertices each with local group $\mathcal C^{i}_{p} \times \mathcal C^{i-1}_{p_2-1}$, which we say are of \textit{type} II;
\item $(p_1-p)^{i-1}$ red vertices each with local group $\mathcal C^{i}_{p} \times \mathcal C^{i-1}_{p_2-1}$;
\item $(p_1-p)^{i}$ green vertices each with local group $\mathcal C^{i}_{p} \times \mathcal C^{i}_{p_2-1}$;
\item $1-1$ edges from each type I purple vertex to each red vertex;
\item $1-1$ edges from each type II purple vertex to each green vertex;
\item edges to each red vertex from (disjoint sets of) $p_1-p$ type II purple vertices; and
\item $1-1$ edges from each type I purple vertex in level $i+1$ to each green vertex in level $i$.
\end{itemize}

\noindent Recall the definition of the complex of groups $\widehat{K'}$ from the proof of Theorem~\ref{buildingthm}. We construct $\mathcal A^p_k$ from $\mathcal H^p_k$ by gluing a multiple of $\widehat{K'}$ to each purple vertex, as in Figure~\ref{ggg}, which depicts the case where $p_1-p=1$, $p_2=2$ and $k$ is arbitrary.

\vspace{2mm}\noindent The fundamental group of each complex of groups in $\{\mathcal A^p_k\}_{k \in \N}$ is a uniform lattice in $\Sigma$ by the same argument as in the proof of Theorem~\ref{buildingthm}. A simple calculation will show that the corresponding sequence of covolumes converges at a rate of $\frac{1}{p(p_2-1)}$ as $k \to \infty$. The sum of reciprocals of orders of the vertex groups (in lowest terms) will therefore have a denominator that is divisible by ever-increasing powers of $p$ as $k$ gets larger.
\end{proof}
\end{Corollary}


\begin{thebibliography}{30}

\bibitem{B}
 {\sc H.~Bass}, 
   {\em Covering theory for graphs of groups}, 
   Journal of Pure and Applied Algebra 89, no. 1 (1993), 3--47.

\bibitem{BL}
 {\sc H.~Bass and \sc A.~Lubotzky}, 
   {\em Tree lattices}, 
   Birkhauser Boston (2001).

\bibitem{BH}
 {\sc M.R.~Bridson and \sc A.~Haefliger}, 
   {\em Metric spaces of non-positive curvature}, 
   Springer-Verlag, Berlin (1999).

\bibitem{BM}
 {\sc M.~Burger and \sc S.~Mozes}, 
   {\em Lattices in product of trees}, 
   Publications Math{\'e}matiques de l'IH{\'E}S 92, no. 1 (2000), 151-194.

\bibitem{D}
 {\sc M.~Davis}, 
   {\em The geometry and topology of Coxeter groups}, 
   Vol. 32, Princeton University Press (2008).

\bibitem{F}
 {\sc B.~Farb, \sc G.C.~Hruska and \sc A.~Thomas}, 
   {\em Problems on automorphism groups of nonpositively curved polyhedral complexes and their lattices}, 
   in Geometry, Rigidity, and Group Actions, B. Farb and D. Fisher (eds), the University of Chicago Press, Chicago (2011).

\bibitem{FH}
 {\sc B.~Farb and \sc G.C.~Hruska}, 
   {\em Commensurability invariants for nonuniform tree lattices}, 
  Israel Journal of Mathematics 152, no. 1 (2006), 125-142..

\bibitem{HP}
 {\sc F.~Haglund and \sc F.~Paulin}, 
   {\em Constructions arborescentes d'immeubles}, 
   Annals of Mathematics 325 (2003), 137-164.
 
\bibitem{L}
 {\sc A.~Lubotzky}, 
   {\em Tree lattices and lattices in Lie groups}, 
   Combinatorial and geometric group theory, Edinburgh (1993), 217--232.

\bibitem{La}
 {\sc N.~Lazarovich}, 
   {\em Uniqueness of homogeneous $\CAT(0)$ polygonal complexes}, 
   Geometriae dedicata 168, no. 1 (2014), 397-414.

\bibitem{M}
 {\sc G.~ Moussong}, 
   {\em Hyperbolic Coxeter groups}, 
  PhD Dissertation, Ohio State University (1987).

\bibitem{R}
 {\sc G.~ Rosenberg}, 
   {\em Towers and covolumes of tree lattices}, 
  PhD Dissertation, Columbia University (2001).

\bibitem{S}
 {\sc J.P.~Serre}, 
   {\em Trees}, 
   Springer Monographs in Mathematics, Springer-Verlag, Berlin (2003).

\bibitem{S1}
 {\sc J.P.~Serre}, 
   {\em Cohomologie des groupes discrets}, 
   Annals of Mathematics Vol. 70, Princeton University Press, Princeton (1971), 77--169.

\bibitem{Si}
 {\sc C.L.~Siegel}, 
   {\em Some remarks on discontinuous groups}, 
   Annals of Mathematics Vol. 46,  no. 4, Princeton University Press, Princeton (1945), 708-718.

 \bibitem{T1}  
 {\sc A.~Thomas}, 
   {\em Covolumes of uniform lattices acting on polyhedral complexes}, 
   Bulletin of the London Mathematical Society 39(\textbf{1})  (2007), 103--111.
 
 \bibitem{T2}
  {\sc A.~Thomas},  
   {\em Existence, covolumes and infinite generation of lattices for Davis complexes}, 
  Groups, Geometry and Dynamics 6(\textbf{4}) (2012), 765--801.

\bibitem{T3}
  {\sc A.~Thomas},  
   {\em Lattices acting on right-angled buildings}, 
  Algebraic and Geometric Topology 6, no. 3 (2006), 1215-1238.

\bibitem{T4}
  {\sc A.~Thomas and A. Kubena},  
   {\em Density of commensurators for uniform lattices of right-angled buildings}, 
  arXiv preprint arXiv: 0812.2280 (2008).

\bibitem{W1}
 {\sc G.~White}, 
   {\em Automorphisms of geometric structures associated to Coxeter groups}, 
    preprint, arXiv:1202.6441 (2012).

\end{thebibliography}
\end{document}